\documentclass[10Ishpt,a4paper,twoside]{amsart}
\usepackage[usenames, dvipsnames]{xcolor}
\usepackage{bm}
\usepackage{amsmath, latexsym, amsfonts, amssymb, amsthm}
\usepackage{graphicx,dsfont,epsfig,caption,wrapfig,subfig}
\usepackage{pstricks}
\usepackage[pagebackref]{hyperref}
\newcommand{\mapolicebackref}[1]{
    \hspace*{\fill} \mbox{\textit {\small #1}}
    }
\renewcommand*{\backref}[1]{}
\renewcommand*{\backrefalt}[4]{%
\ifcase #1 \mapolicebackref{pas de citations}
    \or \mapolicebackref{(cited at page #2)}
    \else \mapolicebackref{(cited at pages #2)}
\fi
}
\hypersetup{
    linktoc=page,
    linkcolor=red,          
    citecolor=blue,        
    filecolor=blue,      
     urlcolor=blue,
    colorlinks=true           
}

\usepackage{datetime}
\usepackage{tikz}
 \usetikzlibrary[patterns]
\usepackage{movie15}
\hypersetup{
    linktoc=page,
    linkcolor=red,          
    citecolor=blue,        
    filecolor=blue,      
     urlcolor=blue,
    colorlinks=true           
}

\numberwithin{equation}{section}

\newtheorem{theorem}{Theorem}[section]
\newtheorem{lemma}[theorem]{Lemma}
\newtheorem{proposition}[theorem]{Proposition}
\newtheorem{corollary}[theorem]{Corollary}

\newtheorem{remark}[theorem]{Remark}

\theoremstyle{definition}

\definecolor{remi}{rgb}{0,0,0}

\renewcommand{\tilde}{\widetilde}          
\DeclareMathSymbol{\leqslant}{\mathalpha}{AMSa}{"36} 
\DeclareMathSymbol{\geqslant}{\mathalpha}{AMSa}{"3E} 
\DeclareMathSymbol{\eset}{\mathalpha}{AMSb}{"3F}     
\renewcommand{\leq}{\;\leqslant\;}                   
\renewcommand{\geq}{\;\geqslant\;}                   

\newcommand{\R}{\mathbb{R}}
\newcommand{\Z}{\mathbb{Z}}
\newcommand{\N}{\mathbb{N}}
\newcommand{\Q}{\mathbb{Q}}

\newcommand{\E}{\mathds{E}}
\renewcommand{\P}{\mathds{P}}
\newcommand{\ind}{\mathds{1}}

\newcommand{\caA}{{\mathcal A}}
\newcommand{\caB}{{\mathcal B}}
\newcommand{\caC}{{\mathcal C}}

\newcommand{\caE}{{\mathcal E}}

\newcommand{\caH}{{\mathcal H}}

\newcommand{\caK}{{\mathcal K}}

\newcommand{\caO}{{\mathcal O}}

\newcommand{\caT}{{\mathcal T}}

\newcommand{\C}{\mathbb{C}}

\renewcommand{\P}{\mathds{P}}

\newcommand{\hf}{\frac{_1}{^2}}

\usepackage{xparse}

\def\bi{\begin{itemize}}
\def\ei{\end{itemize}}

\def\bnum{\begin{enumerate}}
\def\enum{\end{enumerate}}
\def\<#1{\langle #1 \rangle}

\title{The DOZZ formula}

\begin{document}

\title[The DOZZ formula]{Integrability of Liouville theory: proof of the DOZZ Formula}

\author[Antti Kupiainen]{Antti Kupiainen$^{1}$}
\address{University of Helsinki, Department of Mathematics and Statistics,
         P.O. Box 68 , FIN-00014 University of Helsinki, Finland}
\email{antti.kupiainen@helsinki.fi}

\author[R\'emi Rhodes]{R\'emi Rhodes$^{2}$}
\address{ Aix-Marseille Universit{\'e}, CNRS, Centrale Marseille,  I2M UMR 7373,  Technop\^ole Ch\^ateau-Gombert,
39 rue F. Joliot Curie, 13453 Marseille, France}
\email{remi.rhodes@univ-amu.fr}

\author[Vincent Vargas]{ Vincent Vargas$^{2}$}
\address{ENS Ulm, DMA, 45 rue d'Ulm,  75005 Paris, France}
\email{Vincent.Vargas@ens.fr}

\footnotetext[1]{Supported by the Academy of Finland and ERC Advanced Grant 741487,$^2$Research supported in part by ANR grant Liouville (ANR-15-CE40-0013)}



 \begin{abstract}
Dorn and Otto (1994) and independently Zamolodchikov and Zamolodchikov (1996) proposed a remarkable explicit expression, the so-called DOZZ formula,  for the 3 point structure constants of Liouville Conformal Field Theory (LCFT), which is expected to describe the scaling limit of large planar maps properly embedded into the Riemann sphere. In this paper we give a proof of the DOZZ formula based on a rigorous probabilistic construction of LCFT in terms of Gaussian Multiplicative Chaos given earlier by F. David and the authors. This result is a fundamental step in the path to prove integrability of LCFT, i.e. to mathematically justify the methods of Conformal Bootstrap used by physicists. From the purely probabilistic point of view, our proof constitutes the first nontrivial rigorous integrability result on Gaussian Multiplicative Chaos measures.  
 \end{abstract}

\maketitle
\tableofcontents


\begin{center}
\end{center}
\footnotesize



\normalsize

\section{Introduction}

A. Polyakov introduced Liouville Conformal Field theory (LCFT hereafter) in his 1981 seminal paper \cite{Pol} 
where he proposed a (non-rigorous) way to put a measure on the set of  Riemannian metrics over a fixed two dimensional manifold; in this context, an integral with respect to the measure is called a functional integral. 
Ever since, the work of Polyakov has echoed in various branches of physics and mathematics, ranging from string theory to probability theory through geometry. In the context of $2D$ quantum gravity models, Polyakov's approach  is conjecturally equivalent to the scaling limit of Random Planar Maps (RPM for short), which are natural probability measures over finite size triangulations of a fixed Riemann surface (see \cite{JFL} for an introduction and further references). In the case of uniform RPM, the proof of this equivalence has culminated in the series of works \cite{MS1,MS2,MS3} by Miller and Sheffield (one may also consult \cite[Appendix 5.3]{DKRV} for a statement of the general conjecture).

Motivated by an attempt to ``solve" LCFT 
Belavin, Polyakov and Zamolodchikov (BPZ hereafter) formulated in their 1984 paper
\cite{BPZ} the general structure of Conformal Field Theory (CFT hereafter). In the BPZ approach the basic objects of CFT are correlation functions of random fields and solving CFT  consists in deriving explicit expressions for them. BPZ proposed to construct  the correlation functions of a CFT  recursively from two inputs:  the {\it  spectrum} and the {\it three point structure constants}. Although we will not define the spectrum  in this paper, let us just note that the spectrum encodes the algebraic structure of the CFT that allows one to determine higher order correlation functions knowing the three point correlation functions (see Section \ref{sub:bootstrap}).  This recursive procedure to find higher point correlation functions is called {\it Conformal Bootstrap}.
Though BPZ were able to find the spectrum and structure constants for a large class of CFT's  (e.g. the Ising model)  LCFT was not one of them\footnote{Following their work \cite{BPZ}, Polyakov qualified CFT as an ``unsuccessful attempt to solve the Liouville model" and did not at first want to publish his work, see \cite{Pol1}.}.  The spectrum of LCFT was soon conjectured in \cite{ct, bct, gn} but the structure constants remained a puzzle.

A decade later, Dorn and Otto \cite{DoOt} and independently Zamolodchikov and Zamolodchikov \cite{ZaZaarxiv, ZaZa} (DOZZ hereafter)  proposed  a remarkable formula for the structure constants of LCFT the so-called DOZZ formula. Even by the physicists' standards the derivation was lacking rigor. To quote Zamolodchikov and Zamolodchikov \cite{ZaZaarxiv}: 
``It should be stressed that the arguments of this section have nothing to do with a derivation. These are rather some motivations and we consider the expression proposed as a guess which we try to support in the subsequent sections." Ever since these papers the derivation of the DOZZ formula from the original  (heuristic)  functional integral definition of LCFT given by Polyakov has remained a controversial open problem, even on the physical level of rigor. The subsequent derivations of the DOZZ formula in the physics literature were based on general assumptions of CFT combined with assumptions on the spectrum of LCFT \cite{Tesc,Tesc1,Rib} or attempts to show its integrability \cite{Tesc2}. The approach of Teschner \cite{Tesc} in fact plays an important role in our proof. 

Recently the present authors together with F. David gave a rigorous probabilistic construction of Polyakov's LCFT functional integral \cite{DKRV} and the correlation functions that are the basic objects in the BPZ approach.  Subsequently in \cite{KRV} we proved identities for the correlation functions postulated in the work of BPZ (conformal Ward identities and BPZ equations). 

 The approach in  \cite{DKRV} is based on the   probabilistic theory of Gaussian Multiplicative Chaos (GMC) which enables to define random measures by exponentiating the 2d Gaussian Free Field (GFF). The terminology of  Gaussian Multiplicative Chaos goes back to Kahane in the eighties \cite{cf:Kah} and is concerned with measures constructed by exponentiating log-correlated fields. This theory is neither restricted to the framework of conformal invariance nor to dimension 2. It enables to define (random) measures formally given by $e^{\gamma X(x)} \sigma(dx)$ where $\gamma$ is a parameter, $X$ a log-correlated field and $\sigma$ a Radon measure on some subset of $\R^d$ with $d \geq 1$
. Therefore it generalizes previous works by Albeverio, Gallavotti and Hoegh-Krohn \cite{Albeverio1,Albeverio2} who initiated the study of two-dimensional exponential interactions in the context of Constructive Field Theory and works  on  multiplicative cascades in the realm of fully developed turbulence (see \cite{review} for references). 

The works 
\cite{DKRV} and  \cite{KRV} 
 provide a probabilistic setup to address the issues of conformal bootstrap and in particular the DOZZ formula.
In this paper we address the second problem: we prove that the probabilistic expression given in  \cite{DKRV} for the structure constants is indeed given by the DOZZ formula. Our result should be considered as an integrability result for LCFT and in particular for the specific GMC measure defined in two dimensions  by exponentiating the GFF on the Riemann sphere.   As such it constitutes the first nontrivial rigorous integrability result in GMC theory. Let us mention as supplementary materials  the manuscripts \cite{KRVjhep,varg} which summarize the content of this paper.

Many integrability formulas for GMC theory (in the one dimensional context) have been conjectured in statistical physics in the study  of  disordered systems. In particular an explicit formula for the moments of the  total mass of the GMC  measure on the circle (based on exponentiating the free boundary GFF) was proposed by Fyodorov-Bouchaud  \cite{fybu} (for generalizations to  other 1d geometries like the segment see the work by Fyodorov-Le Doussal-Rosso \cite{FLeR} and Ostrovsky \cite{ostrovsky2, ostrovsky}). It turns out that their formula is a particular case of the conjectured  one point  bulk  structure constant for LCFT on the unit disk with boundary (these formulas can be found in Nakayama's review  \cite{nakayama}). The recent work of Remy \cite{remy} demonstrates that the approach in this paper  can be adapted to the case of the disk to give a proof  of the Fyodorov-Bouchaud formula. These methods were further extended by Remy-Zhu  to the case of an interval in \cite{zhu}. More generally, we believe the methods developed in this paper and the previous companion paper \cite{KRV} will lead to numerous new integrability results in the field of GMC.    

It should be noted that the LCFT structure constants and the DOZZ formula have a wide range of applications in CFT. Indeed, it has been argued  \cite{Rib} that   LCFT seems to be a universal CFT: e.g. the minimal model structure constants (e.g. the Ising model, tri-critical Ising model and the 3 states Potts model) originally found by BPZ may be recovered from the DOZZ formula by analytic continuation. Furthermore there is strong numerical evidence \cite{bootstrap} that LCFT is essentially the unique CFT for central charge $c>1$: the conformal bootstrap equations seem to have the DOZZ structure constants as their only solution. In another spectacular development, the LCFT structure constants show up in a seemingly completely different setup of four dimensional gauge theories via the so-called AGT correspondence  \cite{AGT} (see the work by Maulik-Okounkov \cite{MO} and Schiffmann-Vasserot \cite{SV} for the mathematical implications in quantum cohomology  of these ideas). 

In the remaining part of this introduction, we briefly review the functional integral approach to LCFT and state the DOZZ formula. 

\subsection{LCFT correlation functions }\label{pathi}

A rigorous formulation of LCFT will be given later  (see Section \ref{sec:backgr}). Heuristically Polyakov's formulation of LCFT on the Riemann sphere $\hat \C$ is the study of  conformal metrics on $\hat \C$  of the form $e^{\gamma\phi(z)}|dz|^2$ where $z$ is the standard complex coordinate; in this context, let $d^2z$ denote the  Lebesgue measure. $\phi(z)$ is a random function (in fact, $\phi(z)$ turns out to be a random distribution in the sense of Schwartz once properly defined mathematically). Expectations of suitable functions of $\phi$ are given by the formal integral\footnote{ Restricting $F$ to indicator functions indeed gives rise to a measure on some appropriate functional space. We use brackets and not $\E$ for the positive linear functional  \eqref{liouvlawintro} since it turns out that the measure  $e^{-S_L(\phi)}D\phi$ once rigorously defined is not normalizable into a probability measure.}
\begin{equation}\label{liouvlawintro}
\langle F\rangle :=\int F(\phi)e^{-S_L(\phi)}D\phi 
 \end{equation}   
where $S_L$ is the {\it Liouville Action functional}
\begin{equation}\label{liouvlawintro1}
S_L(\phi)=\frac{_1}{^{\pi}}\int_\C(|\partial_z\phi(z)|^2+\pi\mu e^{\gamma\phi(z)})d^2z.
 \end{equation} 
The formal ``functional integral'' \eqref{liouvlawintro}, once rigorously defined,  gives rise to  a Conformal Field Theory, the LCFT, which is the topic of this paper.

LCFT has two parameters $\gamma \in (0,2)$ and $\mu>0$.The parameter $\mu$ 
is called the cosmological constant and for LCFT $\mu$ has to be strictly positive. The case $\mu=0$ corresponds to the Free Field theory, which is a different  Conformal Field Theory with different structure constants. The precise value of $\mu$ in LCFT plays no specific role since the dependence on $\mu$ 
is governed by a scaling relation, see \cite{DKRV}.  On the other hand, the parameter $\gamma$ encodes the conformal structure of the theory; more specifically, one can show that the central charge\footnote{In this article this concept will not appear and hence we  
refer to the works \cite{DKRV}, \cite{KRV} for an account on the central charge.}  of the theory is $c_L=1+6 Q^2$ with 
 \begin{equation}\label{valueQ}
 Q= \frac{2}{\gamma}+ \frac{\gamma}{2}.
 \end{equation} 
 
 The basic objects of interest in LCFT are in physics terminology {\it vertex operators}
  \begin{equation}\label{KPZformula1}
V_{\alpha}(z)  = e^{\alpha\phi(z)}
\end{equation} 
where $\alpha$ is a complex number and their correlation functions 
 $ \langle \prod_{k=1}^N V_{\alpha_k}(z_k)     \rangle
$. 
Their definition involves a regularization and renormalization procedure and they were constructed rigorously in \cite{DKRV} for $N \geq 3$ and for real $\alpha_i$ satisfying certain conditions. The construction of the correlations in \cite{DKRV} is probabilistic and based on interpreting $e^{-\frac{1}{^{\pi}}\int_\C |\partial_z\phi(z)|^2d^2z} D\phi$ in terms of a suitable Gaussian Free Field (GFF) probability measure: see subsection \ref{sec:backgr} below for precise definitions and an explicit formula for the correlations 
 in terms of   the GMC associated to the GFF. 
 
 In particular it was proved  in \cite{DKRV} that these correlation functions 
are {\it conformal tensors}. More precisely, if $z_1, \dots, z_N$ are $N$ distinct points in $\C$  then for a M\"obius map $\psi(z)= \frac{az+b}{cz+d}$ (with $a,b,c,d \in \C$ and $ad-bc=1$) 
 \begin{equation}\label{KPZformula}
\langle \prod_{k=1}^N V_{\alpha_k}(\psi(z_k))    \rangle=  \prod_{k=1}^N |\psi'(z_k)|^{-2 \Delta_{\alpha_k}}     \langle \prod_{k=1}^N V_{\alpha_k}(z_k)     \rangle
\end{equation}  
where  $\Delta_{\alpha}=\frac{\alpha}{2}(Q-\frac{\alpha}{2})$ is called the conformal weight.
This global conformal symmetry fixes the three point correlation functions up to a constant: 
\begin{equation}\label{firstdefstructure}
 \langle      \prod_{k=1}^3 V_{\alpha_k}(z_k)  \rangle 
 =  |z_1-z_2|^{ 2 \Delta_{12}}  |z_2-z_3|^{ 2 \Delta_{23}} |z_1-z_3|^{ 2 \Delta_{13}}C_\gamma(\alpha_1,\alpha_2,\alpha_3)
\end{equation}
with $\Delta_{12}= \Delta_{\alpha_3}-\Delta_{\alpha_1}-\Delta_{\alpha_2}$, etc. The constants $C_\gamma(\alpha_1, \alpha_2, \alpha_3)$ are called the three point structure constants and they have an explicit expression in terms of the GMC associated to the GFF, see Section \ref{Structure constants and four point functions}. They are also the building blocks of LCFT in the conformal bootstrap approach, see subsection \ref{sub:bootstrap}. We should also note that the law of the area measure $e^{\gamma \phi(z)}d^{2}z$ normalized to unit total area when $\phi$  is sampled from   the normalized measure $ \prod_{k=1}^3 V_{\gamma}(z_k) e^{-S_L(\phi)}D\phi$ coincides with the law of the (unit area) three point quantum sphere as defined by Duplantier, Miller and Sheffield \cite{DMS} (recall that quantum spheres are equivalence classes of random measures on the sphere with two marked points $0$ and $\infty$; one can construct a three point quantum sphere by sampling a point $z$ according to the quantum sphere and  taking the image of the quantum sphere by the unique M\"obius maps that sends the points $0,z,\infty$ to the fixed points $z_1, z_2,z_3$): see \cite{aru} for the equivalence.

\subsection{The DOZZ formula}

As mentioned above, an explicit expression for the LCFT structure constants was proposed  in \cite{DoOt,ZaZa} . Subsequently it was observed by Teschner \cite{Tesc} that this formula may be derived by applying the bootstrap framework to  special four point functions (see section \ref{bpzeq}). He argued that this leads to the following remarkable periodicity relations for the structure constants:
\begin{align}\label{3pointconstanteqintro}
C_\gamma(\alpha_1+\frac{_\gamma}{^2},\alpha_2,\alpha_3)&=- \frac{_1}{^{\pi \mu}} \mathcal{A}_\gamma(\frac{_\gamma}{^2},\alpha_1,\alpha_2,\alpha_3)C_\gamma(\alpha_1-\frac{_\gamma}{^2},\alpha_2,\alpha_3)\\
C_\gamma(\alpha_1+\frac{_2}{^\gamma},\alpha_2,\alpha_3)&=- \frac{_1}{^{\pi \tilde{\mu}}} \mathcal{A}_\gamma(\frac{_2}{^\gamma},\alpha_1,\alpha_2,\alpha_3)C_\gamma(\alpha_1-\frac{_2}{^\gamma},\alpha_2,\alpha_3) \label{3pointconstanteqintrodual}
\end{align}
with $\tilde{\mu}= \frac{(\mu \pi l(\frac{\gamma^2}{4})  )^{\frac{4}{\gamma^2} }}{ \pi l(\frac{4}{\gamma^2})}$ and 
\begin{align}\label{Aformula}
\mathcal{A}_\gamma(\chi,\alpha_1,\alpha_2,\alpha_3)=\frac{l(-\chi^2)  l(\chi\alpha_1) l(\chi\alpha_1-\chi^2 )  l(\frac{\chi}{2} (\bar{\alpha}-2\alpha_1- \chi) )   }{l( \frac{\chi}{2} (\bar{\alpha}-\chi - 2Q)  ) l( \frac{\chi}{2} (\bar{\alpha}-2\alpha_3-\chi ))  l( \frac{\chi}{2} (\bar{\alpha}-2\alpha_2-\chi )) }
\end{align}
where $\bar{\alpha}= \alpha_1+\alpha_2+\alpha_3$ and
\begin{equation}\label{Defl}
l(x)=\Gamma(x)/\Gamma(1-x).
\end{equation} 

The equations \eqref{3pointconstanteqintro}, \eqref{3pointconstanteqintrodual} have a meromorphic solution which is the DOZZ formula. It is expressed in terms of a special function $\Upsilon_{\frac{\gamma}{2}}(z)$ defined for $0<\Re (z)< Q$ 
by the formula\footnote{The function 
has a simple construction in terms of standard double gamma functions: see the reviews \cite{nakayama,Rib,Tesc1} for instance.}
\begin{equation}\label{def:upsilon}
\ln \Upsilon_{\frac{\gamma}{2}} (z)  = \int_{0}^\infty  \big ( (\tfrac{Q}{2}-z  )^2  e^{-t}-  \frac{( \sinh( (\frac{Q}{2}-z )\frac{t}{2}  )   )^2}{\sinh (\frac{t \gamma}{4}) \sinh( \frac{t}{\gamma} )}    \big ) \frac{dt}{t}.
\end{equation}
The function $\Upsilon_{\frac{\gamma}{2}}$  can be analytically continued to $\C$ because it satisfies remarkable functional relations: see formula \eqref{shiftUpsilon} in the appendix. It has no poles in $\C$ and the zeros of $\Upsilon_{\frac{\gamma}{2}}$ are simple (if $\gamma^2 \not \in \Q$) and given by the discrete set $(-\frac{\gamma}{2} \N-\frac{2}{\gamma} \N) \cup (Q+\frac{\gamma}{2} \N+\frac{2}{\gamma} \N )$. With these notations, the DOZZ formula (or proposal) $C_\gamma^{{\rm DOZZ}}(\alpha_1,\alpha_2,\alpha_3)$ is the following expression
\begin{equation}\label{DOZZformula}
C_\gamma^{{\rm DOZZ}}(\alpha_1,\alpha_2,\alpha_3)= (\pi \: \mu \:  l(\tfrac{\gamma^2}{4})  \: (\tfrac{\gamma}{2})^{2 -\gamma^2/2} )^{\frac{2 Q -\bar{\alpha}}{\gamma}}   \frac{\Upsilon_{\frac{\gamma}{2}}'(0)\Upsilon_{\frac{\gamma}{2}}(\alpha_1) \Upsilon_{\frac{\gamma}{2}}(\alpha_2) \Upsilon_{\frac{\gamma}{2}}(\alpha_3)}{\Upsilon_{\frac{\gamma}{2}}(\frac{\bar{\alpha}-2Q}{2}) 
\Upsilon_{\frac{\gamma}{2}}(\frac{\bar{\alpha}}{2}-\alpha_1) \Upsilon_{\frac{\gamma}{2}}(\frac{\bar{\alpha}}{2}-\alpha_2) \Upsilon_{\frac{\gamma}{2}}(\frac{\bar{\alpha}}{2} -\alpha_3)   } .
\end{equation} 

The main result of the present paper is to show the first important equality between LCFT in the functional integral formulation (rigorously defined in \cite{DKRV} via probability theory) and the conformal bootstrap approach, namely to prove that for $\gamma \in (0,2)$ and appropriate $\alpha_1,\alpha_2,\alpha_3$  the structure constants $C_\gamma(\alpha_1,\alpha_2,\alpha_3)$  in \eqref{firstdefstructure} are equal to $C_\gamma^{{\rm DOZZ}}(\alpha_1,\alpha_2,\alpha_3)$ defined by \eqref{DOZZformula}. 

Our proof is based on  deriving the equations \eqref{3pointconstanteqintro}, \eqref{3pointconstanteqintrodual} for the probabilistically defined $C_\gamma$. An essential role in this derivation is an identification in probabilistic terms of  the {\it reflection coefficient} of LCFT. It has been known for a long time \cite{ZaZa,Tesc} that in  LCFT the following reflection relation should hold in some sense:
\begin{equation}\label{rrel}
V_\alpha=R(\alpha) V_{2Q-\alpha}.
\end{equation} 
Indeed the DOZZ formula is compatible with the following form of \eqref{rrel}  \cite{ZaZa}:
\begin{equation}\label{rrel1}
C_\gamma^{{\rm DOZZ}}(\alpha_1,\alpha_2,\alpha_3)=
R^{{\rm DOZZ}}(\alpha_1)C_\gamma^{{\rm DOZZ}}(2Q-\alpha_1,\alpha_2,\alpha_3) 
\end{equation} 
with 
\begin{equation}\label{defRDOZZ}
R^{{\rm DOZZ}}(\alpha)=- (\pi \: \mu \:  l(\tfrac{\gamma^2}{4})  )^{\frac{2 (Q -\alpha)}{\gamma}} \frac{\Gamma(-\frac{\gamma (Q-\alpha)}{2})}{\Gamma(\frac{\gamma (Q-\alpha)}{2})} \:  \frac{\Gamma(-\frac{2 (Q-\alpha)}{\gamma })}{\Gamma(\frac{2 (Q-\alpha)}{\gamma})}.
\end{equation}
The mystery relation \eqref{rrel} lies in the fact that the probabilistically defined $C_\gamma(\alpha_1,\alpha_2,\alpha_3)$ {\it vanish} if any of the  $\alpha_i\geq Q$ whereas they are nonzero for $\alpha_i<Q$, see Section \ref{lcftcorr}. 

In our proof  $R(\alpha) $ emerges from the analysis of the {\it tail behavior} of a GMC observable.  
We prove that it is also given by the following limit
\begin{equation}\label{defRDOZZlim}
4 R(\alpha)= \underset{\epsilon \to 0}{\lim} \: \epsilon \: C_{\gamma}(\epsilon,\alpha,\alpha)
\end{equation}
i.e. $R(\alpha)$ has an interpretation in terms of a renormalized two-point function. We will show that for those values of $\alpha$ such that $R(\alpha)$ makes sense from the functional integral perspective, i.e. $\alpha\in (\tfrac{\gamma}{2},Q)$, 
$$R(\alpha)=R^{{\rm DOZZ}}(\alpha).$$
It turns out that some material related to the coefficient $R(\alpha)$ already appears in the beautiful work by Duplantier-Miller-Sheffield \cite{DMS}: within this framework, the reflection coefficient $R(\alpha)$ can naturally be interpreted as the partition function of the theory underlying the quantum sphere. We will not elaborate more on this point as no prior knowledge of the work by  Duplantier-Miller-Sheffield is required to understand the sequel (see \cite{aru,tail} for an account of the relation between \cite{DKRV} and \cite{DMS}). More precisely, the required background to understand $R(\alpha)$ will be introduced in subsection \ref{Reflection coefficient} below. 

Finally, let us stress that the DOZZ formula \eqref{DOZZformula}   is invariant under the substitution of parameters
 \begin{equation*}
 \frac{\gamma}{2} \leftrightarrow \frac{2}{\gamma}, \quad  \mu  \leftrightarrow \tilde{\mu}= \frac{(\mu \pi \ell(\frac{\gamma^2}{4})  )^{\frac{4}{\gamma^2} }}{ \pi \ell(\frac{4}{\gamma^2})}.
\end{equation*}
This  duality symmetry is at the core of the DOZZ controversy. Indeed this symmetry is not manifest in the Liouville action  
 functional \eqref{liouvlawintro1} though duality was axiomatically assumed by Teschner \cite{Tesc1} in his argument, especially to get \eqref{3pointconstanteqintrodual}. It was subsequently argued that this duality could come from the presence in the action \eqref{liouvlawintro1} of an additional ``dual" potential of the form $e^{\tfrac{2}{\gamma}\phi}$ with cosmological constant $\tilde{\mu}$ in front of it. As observed by Teschner \cite{Tesc1}, this dual cosmological constant may take  negative (even infinite) values, which makes clearly no sense from the functional integral perspective. That is why the derivation of the DOZZ formula from the LCFT functional integral \eqref{liouvlawintro} has remained shrouded in mystery for so long\footnote{Indeed, there are numerous reviews and papers within the physics literature on the functional integral approach of LCFT and its relation with the bootstrap approach but they offer different perspectives and conclusions (see \cite{HaMaWi,OPS,seiberg} for instance).}.
 

\subsection{Organization of the paper}
In the next section, we introduce the probabilistic expressions of the LCFT correlation functions and structure constants and state the main result of the paper: Theorem \ref{theoremDOZZ}. We also discuss briefly the conformal bootstrap conjecture and prospects for a probabilistic approach to proving it. In Section \ref{Reflection coefficient} we introduce the probabilistic definition of the reflection coefficient which is the central concept in our proof and state the main Theorem   \ref{Rtheor} on it. Section \ref{sec:further} gathers some further results from \cite{KRV} on differential equations (the BPZ equations) satisfied by certain four point functions and their consequences. Since the proof of our main result is a combination of several sub-results with interdependencies we present  in Section \ref{sec:bone} the outline of the argument together with a chart of the logical structure. The rest of the paper is devoted to  the proof of the main results, Theorem \ref{theoremDOZZ} and Theorem \ref{Rtheor}.  In section \ref{sec:analytic}, we show that the correlation functions of vertex operators are analytic functions of their arguments $\alpha_k$. 
Section  \ref{tailestimates} is devoted to the study of tail estimates of GMC and their connection with the reflection coefficient. In Section \ref{sec:proofreflection} we prove a Lemma relating the reflection coefficient to the structure constants. In Section \ref{bpzeq}  we study the asymptotics of four point functions when two of their arguments approach each other (``fusion" rules in the physics jargon). This section is the technical core of the paper and the key input in the probabilistic identification of the reflection coefficient. Finally, in Sections \ref{proofreflection} and \ref{proofDOZZ}, Theorems  \ref{Rtheor} and  \ref{theoremDOZZ} are proved.

\subsubsection*{Acknowledgements} The authors wish to thank Fran\c{c}ois David, Sylvain Ribault and Raoul Santachiara for fruitful discussions on Liouville field theory and the conformal bootstrap approach. The authors would also like to thank the anonymous referees for their careful reading of a prior version of this paper; their numerous comments have certainly improved the paper's readability.

\section{Probabilistic Formulation of LCFT  and the Main Result}\label{sec:main}
In this section, we recall the precise definition of the Liouville correlation functions 
as given in  \cite{DKRV} 
and state the main result on the DOZZ formula.

\subsubsection{Conventions and notations.} In what follows,  $z$,  $x,y$ and   $z_1,\dots, z_N$ all denote complex variables. We use the standard notation for complex derivatives $\partial_x= \frac{1}{2}(\partial_{x_1}- i \partial_{x_2})$ and $\partial_{\bar{x}}= \frac{1}{2} (\partial_{x_1}+ i \partial_{x_2})$ for $x=x_1+ix_2$. The  Lebesgue measure on $\C$ (seen as $\R^2$) is denoted by $d^2x$ . We will also denote $|\cdot|$ the norm in $\C$ of the standard Euclidean (flat) metric and for all $r>0$ we will denote by $B(x,r)$ the Euclidean ball of center $x$ and radius $r$. 

\subsection{Gaussian Free Field and Gaussian Multiplicative Chaos}\label{sec:backgr}
The probabilistic definition of   the  integral \eqref{liouvlawintro} goes by expressing it 
in terms of the Gaussian Free Field (GFF).  The   setup is  the Riemann Sphere  $\hat\C= \C \cup \lbrace \infty \rbrace$ equipped with a conformal metric $g(z)|dz|^2$. 
The correlation functions of LCFT will then depend on the metric but they have simple transformation properties under the change of  $g$, the so-called Weyl anomaly formula.  We refer the reader to \cite{DKRV} for this point and proceed here by just stating a formulation that will be useful for the purposes of this paper. 

We define the GFF $X(z)$ 
as the centered Gaussian random field with covariance (see \cite{Dub0,She07} for background on the GFF)
\begin{equation}\label{hatGformula}
\E [ X(x)X(y)] 
=\ln\frac{1}{|x-y|}+\ln|x|_++\ln|y|_+:=G(x,y)
\end{equation}
where we use the notation $ |z|_+=|z|$ if $|z|\geq 1$ and  $ |z|_+=1$ if $|z|\leq 1$. 
\begin{remark}
In the terminology of  \cite{DKRV}, consider the metric   $g(z)=|z|_+^{-4}$ with   scalar  curvature $R_{g}(z):=-4g^{-1}\partial_z\partial_{\bar z}\ln g(z)=8\pi \nu$ with $\nu$ the uniform probability measure on the equator $|z|=1$. Then $X$ is the GFF with zero average on the equator: $\int Xd\nu
=0$. 
\end{remark}

For LCFT we need to consider the exponential of $X$. Since $X$
 is distribution valued a renormalization procedure is needed.  Define the circle average of $X$ by 
\begin{equation}\label{circleav}X_{r}(z):=\frac{1}{2\pi i}\oint_{|w|=e^{-r}} X(z+w)\frac{dw}{w}
\end{equation}
and consider the measure 
\begin{equation}\label{Meps1}
M_{\gamma,r}(d^2x):=
e^{\gamma X_{r}(x)-\frac{\gamma^2}{2} \E[ X_{r}(x)^2] } 
|x|_+^{-4}d^2x.
\end{equation}
Then, 
for $\gamma\in [0,2)$,   we have the convergence  in probability
\begin{equation}\label{law}
M_{\gamma}=\lim_{r\to \infty}M_{\gamma,r}
 \end{equation}
 and convergence is in the sense of weak convergence of measures. This limiting measure is non trivial and is  GMC associated to the field $X$ with respect to the measure $|x|_+^{-4}d^2x$ (see Berestycki's work \cite{Ber} for an elegant and elementary approach to GMC and references). 
 
\begin{remark}\label{indebm} For later purpose we state a useful property of the circle averages. First, $X_0(0)=0$, the processes $r\in\R_+\to X_r(0)$ and $r\in\R_+\to X_{-r}(0)$ are two independent Brownian motions starting from $0$. For $z$ center of a unit ball contained in $B(0,1)^c$  the process $r\in\R_+\to X_r(z)-X_0(z)$ is  also a Brownian motion starting at $0$ and  for distinct points $(z_k)_{1 \leq k \leq N}$ such that the balls $B(z_k,1)\subset B(0,1)^c$ are disjoint the processes $r\mapsto X_r(z_i)-X_0(z_i)$ are mutually independent and independent of the sigma algebra $\sigma\{X(z);z\in [\cup_{k=1}^NB(z_k,1)]^c\}$. This results from a simple check of covariances.
\end{remark}


\subsection{Liouville correlation functions}\label{lcftcorr} We may now give the probabilistic definition of the  integral \eqref{liouvlawintro}\footnote{The global constant $2$ 
  is included to match with the physics literature normalization which is based on the DOZZ formula \eqref{DOZZformula}.}:
\begin{equation}\label{liouvlaw1}
\langle F\rangle :=2 \int_\R e^{-2Q c}  \E\big[F(X 
-2Q\ln|z|_++c)    e^{- \mu e^{\gamma c} M_{\gamma}(\C)    } \big] \: dc
 \end{equation}
where $\E$ is expectation over the GFF. We refer the reader to \cite{DKRV} (or to \cite{KRV} for a brief summary) for the explanation of the connection between \eqref{liouvlawintro} and \eqref{liouvlaw1}. Briefly, the variable $c$ is essential and stems from the fact that in \eqref{liouvlawintro} we need to integrate over all $\phi$ and not only the GFF $X$ which is constrained by the relation $X_0(0)=0$. The origin of the factor $e^{-2Q c} $  is topological and depends on the fact that we work on the sphere $\hat\C$.  The random variable $M_{\gamma}(\C) $ is almost surely finite because $\E M_{\gamma}(\C) =
\int_\C |z|_+^{-4} d^2z<\infty
$. This implies that $\langle\cdot\rangle$ is not normalizable: $\langle 1\rangle=\infty$. 

 The class of $F$ for which \eqref{liouvlaw1} is defined includes suitable vertex operator correlation functions  once these are properly renormalized. We set for $\alpha\in\R$ and $z\in\C$
\begin{equation}\label{Vdefi}
V_{\alpha, \epsilon}(z)= 
e^{\alpha c} e^{\alpha X_{\epsilon}(z)-\frac{\alpha^2}{2} \E[ X_{\epsilon}(z)^2] }|z|_+^{-4\Delta_\alpha}
\end{equation}
where we recall $\Delta_{\alpha}= \frac{\alpha}{2} (Q-\frac{\alpha}{2})$. The point $z$ will often be referred to as {\it insertion point} (or just insertion) and $\alpha$ as insertion weight (or just weight).
Let $z_i\in\C$, $i=1,\dots,N$ with $z_i \not = z_j $ for all $i\neq j$. It was shown in \cite{DKRV} that the  limit
\begin{equation}\label{lcorre}
 \langle \prod_{k=1}^N  V_{\alpha_k}(z_k) \rangle :=   
  \underset{\epsilon \to 0}{\lim}  \:  \langle \prod_{k=1}^N  V_{\alpha_k,\epsilon}(z_k)  \rangle
\end{equation}
exists, is finite and nonzero if and only if the following {\it Seiberg bounds } originally introduced in \cite{seiberg} hold:
\begin{equation}\label{TheSeibergbounds}
\sum_{k=1}^N\alpha_k>2Q, \quad \quad   \alpha_k<Q, \; \;  \forall k.
\end{equation}
The first condition guarantees that the limit is finite and the second that it is non vanishing. Indeed, if there exists $k$ such that $\alpha_k\geq Q$ then the limit is zero.  Note that these bounds imply that for a nontrivial correlation we need at least {\it three} vertex operators; therefore, we have $N \geq 3$ in the sequel. 
The correlation function \eqref{lcorre}  satisfies the    conformal invariance property \eqref{KPZformula}.

The correlation function can be further simplified by performing the change of variables $y=\mu e^{\gamma c}M_\gamma(\C)$ in the $c$-integral (see  \cite[section 3]{DKRV}):
  \begin{align}\label{Z00}
   \langle \prod_{k=1}^N  V_{\alpha_k}(z_k)  \rangle& =  2 \mu^{-s} 
   \gamma^{-1}\Gamma(s)\lim_{\epsilon\to 0}\E  \left [  \prod_{k=1}^N e^{{\alpha_k} X_{\epsilon}(z_k)- \frac{\alpha_k^2}{2}\E X_{\epsilon}(z_k)^2}   |z_k|_+^{-4\Delta_{\alpha_k}}
   M_\gamma(\C) ^{-s}    \right ]
\end{align}
where 
\begin{equation}\label{defs}
s=\frac{\sum_{k=1}^N\alpha_k-2Q}{\gamma}.\footnote{One should notice that $s$ depends on $N$ and $\alpha_1, \cdots, \alpha_N$ (and also $\gamma$). However, in what follows, we will use the generic notation $s$ in the definition of the correlation functions as it should be clear from the context the value of $s$ that we are considering.}
\end{equation}
Using the Cameron-Martin theorem\footnote{In the sequel, the Cameron-Martin theorem will refer to Theorem \ref{th:Girsanov} or Corollary \ref{coro:Girsanov}.}  (we apply Theorem \ref{th:Girsanov} of the Appendix with the Gaussian variable $\mathcal{X}= \sum_{k=1}^N \alpha_k X_{\epsilon}(z_k)$; see also  \cite[Th 3.4 and 3.5]{DKRV}) we may trade the vertex operators to a shift of $X$ to obtain an expression in terms of the multiplicative chaos:
\begin{equation} \label{Z1}
 \langle    \prod_{k=1}^N V_{\alpha_k}(z_k)  \rangle =2 \mu^{-s} \gamma^{-1}\Gamma(s)
 \prod_{i < j} \frac{1}{|z_i-z_j|^{\alpha_i \alpha_j}}\E \left [  \left (  \int_{\C}  F(x,{\bf z}) M_\gamma(d^2x)  \right )^{-s}  \right ] 
\end{equation}
where 
\begin{equation}\label{coulomb}
F(x,{\bf z})=\prod_{k=1}^N \left ( \frac{ |x|_+}{|x-z_k|}  \right )^{\gamma \alpha_k} .
\end{equation}
Thus, up to explicit factors the Liouville correlations are reduced to the study of the random variable $\int_{\C} F(x, {\bf z}) M_\gamma(d^2x) $. In particular, the Seiberg bounds $\alpha_k<Q$ for all $k$ are the condition of integrability of $F$ against the chaos measure  $M_\gamma$ (see  \cite{DKRV}). Furthermore the expression \eqref{Z1} allows us to extend  the definition of the correlation functions to those values of $s \leq 0$ such that the expectation in \eqref{Z1} makes sense: it was  shown in \cite[Lemma 3.10]{DKRV} that 
\begin{equation}  \label{ThextendedSeibergbounds0}
0 < \E \left [  \left (  \int_{\C}  F(x,{\bf z}) M_\gamma(d^2x)  \right )^{-s}  \right ]  < \infty
\end{equation}
provided the following {\it extended Seiberg's bounds} are satisfied
\begin{equation}\label{ThextendedSeibergbounds}
- s
< \frac{_4}{^{\gamma^2}} \wedge \min_{1 \leq k \leq N}  \frac{_2}{^\gamma}(Q-\alpha_k), \quad \quad \alpha_k<Q, \; \; \forall k
\end{equation}
with $s$ given by \eqref{defs}. The standard $\Gamma$ function has poles on the non-positive integers hence for $s=-n$ with $n$ integer and satisfying   \eqref{ThextendedSeibergbounds}, we simply set the correlations to be equal to infinity.

 Under condition \eqref{ThextendedSeibergbounds}, it is also natural to define the so-called {\it  unit volume correlations}  by 
 \begin{equation}\label{unitvol}
  \langle    \prod_{k=1}^N V_{\alpha_k}(z_k)  \rangle_{\mathrm{uv}}= \mu^s  \frac{ \langle    \prod_{k=1}^N V_{\alpha_k}(z_k)  \rangle}{\Gamma(s)}
 \end{equation}      
i.e. we divide by the $\Gamma$ function to remove the mentioned poles; hence $ \langle    \prod_{k=1}^N V_{\alpha_k}(z_k)  \rangle_{\mathrm{uv}}$ is well defined under condition \eqref{ThextendedSeibergbounds}. An important ingredient in our proof of the DOZZ formula is Theorem \ref{analytic_hyperbolic} which says that   these correlation functions have an analytic continuation in the $\alpha_i$'s to a complex neighborhood of the region allowed by the bounds \eqref{ThextendedSeibergbounds}.
\begin{remark} \label{analingamma} The DOZZ formula for the structure constants is analytic not only in $\alpha_i$ but also in $\gamma$. A direct proof of analyticity of the probabilistic correlation functions in $\gamma$ seems difficult. However, it is an easy exercise in Multiplicative Chaos theory to prove their {\it continuity} in $\gamma$, a fact we will need in our argument. Actually, it is not hard to prove that they are $C^\infty$ in $\gamma$ but we will omit this as it is not needed in our argument.
\end{remark}

\subsection{Structure constants}\label{sub:struct}
The structure constants  $C_\gamma$ in \eqref{firstdefstructure} can be recovered as the following limit
\begin{equation}\label{Climit}
C_\gamma(\alpha_1,\alpha_2,\alpha_3)=\lim_{z_3\to\infty} |z_3|^{4 \Delta_3} \langle       V_{\alpha_1}(0)  V_{\alpha_2}(1) V_{\alpha_3}(z_3) \rangle,
\end{equation}
where here and thereafter we use the shortcut notation $\Delta_j:=\Delta_{\alpha_j}$. 
Combining \eqref{Z1}  with  \eqref{Climit} we get 
\begin{equation}\label{expression3pointstruct}
C_\gamma(\alpha_1,\alpha_2,\alpha_3)= 
2 \mu^{-s}  \gamma^{-1}\Gamma(s)  \E (\rho(\alpha_1,\alpha_2,\alpha_3)^{-s})
\end{equation}
where $s=(\sum_{i=1}^3\alpha_i-2Q)/\gamma$ and
\begin{equation*}
\rho(\alpha_1,\alpha_2,\alpha_3)=  
 \int_{\C}  
  \frac{|x|_+^{\gamma(\alpha_1+\alpha_2+\alpha_3)}}{ |x|^{\gamma \alpha_1}  |x-1|^{\gamma \alpha_2}  }    
  M_\gamma(d^2x).
\end{equation*}
Furthermore, using \eqref{KPZformula} and \eqref{firstdefstructure} we see that $C_{\gamma}$  is a symmetric function of the variables $\alpha_1,\alpha_2,\alpha_3 $.

\subsection{Statement of the main result}


The main result of this paper is the following identity:

\begin{theorem}\label{theoremDOZZ}
Let $\alpha_1,\alpha_2,\alpha_3$ satisfy the bounds \eqref{ThextendedSeibergbounds} with $N=3$. The following equality  holds $$C_\gamma(\alpha_1,\alpha_2,\alpha_3)=C_\gamma^{{\rm DOZZ}}(\alpha_1,\alpha_2,\alpha_3).$$
\end{theorem}

From the purely probabilistic point of view, Theorem \ref{theoremDOZZ} can be interpreted as a far reaching integrability result on GMC on the Riemann sphere; indeed recall that $C_\gamma(\alpha_1,\alpha_2,\alpha_3)$ has an expression in terms of a fractional moment of some form of GMC: see formula \eqref{expression3pointstruct}. There are numerous integrability results on GMC in the physics literature (see the introduction); to the best of our knowledge, Theorem \ref{theoremDOZZ} is the first rigorous nontrivial integrability result on GMC; as argued in the introduction, we believe the techniques of this paper and the companion paper \cite{KRV} will enable to prove many other integrability results for GMC.

\subsection{Further work: conformal bootstrap}\label{sub:bootstrap}

 Theorem \ref{theoremDOZZ} is also an integrability result on LCFT. Based on general principles of conformal field theories as spelled out by Belavin, Polyakov and Zamolodchikov one expects that the correlation function  \eqref{Z1} also has a (semi)  explicit expression in terms of the structure constants. The recursive procedure to obtain this expression is called conformal bootstrap. It postulates a recursion relating an $N$-point correlation function to $N-1$-point correlation functions with coefficients involving the structure constants. Applying bootstrap to the LCFT four point function one obtains the following conjecture \cite{ZaZa}
\begin{multline}\label{eq:cb}
  \langle  V_{\alpha_1} (0) V_{\alpha_2} (z) V_{\alpha_3} (1)  V_{\alpha_4} (\infty) \rangle   \\
=  \int_{\R} |z|^{2(\triangle_{Q+iP}-\triangle_{\alpha_1}-\triangle_{\alpha_2})} C_\gamma(\alpha_1, \alpha_2, Q-iP)  C_\gamma(Q+iP, \alpha_3, \alpha_4  )  |  \mathcal{F}_{P,\{\alpha_i\}} (z) |^2   \, \frac{dP}{8\pi}
\end{multline}
where  $\mathcal{F}_{P,\{\alpha_i\}} (z)$ are  meromorphic functions (the so-called universal conformal blocks) which depend only on the parameters $\alpha_i,  P$ and $\gamma$ through  the central charge of LCFT $c_L=1+6 Q^2$.  The integral over $P$ is here the standard Lebesgue integral over $\R$.  Note that the structure constants in this expression are evaluated at complex weights $Q\pm iP$ and have to be interpreted in terms of analytic continuation from the the real weights. 
 Indeed, our proof constructs this continuation and shows it is given by the DOZZ formula. In the physics terminology these complex weights determine the {\it spectrum} of LCFT. This means in particular that one expects that to LCFT there corresponds a canonical Hilbert space $\caH$ and a unitary representation of the Virasoro algebra with central charge $c_L=1+6 Q^2$ on $\caH$.
This representation is expected to reduce to a direct integral of highest weight representations indexed by $P$. The bootstrap conjecture then formally follows from representation theory.

On the mathematical level \eqref{eq:cb} remains a conjecture. However   Baverez and Wong \cite{baverez} were able to prove that  it holds at the level of leading asymptotics when $z\to 0$ (see also the discussion in \cite{baverez} on the relevance of this asymptotic in the context of the scaling limit of large random planar maps). The canonical Hilbert space  $\caH$ can also  be constructed  using the Osterwalder-Schrader reconstruction theorem \cite{OS1,OS2},  see  \cite{cargese} for lecture notes on this. The  representation of the Virasoro algebra on  $\caH$ should then follow using the results of \cite{KRV}. However, a probabilistic understanding of the highest weight vectors  is a challenge as they seem to involve vertex operators with complex weights  $\alpha$ whereas the probabilistic approach naturally deals with real  $\alpha$. Also, the main application of LCFT to Liouville Quantum Gravity involves real values for $\alpha$. In the theory of Liouville Quantum Gravity, the scaling limits of e.g. Ising correlations on a random planar map are given in terms of Liouville correlations with  real  $\alpha$'s and regular planar Ising CFT correlations via the celebrated KPZ relation \cite{cf:KPZ}: for an explicit mathematical conjecture, see \cite{DKRV2,cargese}. Thus the probabilistic and bootstrap approaches are in an interesting way complementary.
  The bootstrap idea has been extremely successful  since the work by BPZ in \cite{BPZ}  and has led to spectacular progress even in three dimensions, e.g. in case of the 3d Ising model \cite{Ry1,Ry2}. A proof of \eqref{eq:cb} would be the first mathematical justification of this idea in a nontrivial and interesting CFT  and we consider it to be a major challenge to probabilists.

\section{Theorem on the reflection coefficient} \label{Reflection coefficient} 
A key ingredient in our derivation of the DOZZ formula is the reflection coefficient. We will see later that it plays a prominent role in the analyticity properties of correlation functions. Briefly, the reason is that expectations of the type \eqref{ThextendedSeibergbounds0} are analytic in $s$ over a region determined by the tail asymptotics of the random variable $ \int_{\C}  F(x,{\bf z}) M_\gamma(d^2x)$, which is in turn completely determined by the behaviour of this integral close to the ``worst" singularity of $F$.
The reflection coefficient enters in the description of the tail of such random variables.

\subsection{Tail behaviour of chaos integrals} To motivate the definitions let us consider the random variable  
\begin{equation} \label{Ialpha}
I(\alpha):=\int_{B(0,1)}|x|^{-\gamma\alpha}M_\gamma(d^2x).
\end{equation}
In the case $\alpha \in (\frac{\gamma}{2},Q)$, the reflection coefficient enters in the tail behaviour of $ I(\alpha)$ whose mass is concentrated around $0$\footnote{When studying the tail behaviour for $\alpha< \frac{\gamma}{2}$, the mass of $I(\alpha)$ is distributed on $B(0,1)$ and not concentrated around a point.} and which is a power law  as we now explain.  To study this we recall basic material introduced in \cite{DMS} and in particular we consider the polar decomposition of the chaos measure. Let $X_s:=X_s(0)$ be the circle average \eqref{circleav}. We have
\begin{equation*} 
X(e^{-s} e^{i \theta})=X_s+Y(s,\theta)
\end{equation*}
where $X_s$ is a standard Brownian Motion starting from the origin at $s=0$ and $Y$ is an independent field with  covariance 
\begin{equation}\label{covlateral}
\E[  Y(s,\theta) Y(t,\theta') ]  = \ln \frac{e^{-s}\vee e^{-t}}{|e^{-s}e^{i \theta} - e^{-t} e^{i \theta'} |}.
\end{equation}
Following \cite{DMS}, we call the field $Y$ the lateral noise. We also introduce the chaos measure with respect to $Y$ 
\begin{equation}\label{chaosmeasure}
N_\gamma(ds d\theta)  =  e^{\gamma Y(s,\theta)-\frac{ \gamma^2 \E[Y(s,\theta)^2]}{2}}  ds d\theta.
\end{equation}
Then we get
\begin{equation}\label{ialphadef}
I(\alpha)\stackrel{law}=\int_0^\infty  e^{\gamma (B_s-(Q-\alpha)s)}Z_s ds 
\end{equation}
with
\begin{equation}\label{DefZ}
Z_s = \int_0^{2 \pi} e^{\gamma Y(s,\theta)-\frac{\gamma^2 \E[Y(s,\theta)^2]}{2}} d\theta.  
\end{equation}
This is a slight abuse of notation since the process $Z_s$ is not a function (for $\gamma \geq \sqrt{2}$) but rather a generalized function. With this convention, notice that $Z_s ds$ is stationary i.e. for all $t$ the equality $Z_{t+s}=Z_s$  holds in distribution.

 It satisfies for all bounded intervals $I$ (see \cite{review})
\begin{equation}\label{zmoment}
\E \left [ \left ( \int_I Z_s ds \right )^p \right ]<\infty, \quad -\infty<p<\frac{4}{\gamma^2}.
\end{equation}

The following decomposition lemma due to Williams (see \cite{Williams})  will be useful in the study of $I(\alpha)$:

\begin{lemma}\label{lemmaWilliams}
Let $(B_s-\nu s)_{s \geq 0}$ be a Brownian motion with negative drift, i.e. $\nu >0$ and let $M= \sup_{s \geq 0} (B_s-\nu s)$.  Then conditionally on $M$ the law of the path   $(B_s-\nu s)_{s \geq 0}$ is given by the joining of two  independent paths:

\begin{itemize}
\item
A Brownian motion $((B_s^1+\nu s))_{s \leq \tau_M}$ with positive drift $\nu >0$ run until its hitting time $\tau_M$ of $M$. 
\item
$(M+{B}^2_t- \nu t)_{t\geq 0}$ where ${B}^2_t- \nu t$  is a  Brownian motion with negative drift conditioned to stay negative.
\end{itemize}

Moreover, one has the following time reversal property for all $C>0$ (where $\tau_C$ denotes the hitting time of $C$) 
\begin{equation*}
(B_{\tau_C-s}^1+\nu (\tau_C-s)-C)_{s  \leq \tau_C}\stackrel{law}=  (\tilde{B}_s- \nu s)_{s \leq L_{-C}}
\end{equation*}
where $(\tilde{B}_s- \nu s)_{s \geq 0}$ is a Brownian motion with drift $-\nu$ conditioned to stay negative and $ L_{-C}$ is the last time  $(\tilde{B}_s- \nu s)$ hits $-C$.  

\end{lemma}

\begin{remark}\label{fundamentalremark}
As a consequence of the above lemma, one can also deduce that the process $(\tilde{B}_{L_{-C}+s}- \nu (L_{-C}+s)+C)_{s \geq 0}$ is equal in distribution to $(\tilde{B}_s- \nu s)_{s \geq 0}$. 
\end{remark}
 
 This lemma motivates defining the process $ \mathcal{B}^\alpha_s$ 
\begin{equation*}
 \mathcal{B}^\alpha_s = \left\{
 \begin{array}{ll}
  B^\alpha_{-s} & \text{if } s < 0\\
    \bar{B}^\alpha_{s} & \text{if } s >0 \end{array} \right.
\end{equation*}
where $B^{\alpha}_s,\bar B^{\alpha}_s$ are  two independent Brownian motions with negative drift $\alpha-Q$ and conditioned to stay negative. We may apply Lemma \ref{lemmaWilliams} to \eqref{ialphadef}. Let $M=\sup_{s\geq 0}(B_s-(Q-\alpha)s)$ and $L_{-M}$ be the last time $(B^\alpha_s)_{s \geq 0}$ hits $-M$. Then
\begin{align}\label{willaplied}
&\int_0^\infty e^{\gamma (B_s-(Q-\alpha)s)}Z_s ds\stackrel{law}=e^{\gamma M}\int_{-L_{-M}}^\infty e^{\gamma\caB_s^{\alpha}}Z_{s+L_{-M}}ds\stackrel{law}=e^{\gamma M}\int_{-L_{-M}}^\infty e^{\gamma\caB_s^{\alpha}}Z_{s} ds
\end{align}
where we used stationarity of the process $Z_s$ (and independence of $Z_s$ and $B_s$). We will prove in section \ref{tailestimates} that the tail behaviour of $I(\alpha)$ coincides with that of
\begin{equation*}
J(\alpha)= e^{\gamma M }   \int_{-\infty}^\infty  e^{\gamma \caB^\alpha_s}  Z_s ds.
\end{equation*}
The distribution of $M$ is well known (see section 3.5.C in the textbook \cite{KaraSh} for instance):
\begin{equation}\label{tailofmax}
 \P (e^{\gamma M}>x ) = \frac{1}{x^{\frac{2(Q-\alpha)}{\gamma}}} ,  \quad  x\geq 1
\end{equation}
 which implies
\begin{equation}\label{tailofmax1}
\P   ( J(\alpha) >x  ) \underset {x \to \infty}{\sim}   \E  \left  [  \left ( \int_{- \infty} ^\infty  e^{\gamma \caB^\alpha_s}  Z_s ds  \right )^{\frac{2(Q-\alpha)}{\gamma}} \right ]\, x^{-\frac{2(Q-\alpha)}{\gamma}}.
\end{equation}
This is the tail behaviour that we prove for $I(\alpha)$ and its generalizations in section \ref{tailestimates}. Define the {\it unit volume reflection coefficient} $\bar{R}(\alpha)$  for $\alpha \in (\frac{\gamma}{2},Q)$ by the following formula
\begin{equation}\label{defunitR}
 \bar{R}(\alpha)=  \E \left [ \left ( \int_{-\infty}^\infty  e^{   \gamma \mathcal{B}_s^\alpha  } Z_s ds  \right )^{\frac{2}{\gamma}  (Q-\alpha)} \right ].
\end{equation}
 $\bar{R}(\alpha)$  is  indeed well defined as can be seen from the following lemma, the proof of which is postponed to  Appendix \ref{app:momR} (see also section 4 in \cite{DMS} for the case $\alpha=\gamma$).

\begin{lemma}\label{defmomR}
Let $\alpha \in (\frac{\gamma}{2},Q)$. Then
\begin{equation}\label{positivemomentR}
\E \left [ \left ( \int_{-\infty}^\infty  e^{   \gamma \mathcal{B}_s^\alpha  } Z_s ds  \right )^{p}  \right ] < \infty
\end{equation}
for all $-\infty<p < \frac{4}{\gamma^2}$.
\end{lemma}

The full reflection coefficient is now defined for all $\alpha \in (\frac{\gamma}{2},Q) \setminus \cup_{n \geq 0} \lbrace \frac{2}{\gamma} -\frac{n}{2} \gamma \rbrace $ by 
\begin{equation}\label{deffullR}  
R(\alpha)= \mu^{\frac{2(Q-\alpha)}{\gamma}}   \Gamma (-\tfrac{2(Q-\alpha)}{\gamma})  \tfrac{2(Q-\alpha)}{\gamma} \bar R(\alpha). 
\end{equation}
The function $R(\alpha)$ has a divergence at the points $\frac{2}{\gamma}- \frac{n}{2} \gamma $ with $n \geq 0$ because of the $\Gamma$ function entering the definition. Its connection to the structure constants is the following (see the proof in Section \ref{sec:proofreflection}):

\begin{lemma}\label{defR}
For all $\alpha \in (\frac{\gamma}{2},Q)  \setminus \cup_{n\geq 0} \lbrace \frac{2}{\gamma} -\frac{n}{2} \gamma \rbrace $, the following limit holds
\begin{equation*}
\underset{\epsilon \to 0}{\lim} \:  \epsilon \: C_\gamma(\epsilon, \alpha,\alpha)= 4 R(\alpha).
\end{equation*}
\end{lemma}
Hence the reflection coefficient should be seen as a 2 point correlation function. Let us mention that the fact that some form of 2 point correlation function should exist in LCFT goes back to Seiberg \cite{seiberg}.  
\subsection{Main result on the reflection coefficient}

The second main result of this paper is the following exact formula for  the reflection coefficient (recall \eqref{defRDOZZ}): 
\begin{theorem}\label{Rtheor}
For all $\alpha \in (\frac{\gamma}{2},Q)$ one has
\begin{equation}\label{keyformula}
R(\alpha)= R^{{\rm DOZZ}}(\alpha).
\end{equation}
\end{theorem}

\section{BPZ equations and their consequences} \label{sec:further}
In this section we collect some previous results from the companion paper \cite{KRV} that will be used in the proof of Theorems \ref{theoremDOZZ} and \ref{Rtheor}.

\subsection{Structure constants and four point functions}\label{Structure constants and four point functions}
We complete the description of the three point structure constants \eqref{expression3pointstruct} by introducing the unit volume three point structure constants defined by the formula
\begin{equation}\label{unitvolumethreepoint}   
\bar{C}_\gamma(\alpha_1,\alpha_2,\alpha_3)= \mu^s  \frac{C_\gamma(\alpha_1,\alpha_2,\alpha_3)}{\Gamma(s)}.
\end{equation}
 where $s=(\sum_{i=1}^3\alpha_i-2Q)/\gamma$.
 The four point function (equation \eqref{Z1} with $N=4$) is fixed by the M\"obius invariance \eqref{KPZformula} up to a single function depending on the cross ratio of the points. For later purpose we label the insertion points from $0$ to $3$ and consider the weights $\alpha_1, \alpha_2, \alpha_3$ fixed:
\begin{align}
 \nonumber \langle    
   \prod_{k=0}^3 V_{\alpha_k}(z_k)  \rangle 
 & = |z_3-z_0|^{- 4 \Delta_{0}}  |z_2-z_1|^{ 2 (\Delta_3-\Delta_2-\Delta_1-\Delta_0)}|z_3-z_1|^{2(\Delta_2+\Delta_0 -\Delta_3 -\Delta_1  )} \\ &\times |z_3-z_2|^{2 (\Delta_1+\Delta_0-\Delta_3-\Delta_2)} G_{\alpha_0}\left ( \frac{(z_0-z_1)(z_2-z_3)}{ (z_0-z_3) (z_2-z_1)}  \right )  \label{confinv}
\end{align}
 where here again we use the shortcut notation $\Delta_j:=\Delta_{\alpha_j}$. 
We can recover $G_{\alpha_0}$ as the following limit
\begin{equation}
G_{\alpha_0}(z)= \lim_{ z_3 \to\infty }|z_3|^{4 \Delta_3} \langle    V_{\alpha_0}  (z)  V_{\alpha_1}(0)  V_{\alpha_2}(1) V_{\alpha_3}(z_3)  \rangle   \label{Glimit0}.
\end{equation}
Combining with \eqref{Z1}  we get
 \begin{equation}
  G_{\alpha_0}(z)=|z|^{-\alpha_0 \alpha_1}   |z-1|^{-\alpha_0 \alpha_2 }  \mathcal{T}_{\alpha_0}(z) 
   \label{Glimit1}
\end{equation}
where, setting $s=\frac{\alpha_0+\alpha_1+\alpha_2+\alpha_3-2Q}{\gamma}$, $ \mathcal{T}_{\alpha_0}
  (z)$ is given by 
\begin{align}\label{Tdefi0}
  \mathcal{T}_{\alpha_0}(z) &=
  2 \mu^{-s}  \gamma^{-1}\Gamma(s)  \E [ \mathcal{S}_{\alpha_0}(z)^{-s}]
 \end{align}
and
\begin{equation}\label{defcalS}
\mathcal{S}_{\alpha_0}(z)= 
\int_{\C} 
\frac{|x|_+^{\gamma\sum_{k=0}^3 \alpha_k}}{ |x-z|^{\gamma\alpha_0}|x|^{\gamma \alpha_1}  |x-1|^{\gamma \alpha_2}  }  M_\gamma(d^2x).
\end{equation}
In this paper we will study the structure constants \eqref{expression3pointstruct} by means of  four point functions \eqref{confinv} with special values of $\alpha_0$.

\subsection{BPZ equations}\label{BPZ equations}
There are two special values of $\alpha_0$ for which the reduced  four point function $ \mathcal{T}_{\alpha_0}(z)$ satisfies a second order differential equation. That such equations are expected in Conformal Field Theory goes back to BPZ \cite{BPZ}. In the  case of LCFT it was proved in \cite{KRV}  that, under suitable assumptions on $\alpha_1,\alpha_2,\alpha_3$, if $\alpha_0\in\{-\frac{\gamma}{2},-\frac{2}{\gamma}\}$  then $\mathcal{T}_{\alpha_0}$ is a   solution of a PDE version of the Gauss hypergeometric equation
\begin{equation}\label{hypergeo}
\partial_{z}^2 \caT_{\alpha_0}(z)+  \frac{({ c}-z({ a}+{ b}+1))}{z(1-z)}\partial_z \caT_{\alpha_0}(z) -\frac{{ a}{ b} }{z(1-z)}\caT_{\alpha_0}(z)=0
\end{equation}
where ${ a},{ b},{ c}$ are given by
\begin{align}\label{defabcfirst}
{ a}&= \frac{_{\alpha_0}}{^2} (Q-2\alpha_0-\alpha_1-\alpha_2-\alpha_3)-\hf,  \quad { b}=\frac{_{\alpha_0}}{^2} (Q-\alpha_1-\alpha_2+\alpha_3)+\hf,  \quad c=1+\alpha_0 (Q-\alpha_1).
\end{align}
This equation has two holomorphic solutions defined on  $\mathbb{C} \setminus \lbrace (-\infty,0) \cup (1,\infty) \rbrace$:
 \begin{equation}\label{Fpmdef}
F_{-}(z)= {}_2F_1({ a},{  b},{ c},z), \quad F_{+}(z)= z^{1-{ c}} {}_2F_1(1+{ a}-{ c},1+{ b}-{ c},2-{ c},z)
\end{equation}
where $_2F_1(a,b,c,z)$ is given by the standard hypergeometric series (which can be extended holomorphically on $\mathbb{C} \setminus  (1,\infty) $). Using the facts that $ \caT_{\alpha_0}(z)$ is real, single valued and $C^2$ in $\mathbb{C} \setminus \{0,1\}$ we proved in \cite{KRV}  (Lemma 4.4) that it is determined up to a multiplicative constant $\lambda\in\R$ as
\begin{equation}\label{Tsolution}
\caT_{\alpha_0}(z)=  
\lambda( | F_{-}(z) |^2+A_\gamma(\alpha_0,\alpha_1 , \alpha_2, \alpha_3) | F_{+}(z) |^2)
\end{equation}
where  the coefficient $A_\gamma(\alpha_0,\alpha_1 , \alpha_2, \alpha_3)$ is given by  
\begin{equation}\label{Fundrelation}
A_\gamma(\alpha_0,\alpha_1 , \alpha_2, \alpha_3)=- \frac{\Gamma(c)^2  \Gamma(1-a)  \Gamma(1-b)  \Gamma(a-c+1)  \Gamma(b-c+1) }{ \Gamma(2-c)^2  \Gamma(c-a)  \Gamma(c-b)  \Gamma (a)  \Gamma(b) }
\end{equation}
where recall that $a,b,c$ are defined in terms of $\alpha_0,\alpha_1 , \alpha_2, \alpha_3$ by \eqref{defabcfirst} and provided
$c \in \R \setminus \Z$ and $c-a-b \in  \R \setminus \Z$.
Furthermore, the constant $\lambda$ is found by using  the expressions \eqref{expression3pointstruct} and \eqref{Tdefi0} (note that $s$ has a different meaning in these two expressions):
\begin{equation}\label{Fundrelation1}\lambda=\caT_{\alpha_0}(0)=C_\gamma(\alpha_1+\alpha_0,\alpha_2,\alpha_3).
\end{equation}
Hence for $\alpha_0\in\{-\frac{\gamma}{2},-\frac{2}{\gamma}\}$  $\mathcal{T}_{\alpha_0}$ is completely
determined in terms of $C_\gamma(\alpha_1+\alpha_0,\alpha_2,\alpha_3)$. 

In the case $\alpha_0=-\frac{\gamma}{2}$ we were able to determine  in \cite[Lemma 4.5]{KRV}  the leading asymptotics of the expression \eqref{Tdefi0} as $z\to 0$ {\it provided }  $\frac{1}{\gamma}+\gamma<\alpha_1+\frac{\gamma}{2}<Q$ : 
\begin{equation}\label{Fundrelation300}
\caT_{-\frac{\gamma}{2}}(z)= C_\gamma(\alpha_1-\frac{\gamma}{2},\alpha_2,\alpha_3)+ B(\alpha_1)C_\gamma(\alpha_1+\frac{\gamma}{2},\alpha_2,\alpha_3) | z |^{2(1-c)}+o( | z |^{2(1-c)})
\end{equation}
where
\begin{equation}\label{Fundrelation2}
B(\alpha)=-\mu  \frac{\pi}{  l(-\frac{\gamma^2}{4}) l(\frac{\gamma \alpha}{2})  l(2+\frac{\gamma^2}{4}- \frac{\gamma \alpha}{2}) } .
\end{equation}
In view of \eqref{Tsolution} and the fact that $2(1-c)<1$ (since $2(1-c)= \gamma (Q-\alpha_1)$ and $\frac{1}{\gamma}+\frac{\gamma}{2}<\alpha_1$), we also have the following expansion around $z=0$ 
\begin{equation*}
\caT_{-\frac{\gamma}{2}}(z)= C_\gamma(\alpha_1-\frac{\gamma}{2},\alpha_2,\alpha_3)+ A_\gamma(-\frac{\gamma}{2}, \alpha_1 , \alpha_2, \alpha_3)C_\gamma(\alpha_1-\frac{\gamma}{2},\alpha_2,\alpha_3) | z |^{2(1-c)}+o( | z |^{2(1-c)}).
\end{equation*}
By unicity of the Taylor expansion around $z=0$ we get
\begin{equation}\label{Fundrelation30}
B(\alpha_1)C_\gamma(\alpha_1+\frac{\gamma}{2},\alpha_2,\alpha_3)=A_\gamma(-\frac{\gamma}{2}, \alpha_1 , \alpha_2, \alpha_3)C_\gamma(\alpha_1-\frac{\gamma}{2},\alpha_2,\alpha_3).
\end{equation}
Let us now register the following relation between $A_\gamma$ and $\mathcal{A}_\gamma$ which stems from straightforward (but lengthy!) algebra
\begin{equation}\label{relationAotherA}
\frac{A_\gamma (-\frac{\gamma}{2}, \alpha_1 , \alpha_2, \alpha_3)}{B(\alpha_1)}= -\frac{1}{\pi \mu} \mathcal{A}_\gamma (\frac{\gamma}{2}, \alpha_1 , \alpha_2, \alpha_3).
\end{equation}   
Let us also register here (by anticipation of the case $\alpha_0=-\frac{2}{\gamma}$) the following analogue dual relation 
\begin{equation}\label{relationAotherAdual}
\frac{A_\gamma (-\frac{2}{\gamma}, \alpha_1 , \alpha_2, \alpha_3)}{\tilde{B}(\alpha_1)}= -\frac{1}{\pi \tilde{\mu}} \mathcal{A}_\gamma (\frac{2}{\gamma}, \alpha_1 , \alpha_2, \alpha_3)
\end{equation}  
where $\tilde{\mu}= \frac{(\mu \pi l(\frac{\gamma^2}{4})  )^{\frac{4}{\gamma^2} }}{ \pi l(\frac{4}{\gamma^2})}$ and 
\begin{equation}\label{Fundrelation2dual}
\tilde{B}(\alpha)=-\tilde{\mu}  \frac{\pi}{  l(-\frac{4}{\gamma^2}) l(\frac{2 \alpha}{\gamma})  l(2+\frac{4}{\gamma^2}- \frac{2 \alpha}{\gamma}) }. 
\end{equation}
 Therefore thanks to \eqref{Fundrelation30} and \eqref{relationAotherA} we get relation \eqref{3pointconstanteqintro} in the case  $\frac{1}{\gamma}+\gamma<\alpha_1+\frac{\gamma}{2}<Q$ and also 
\begin{equation}\label{Fundrelation3}
\caT_{-\frac{\gamma}{2}}(z)= C_\gamma(\alpha_1-\frac{\gamma}{2},\alpha_2,\alpha_3) | F_{-}(z) |^2+ B(\alpha_1)C_\gamma(\alpha_1+\frac{\gamma}{2},\alpha_2,\alpha_3) | F_{+}(z) |^2.
\end{equation}
The restriction   $\frac{1}{\gamma}+\gamma<\alpha_1+\frac{\gamma}{2}$  for $\alpha_1$ was technical in \cite{KRV} and will be removed in section \ref{bpzeq}. The restriction $\alpha_1+\frac{\gamma}{2}<Q$ seems necessary due to the Seiberg bounds as the probabilistic $C_\gamma(\alpha_1+\frac{\gamma}{2},\alpha_2,\alpha_3)$ vanishes otherwise. Understanding  what happens when  $\alpha_1+\frac{\gamma}{2}>Q$  is the key to our proof of the DOZZ formula. Before turning to this we  draw a useful corollary from the results of this section.

\subsection{Crossing relation}\label{Crossing relation}
Let us assume  the validity of the Seiberg bounds for the four point correlation function with weights  $(-\frac{\gamma}{2},\alpha_1,\alpha_2,\alpha_3)$, i.e.  $\sum_{k=1}^3 \alpha_k>2Q+\frac{\gamma}{2}$ and $\alpha_k<Q$ for all $k$.  
We have from the previous subsection
\begin{equation}\label{Fundrelation4}
\caT_{-\frac{\gamma}{2}}(z)= C_\gamma(\alpha_1-\frac{\gamma}{2},\alpha_2,\alpha_3) (| F_{-}(z) |^2+ A_\gamma(-\frac{\gamma}{2},\alpha_1 , \alpha_2, \alpha_3) | F_{+}(z) |^2).
\end{equation}
The hypergeometric equation \eqref{hypergeo} has another basis of holomorphic solutions defined on  $\mathbb{C} \setminus \lbrace (-\infty,0) \cup (1,\infty) \rbrace$: 
\begin{equation}\label{Fpmdef'}
G_{-}(z)= {}_2F_1({ a},{  b},{ c'},1-z), \quad G_{+}(z)=(1- z)^{1-{ c'}} {}_2F_1(1+{ a}-{ c'},1+{ b}-{ c'},2-{ c'},1-z)
\end{equation} 
where $c'=1+a+b-c=1-\frac{\gamma}{2}(Q-\alpha_2)$  (i.e. these are obtained by  interchanging $\alpha_1$ and $\alpha_2$ and replacing $z$ by $1-z$). The two basis are linearly related
\begin{align}
F_{-}(z)&= \frac{\Gamma (c) \Gamma (c- a- b) }{\Gamma (c- a) \Gamma (c- b)} G_-(z)   + \frac{\Gamma (c) \Gamma ( a+ b- c)  }{\Gamma ( a) \Gamma ( b)}  (1-z)^{c-a-b}G_+(z)\label{f+}
\\
F_{+}(z)&= \frac{\Gamma (2- c) \Gamma ( c- a- b) }{\Gamma (1- a) \Gamma (1- b)}  G_-(z)  + \frac{\Gamma (2-c) \Gamma ( a+ b- c)  }{\Gamma ( a- c+1) \Gamma ( b- c+1)} (1-z)^{c-a-b}G_+(z)  \label{f-}
\end{align}
and we get 
\begin{equation}\label{newformula}
\caT_{-\frac{\gamma}{2}}(z)= C_\gamma(\alpha_1-\frac{\gamma}{2},\alpha_2,\alpha_3) (D| G_{-}(z) |^2+ E | G_{+}(z) |^2)
\end{equation}
with explicit coefficients $D,E$ (see \cite{KRV}, Appendix).
On the other hand by studying the asymptotics as $z\to 1$ we get 
\begin{equation}\label{Fundrelation5}
\caT_{-\frac{\gamma}{2}}(z)= C_\gamma(\alpha_1,\alpha_2-\frac{\gamma}{2},\alpha_3)+ B(\alpha_2)C_\gamma(\alpha_1,\alpha_2+\frac{\gamma}{2},\alpha_3) | 1-z |^{2(1-c')}+o( |1- z |^{2(1-c')}).
\end{equation}
More precisely, this asymptotic has been  established in \cite{KRV} under the restriction  $\frac{1}{\gamma}+\gamma<\alpha_2+\frac{\gamma}{2}<Q$, which is empty for $\gamma^2>2$! We anticipate here Section \ref{bpzeq} and Theorem  \ref{theo4point} where the validity of \eqref{Fundrelation5} (or a version of \eqref{Fundrelation5} with extra $1-z$ and $1-\bar{z}$ terms in the expansion when $2(1-c')>1$) will be relaxed to the range of parameters $\gamma<\alpha_2+\frac{\gamma}{2}<Q$, which is non empty whatever the value of $\gamma<2$.

Comparing the $z\to1$ expansion of \eqref{newformula} with \eqref{Fundrelation5} leads to the following crossing symmetry relation: 
\begin{proposition}\label{propcrossing1} 
Let $\alpha_2+\frac{\gamma}{2}<Q$ and  $\alpha_1+ \alpha_2+\alpha_3 -\frac{\gamma}{2}>2Q$. Then  \begin{equation} \label{YetaFundamentalrelation}
 C_\gamma (\alpha_1-\frac{\gamma}{2}, \alpha_2, \alpha_3)= T(\alpha_1,\alpha_2,\alpha_3)  C_\gamma (\alpha_1,\alpha_2 +\frac{\gamma}{2}, \alpha_3)
\end{equation}
where
$T$ is  given by the following formula 
\begin{equation} \label{defT0}
T(\alpha_1,\alpha_2,\alpha_3)= - \mu \pi \frac{l(a) l(b)}{ l(c) l(a+b-c)}    \frac{1}{l(-\frac{\gamma^2}{4}) l(\frac{\gamma \alpha_2}{2})  l(2+\frac{\gamma^2}{4}- \frac{\gamma \alpha_2}{2})}
\end{equation}
with $a,b,c$ given by \eqref{defabcfirst} for $\alpha_0=-\frac{\gamma}{2}$.
\end{proposition}
The statement in the proposition above  should be further restricted to $\gamma<\alpha_2+\frac{\gamma}{2}<Q $ according to the previous discussion. However, here we anticipate Theorem \ref{analytic_hyperbolic} in Section \ref{sec:analytic} to extend by analyticity our statement to the range of parameters as formulated above. 
\begin{remark}
The relations  \eqref{3pointconstanteqintro} and \eqref{YetaFundamentalrelation} were derived in the physics literature \cite{Tesc} by assuming (i) validity of BPZ equations for degenerate field insertions, (ii) that these fields are given by the vertex operators with weights $-\frac{\gamma}{2}, -\frac{2}{\gamma}$  and that they satisfy an appropriate operator product expansion (iii) the diagonal form of the solution \eqref{Tsolution} (iv) crossing symmetry (an essential input in the bootstrap approach). We want to stress that our proof makes no such assumptions, in fact  (i)-(iv) are theorems: (i) follows from integration by parts in the Gaussian measure \cite{KRV}, (ii) follows from an asymptotic analysis of the probabilistically defined four point functions as points are ``fused'' together in Section \ref{bpzeq} (iii) and (iv) follow from an analysis of the  BPZ equations and proof of regularity of solutions \cite{KRV}.
\end{remark}

\section{Strategy and plan of proof} \label{sec:bone}
In this section,  we outline  our strategy for the proof while  giving pointers to the remaining parts of the paper so that the reader can have a better view of the whole structure. We will first explain the proof of Theorem \ref{Rtheor} which gives an explicit expression for the reflection coefficient. This exact expression for $R$ is then used (in an essential way) to derive Theorem \ref{theoremDOZZ} on the DOZZ formula.   

\subsection{Proof of Theorem \ref{Rtheor}: analysis of the reflection coefficient}
The proof of Theorem \ref{Rtheor} is gathered in Section \ref{proofreflection}. Recall that the reflection coefficient $R(\alpha)$ is defined by \eqref{deffullR}.  The proof is based on establishing the following properties:
\begin{description}
\item[ARC1]  the unit volume reflection  coefficient $\bar{R}(\alpha)$ (see \eqref{defunitR} for the definition) defines an analytic function of $\alpha$ over a complex neighborhood of the interval $(\frac{\gamma}{2},Q)$.
\item[ARC2] $R(\alpha)$ satisfies the following $\frac{\gamma}{2}$-shift equation for $\alpha$ close to but smaller than $ Q$
\begin{equation}\label{sketch:shift1}
R(\alpha-\frac{\gamma}{2})= - \mu \pi \frac{R(\alpha)}{ l(-\frac{\gamma^2}{4}) l(\frac{\gamma\alpha}{2}-\frac{\gamma^2}{4})  l(2+\frac{\gamma^2}{4}- \frac{\gamma \alpha}{2})}.
\end{equation}
This relation allows us to extend analytically $R(\alpha)$ to a meromorphic function, {\bf still denoted by $R$}, over a complex neighborhood of the real line; this complex neighborhood contains $\R \times (-\eta,\eta)$ for some $\eta>0$.
\item[ARC3] $R(\alpha)$ satisfies the inversion relation
\begin{equation}\label{sketch:inversion}
R(\alpha)R(2Q-\alpha)=1.
\end{equation}
\item[ARC4] $R(\alpha)$  satisfies the $\frac{2}{\gamma} $-shift equation
\begin{equation}\label{sketch:shift2}
R(\alpha)=   - c_\gamma \frac{R(\alpha+\frac{2}{\gamma})}{l(-\frac{4}{\gamma^2})l(\frac{2 \alpha}{\gamma}) l(2+\frac{4}{\gamma^2}-\frac{2 \alpha}{\gamma})} 
\end{equation}
where $c_\gamma = \frac{\gamma^2}{4} \mu \pi R(\gamma) \not = 0$.
\end{description}
%
According to Liouville's theorem, if a continuous function $f$ has two periods $T_1$ and $T_2$ such that $\frac{T_2}{T_1} \not\in \Q$ then $f$ is a constant function. Therefore, provided $\gamma^2 \not \in \Q $, the two equations \eqref{sketch:shift1} and \eqref{sketch:shift2} fully determine $R$ up to a $\gamma$-depending constant (the inversion relation \eqref{sketch:inversion} is used in the proof of {\bf ARC4}), which we determine easily by computing $R(Q)=-1$ via the probabilistic representation \eqref{defunitR} and \eqref{deffullR}. In particular it determines the value of $c_\gamma$
\begin{equation}\label{shift2a}
 c_\gamma={(\mu \pi l(\tfrac{\gamma^2}{4})  )^{\frac{4}{\gamma^2} }}{l(\tfrac{4}{\gamma^2})}^{-1}.
\end{equation}
On the other hand in the Appendix (see relations \eqref{shift1DOZZ} and \eqref{shift3DOZZ}), we check that $R^{{\rm DOZZ}}$ satisfies \eqref{sketch:shift1} and \eqref{sketch:shift2} with $c_{\gamma}$ given by \eqref{shift2a}.
Therefore we conclude $R$ is equal to $R^{{\rm DOZZ}}$ for $\gamma^2 \not \in \Q $; the general case can be deduced by continuity in $\gamma$ of $R$ and $R^{{\rm DOZZ}}$. In what follows, we now give an idea of how we prove {\bf ARC1}-{\bf ARC4}. 



\subsubsection{Analyticity of correlation functions and the reflection coefficient}
As most of our arguments are based on analyticity properties, we first need to show that the probabilistically defined correlation functions  \eqref{unitvol} are analytic in a complex neighborhood of the real valued parameters $(\alpha_k)_k$ delimited by the extended Seiberg bounds \eqref{ThextendedSeibergbounds}. This is done in Section \ref{sec:analytic} but is restricted to $N$-point correlation functions with $N\geq 3$. The argument is based on the fact that regularized correlation functions are analytic in the parameters $(\alpha_k)_k$ and converge locally uniformly over a complex neighborhood of the extended Seiberg bounds. 

The case of the reflection coefficient (or 2 point correlation function) requires more insight as it is not clear how to choose a regularized version that is analytic in the parameter $\alpha$ and converges nicely towards $R$. The main idea is to interpret the reflection coefficient as the leading order coefficient in the tail expansion of the random variable 
\begin{align}\label{rhooo} \rho(\alpha_1,\alpha_2,\alpha_3)=  
 \int_{\C}  \frac{|x|_+^{\gamma(\alpha_1+\alpha_2+\alpha_3)}}{ |x|^{\gamma \alpha_1}  |x-1|^{\gamma \alpha_2}  }    M_\gamma(d^2x) 
 \end{align}
  involved in the probabilistic representation \eqref{expression3pointstruct} of the 3 point structure constant. The reason for that relies on a general simple argument: assume we are given a positive random variable $X$ with tail asymptotics given by
$$\P(X>t)=\frac{c_{\star}}{t^\beta}+o(t^{-\beta-\delta})$$
for some $\beta,\delta>0$. Then the function $s\mapsto\E[X^s]$ is analytic over a complex neighborhood of $(0,\beta)$. Furthermore it extends to a meromorphic function over a complex neighborhood of $(0,\beta+\delta)$ with  a pole at $s=\beta$, given by $\frac{c_\star s}{\beta-s}$. One can then recover the value of $c_\star$ by taking the limit
$$\lim_{s\to\beta}(\beta-s)\E[X^s]=\beta c_\star.$$
This type of argument will be repeatedly used in the paper (and in fact even pushed further to the next pole of $\E[X^s]$ beyond $s=\beta$).

As the integral in \eqref{rhooo}  can be decomposed as a sum of singular GMC integrals of the type  \eqref{Ialpha} (around $0$,$1$ and $\infty$), a detailed study of the tail of such singular GMC integrals needs first to be carried out. This is the content of Section  \ref{tailestimates} where $R$ emerges in the tail expansion of integrals of the type \eqref{Ialpha}.  The first outcome of this study is the proof of Lemma \ref{defR} in Section \ref{sec:proofreflection}. Actually we even prove a stronger result (Proposition \ref{defR1}) that the reflection coefficient can be recovered from the structure constant as the following limit for $\tfrac{\gamma}{2}<\alpha_2\leq\alpha_3<Q$
\begin{equation}\label{bridge}
\lim_{\alpha_1\searrow \alpha_3-\alpha_2}(\alpha_1-\alpha_3+\alpha_2) C_\gamma (\alpha_1, \alpha_2, \alpha_3)=
a {R}(\alpha_3).
\end{equation}
where $a=2$ if  $\alpha_2<\alpha_3$ and $a=4$ if $\alpha_2=\alpha_3$. Recalling the interpretation by Seiberg \cite{seiberg} of the reflection coefficient as a two point structure constant, let us call \eqref{bridge} the ``$3\to2$"-bridge.

The second step is then to use the ``$3\to2$"-bridge in relations involving three point structure constants in order to produce  relations on the reflection coefficient $R(\alpha)$.  The type of relations we have in mind are the crossing symmetries of the type exposed in Proposition \ref{propcrossing1}. As a first example, Subsection \ref{Step1} explains how we use this bridge in the  crossing symmetry relation of Proposition \ref{propcrossing1} in order to express $\bar{R}(\alpha)$ as a function of three point structure constants for $\alpha\in (\frac{\gamma}{2}, Q)$ 
\begin{equation}\label{keyrelationonR}
\bar{R}(\alpha)= -   \frac{_\pi }{^\gamma}   \frac{   l(\frac{\gamma}{2} \alpha -\frac{\gamma^2}{4}-1)}{   l(1+\frac{\gamma}{2}(\alpha-Q)) l(-\frac{\gamma^2}{4}) l(\frac{\gamma^2}{4})     }  
\bar{C}_\gamma (\alpha,\gamma, \alpha).
\end{equation}
From relation \eqref{keyrelationonR}, which is new even with respect to the physics literature (to the best of our knowledge), we can deduce analyticity of $\bar{R}(\alpha)$ in  $\alpha\in (\frac{\gamma}{2}, Q)$ as stated in item \textbf{ARC1} above because now we know from Section \ref{sec:analytic}   that $\bar{C}_\gamma (\alpha,\gamma, \alpha)$ is analytic in $\alpha$. In conclusion,  analyticity of $\bar{R}(\alpha)$ seems very difficult to prove directly so we rely on relation \eqref{keyrelationonR} and analyticity of the 3 point structure constants to prove it.

\subsubsection{Exploiting the BPZ equations}   
 As explained above the ``$3\to2$"-bridge reduces \textbf{ARC2}-\textbf{ARC4} to deriving relations involving the three point structures constants.  The flavour of the derivation of these relations has already been explained in subsections   \ref{BPZ equations} and \ref{Crossing relation}. But establishing Theorem \ref{Rtheor} involves generalizing the relations explained in subsection \ref{BPZ equations} which were proved in the companion paper \cite{KRV}: this is the content of Section \ref{bpzeq}.

The first task  is to extend the range of parameters for which the relation \eqref{Fundrelation3} was established in \cite[Theorem 2.3]{KRV}. The reason for the restriction to the range of parameters $\tfrac{1}{\gamma}+\gamma<\alpha_1+\frac{\gamma}{2}<Q$ (or equivalently $\tfrac{1}{\gamma}+\tfrac{\gamma}{2}<\alpha_1<\tfrac{2}{\gamma}$) in  \cite[Theorem 2.3]{KRV} was technical: it relies on the asymptotic expansion \eqref{Fundrelation300} of $\caT_{-\frac{\gamma}{2}}(z)$ as $z\to 0$ in order to identify the constants  in front of the hypergeometric functions $|F_+|^2$ and $|F_-|^2$ in the general form 
 \eqref{Tsolution} of solutions to the $\tfrac{\gamma}{2}$-BPZ equation. This leads to \eqref{Fundrelation3}.  Within this range of values of the parameter $\alpha_1$, computing the first two leading terms of the expansion were enough since $2(1-c)= \gamma (Q-\alpha_1)<1$(recall that $c=1-\frac{\gamma}{2}(Q-\alpha_1)$ in this context).  Notice that the admissible set of values $\alpha_1$ satisfying the relation $\tfrac{1}{\gamma}+\tfrac{\gamma}{2}<\alpha_1<\tfrac{2}{\gamma}$ is empty for $\gamma^2\geq 2$! This is clearly not enough.

So, in Theorem \ref{theo4point}, we establish  \eqref{Fundrelation3} for the extended range of parameters $\gamma<\alpha_1+\frac{\gamma}{2}<Q$ (or equivalently $ \tfrac{\gamma}{2}<\alpha_1<\tfrac{2}{\gamma}$ ), which is non empty whatever the value of $\gamma<2$. In the situation $ \tfrac{\gamma}{2}< \alpha_1<\frac{1}{\gamma}+\frac{\gamma}{2}$, we have  $2(1-c)= \gamma (Q-\alpha_1)>1$ and therefore the expansion around $z=0$ that is required to prove \eqref{Fundrelation3} involves extra terms of the form $z$ and $\bar{z}$ (which correspond to the expansion of $|F_-|^2$ around $z=0$). Then analyticity of correlation functions entails the validity of  \eqref{Fundrelation3} for whatever value of $\alpha_1$ such that the correlation functions involved in  \eqref{Fundrelation3} satisfy the Seiberg bounds. As a consequence we obtain the crossing symmetry relation as stated in Proposition  \ref{propcrossing1}. The proof of Theorem \ref{theo4point} requires a refined version of fusion estimates compared to those proved in \cite[Section 5]{KRV}: this is the content of Lemma \ref{fusion4}. Equating residues on both sides of the relation \eqref{YetaFundamentalrelation} with the help of the ``$3\to2$"-bridge produces the $\frac{\gamma}{2}$-shift equation  as claimed in item \textbf{ARC2}; this is proved in Subsection \ref{Step2}.

 \color{black}
 
\begin{figure}[h] 
\centering
\subfloat[Domain of validity of the $\tfrac{\gamma}{2}$-shift equation]{\begin{tikzpicture}[xscale=0.75,yscale=0.75] 
\tikzstyle{cons}=[minimum width=4cm,minimum height=0.8cm,rectangle,rounded corners=5pt,draw,fill=red!75,text=black,font=\bfseries,text centered,text badly centered];
\shade[top color=Cerulean,bottom color=Apricot] (4,1) rectangle +(8,4) node[midway]{Theorem  \ref{theo4point} or \cite[Theorem 2.3]{KRV}};
\node at (8,2){depending on $\gamma^2\geq 2$ or $\gamma^2<2$};
\shade[top color=white,bottom color=PineGreen] (16,1) rectangle +(3,4) node[midway]{Theorem  \ref{theo4point1}};
\tikzstyle{fleche}=[->,>= stealth,thick,color=red!75];
\node[cons] (C) at (12,7) {under the global Seiberg constraint $\alpha_1+\alpha_2+\alpha_3>2Q+\frac{\gamma}{2}$};
\draw[fleche] (C.190) -- (6,4.7);
\draw[fleche] (C.240) -- (10,4.7);
\draw[fleche] (C.380) -- (17.5,4.7);
\node (Q) at (19,1){ $\bullet$};
\draw (19,1) node[below]{$Q$} ;
\node (D) at (16,1){ $\bullet$};
\draw (16,1) node[below]{$Q-\eta$} ;
\node (Q2) at (2.5,5){ $\bullet$};
\draw (2.5,5) node[left]{$Q$} ;
\node (B) at (12,1){ $\bullet$};
\draw (12,1) node[below]{$\frac{2}{\gamma}$} ;
\node (A) at (4,1){ $\bullet$};
\draw (4,1) node[below]{$  \tfrac{\gamma}{2}$} ;
\draw[>=latex,->,thick] (2,1) -- (20,1) node[right]{$\alpha_1$};
\draw[>=latex,->,thick] (2.5,0.5) -- (2.5,5.5) node[above]{$\alpha_2,\alpha_3$};
\draw[dashed,thick] (2.5,5) -- (20,5);
\end{tikzpicture}
}

\subfloat[Domain of validity of the $\tfrac{2}{\gamma}$-shift equation]{\begin{tikzpicture}[xscale=0.75,yscale=0.75] 
\tikzstyle{sommet}=[circle,draw,fill=yellow,scale=0.4] 
\tikzstyle{cons}=[minimum width=4cm,minimum height=0.8cm,rectangle,rounded corners=5pt,draw,fill=red!75,text=black,font=\bfseries,text centered,text badly centered];
\tikzstyle{fleche}=[->,>= stealth,thick,color=red!75];
\shade[top color=white,bottom color= Rhodamine] (16,4) rectangle +(3,1) node[midway]{ Theorem \ref{theo4point2overgamma}};
\node[cons] (C) at (12,2.5) {under the global Seiberg constraint $\alpha_1+\alpha_2+\alpha_3>2Q+\frac{2}{\gamma}$};
\draw[fleche] (C.10) -- (18,3.9);
\node (Q) at (19,1){ $\bullet$};
\draw (19,1) node[below]{$Q$} ;
\node (D) at (16,1){ $\bullet$};
\draw (16,1) node[below]{$Q-\eta$} ;
\node (Q2) at (2.5,5){ $\bullet$};
\draw (2.5,5) node[left]{$Q$} ;
\node (B) at (12,1){ $\bullet$};
\draw (12,1) node[below]{$\frac{2}{\gamma}$} ;
\node (A) at (8,1){ $\bullet$};
\draw (8,1) node[below]{$ \tfrac{1}{\gamma}+\tfrac{\gamma}{2}$} ;
\node (A) at (4,1){ $\bullet$};
\draw (4,1) node[below]{$  \tfrac{\gamma}{2}$} ;
\draw[>=latex,->,thick] (2,1) -- (20,1) node[right]{$\alpha_1$};
\draw[>=latex,->,thick] (2.5,0.5) -- (2.5,5.5) node[above]{$\alpha_2,\alpha_3$};
\draw[dashed,thick] (2.5,5) -- (20,5);
\end{tikzpicture}}

\subfloat[How the shift equations connect different ranges of parameters (global Seiberg constraint is assumed though not indicated)]{\begin{tikzpicture}[xscale=0.75,yscale=0.75] 
\tikzstyle{cons}=[minimum width=4cm,minimum height=0.8cm,rectangle,rounded corners=5pt,draw,fill=red!75,text=black,font=\bfseries,text centered,text badly centered];
\fill[top color=Cerulean,bottom color=Apricot,color=Apricot] (4,1) rectangle +(4,3);
\fill[top color=Cerulean,bottom color=Apricot,color=Apricot] (14,1) rectangle +(4,3);
\fill[top color=white,bottom color=PineGreen] (12,1) rectangle +(1,3);
\fill[top color=white,bottom color=PineGreen] (18,1) rectangle +(1,3);
\fill[top color=white,bottom color=Rhodamine] (10,3) rectangle +(1,1);
\fill[top color=white,bottom color=Rhodamine] (18,3) rectangle +(1,1);
\draw[>=latex,<->,thick] (10.5,3.5) -- (18.5,3.5);
\draw[>=latex,<->,thick] (7,2.5) -- (15,2.5);
\draw[>=latex,<->,thick] (12.5,1.5) -- (18.5,1.5);
\node (Q) at (19,1){ $\bullet$};
\draw (19,1) node[below]{$Q+\eta$} ;
\node (Q3) at (18,1){ $\bullet$};
\draw (18,1) node[below]{$Q$} ;
\node (Q) at (13,1){ $\bullet$};
\draw (13,1) node[below]{$\frac{2}{\gamma}$} ;
\node (D) at (12,1){ $\bullet$};
\draw (12,1) node[below]{$\frac{2}{\gamma}-\eta$} ;

\node (Q) at (11,1){ $\bullet$};
\draw (11,1) node[below]{$\frac{\gamma}{2}$} ;
\node (D) at (10,1){ $\bullet$};
\draw (10,1) node[below]{$\frac{\gamma}{2}-\eta$} ;

\node (Q2) at (2.5,4){ $\bullet$};
\node (Q) at (13,1){ $\bullet$};
\draw (13,1) node[below]{$\frac{2}{\gamma}$} ;
\node (D) at (12,1){ $\bullet$};
\draw (12,1) node[below]{$\frac{2}{\gamma}-\eta$} ;
\node (Q2) at (2.5,4){ $\bullet$};
\draw (2.5,4) node[left]{$Q$} ;
\node (B) at (8,1){ $\bullet$};
\draw (8,1) node[below]{$\tfrac{2}{\gamma}- \tfrac{\gamma}{2}$} ;
\node (B2) at (14,1){ $\bullet$};
\draw (14,1) node[below]{$\gamma$} ;
\node (A) at (4,1){ $\bullet$};
\draw (4,1) node[below]{$0$} ;
\draw[>=latex,->,thick] (2,1) -- (20,1) node[right]{$\alpha_1$};
\draw[>=latex,->,thick] (2.5,0.5) -- (2.5,4.5) node[above]{$\alpha_2,\alpha_3$};
\draw[dashed,thick] (2.5,4) -- (20,4);
\end{tikzpicture}
}
\caption{Domain of validity of the shift equations}
\label{fig:shift}
\end{figure}


\begin{figure}[h] 
\centering
 \begin{tikzpicture}[xscale=0.64,yscale=0.64]  
 \tikzstyle{noeud}=[minimum width=2cm,text width=3cm,minimum height=0.8cm,rectangle,rounded corners=5pt,draw,fill=yellow!75,text=red,font=\bfseries,text centered,text badly centered]
 \tikzstyle{core}=[minimum width=2cm,text width=3cm,minimum height=0.8cm,rectangle,rounded corners=5pt,draw,fill=green!75,text=red,font=\bfseries,text centered,text badly centered]
  \tikzstyle{main}=[minimum width=2cm,text width=3cm,minimum height=0.8cm,rectangle,rounded corners=5pt,draw,fill=red!75,text=white,font=\bfseries,text centered,text badly centered]
  \tikzstyle{sommet}=[circle,draw,fill=black]
\tikzstyle{pas}=[thick]  
\tikzstyle{fleche}=[->,>= stealth,thick]
\node[noeud] (A) at (0,11) {Analyticity  Section \ref{sec:analytic}};
\node[noeud] (T) at (13,11) {Tail analysis  Section \ref{tailestimates}};
\node[noeud] (B) at (7,8.8) {``$3\to 2$"-bridge  Section \ref{sec:proofreflection}};
\node[noeud] (F1) at (7,4) {$\tfrac{\gamma}{2}$-fusion  with resulting crossing relation Th   \ref{theo4point} \textcolor{red}{\& Pr}  \ref{propcrossing1again}};
\node[noeud] (FR) at (1,4.5) {$\tfrac{\gamma}{2}$-fusion with reflection and resulting crossing relation Th   \ref{theo4point1} \textcolor{red}{\& Pr} \ref{2overgamma}};
\node[noeud] (F2) at (13,4.5) {$\tfrac{2}{\gamma}$-fusion  and resulting crossing relation Th\ref{theo4point2overgamma} \textcolor{red}{\& Pr} \ref{2overgamma}};
\node[core] (C1) at (0,0) {Analyticity of $\bar{R}$ Subsection \ref{Step1}};
\node[core] (C3) at (7,0) {Gluing Lemma  Subsection  \ref{Gluing}};
\node[core] (C2) at (17,0) {$\tfrac{\gamma}{2}$-shift for $R$ Subsection  \ref{Step2}};
\node[core] (C4) at (11,-3) {Inversion for $R$ Subsection \ref{Step3}};
\node[core] (C5) at (7,-6) {$\tfrac{2}{\gamma}$-shift for $R$ Subsection \ref{Step4}};
\node[main] (PR) at (8,-9) {\Large Value of $R$};

\draw[pas] (F1) |- (12.5,2)edge[out=90,in=90] (13.5,2);
\draw[fleche] (13.5,2) -| (C2);
\draw[pas] (F1) |- (4.5,2)edge[out=90,in=90] (3.5,2);
\draw[pas] (3.5,2) -- (1.5,2)edge[out=90,in=90] (0.5,2);
\draw[fleche] (0.5,2) -| (C1);

\draw[fleche] (T) |- (12,10) -| (B);
\draw[fleche] (12,10) -- (19,10) -- (19,3) -| (C2);
\draw[pas] (A) |- (3.9,7.5);
\draw[pas] (3.9,6.5) -- (3.9,7.5);
\draw[pas] (B) |- (3.9,7.5);
\draw[pas] (3.9,6.5) -- (3.9,1.5);
\draw[fleche] (3.9,6.5) -- (-3,6.5) |- (C1);
\draw[pas] (F1) |- (3.9,1.5);
\draw[pas] (FR) |- (3.9,1.5);
\draw[fleche] (F1) -- (C3);
\draw[fleche] (3.9,6.5) -| (C2);
\draw[pas] (C1) |-(4,-1.5);
\draw[pas] (C3) |-(4,-1.5);
\draw[pas] (F2) |-(4,-1.5);
\draw[fleche] (4,-1.5) -| (C4);
\draw[pas] (4,-1.5) |- (8,-4.7);
\draw[pas] (C4) |- (8,-4.7);
\draw[fleche]  (8,-4.7) -| (C5);
\draw[pas] (C1) |- (0,-1) -- (0,-6) |- (8,-7.3) ;
\draw[pas] (C5) |- (8,-7.3);
\draw[pas] (C4) |- (12,-4.3) |- (8,-7.3);
\draw[pas] (C2) |- (8,-7.3);
\draw[fleche]  (8,-7.3) -| (PR);
\end{tikzpicture}
\caption{Diagram of the proof of Theorem \ref{Rtheor}} 
\label{fig:diagram}
\end{figure}

Another important task is to understand the analog  of   \eqref{Fundrelation3} when violating the Seiberg bounds, in particular when $\alpha_1< Q$ but $\alpha_1+\frac{\gamma}{2}\geq Q$. Mass concentration effects, like those involved in computing the tail of singular GMC integrals, will make the reflection coefficient play a prominent role in this context.  In Theorem  \ref{theo4point1}, we will show that for $\alpha_1$ close to $Q$ (but smaller  then $Q$) 
\begin{equation}\label{sketch:BPZ1}
\mathcal{T}_{-\frac{\gamma}{2}}(z) = C_\gamma(\alpha_1-\frac{\gamma}{2}, \alpha_2,\alpha_3)  |F_{-}(z)|^2  + R(\alpha_1) C_\gamma(2Q-\alpha_1-\frac{\gamma}{2}, \alpha_2,\alpha_3)  |F_{+}(z)|^2.
\end{equation}
As an output we prove the gluing Lemma \ref{gluinglemma}, which roughly states that the mapping
\begin{equation}
S(\alpha):= \begin{cases} \label{gluinglemmaintro}
C_\gamma(\alpha, \alpha_2,\alpha_3), \; \text{if} \: \alpha < Q \\
R(\alpha) C_\gamma(2Q-\alpha, \alpha_2,\alpha_3), \; \text{if} \: \alpha > Q \\
\end{cases}
\end{equation}
is   holomorphic in a neighborhood of $Q$. Using the ``$3\to2$"-bridge, this lemma will be instrumental in proving the inversion relation \textbf{ARC3} and the $\frac{2}{\gamma}$-shift equation \textbf{ARC4}.

Finally our final task is to investigate the consequences of the $\frac{2}{\gamma}$-BPZ equation \eqref{hypergeo}. By studying asymptotics as $z\to 0$ in \eqref{Tdefi0}, we show in  Theorem \ref{theo4point2overgamma} that for $\alpha_1$ close to $Q$ (but smaller than $Q$)
\begin{equation}\label{sketch:BPZ2}
\mathcal{T}_{-\frac{2}{\gamma}} (z) = C_\gamma(\alpha_1-\frac{_2}{^\gamma}, \alpha_2,\alpha_3)  |F_{-}(z)|^2  + R(\alpha_1) C_\gamma(2Q-\alpha_1-\frac{_2}{^\gamma}, \alpha_2,\alpha_3)  |F_{+}(z)|^2.
\end{equation}
 This produces new crossing relations as stated in Proposition \ref{2overgamma}. In subsection \ref{Step3}, we prove  the inversion relation stated in item \textbf{ARC3} by combining  the crossing relation in Proposition \ref{2overgamma} with the gluing lemma.

The $\frac{2}{\gamma}$-shift equation stated in item \textbf{ARC4} is  established in subsection \ref{Step4} first by continuing analytically the crossing relation \eqref{Firstinverserelation} to some larger set of values with the help of the gluing Lemma, and then by
equating residues in  both sides of the resulting relation with the help of the inversion relation.

\subsection{Proof of Theorem \ref{theoremDOZZ}: the DOZZ formula}   \label{proofDOZZ}
In order to prove the DOZZ formula, we first want to prove the shift equations  \eqref{3pointconstanteqintro} and \eqref{3pointconstanteqintrodual}.  As explained in subsection \ref{BPZ equations}, the shift relation   \eqref{3pointconstanteqintro} is a consequence of the identity 
\begin{equation}\label{remindshiftgamma}
\caT_{-\frac{\gamma}{2}}(z)= C_\gamma(\alpha_1-\frac{\gamma}{2},\alpha_2,\alpha_3) | F_{-}(z) |^2+ B(\alpha_1)C_\gamma(\alpha_1+\frac{\gamma}{2},\alpha_2,\alpha_3) | F_{+}(z) |^2.
\end{equation}
In order to derive the other shift equation \eqref{3pointconstanteqintrodual}, we need to exploit the  $\frac{2}{\gamma}$-BPZ equation \eqref{hypergeo}. We show in Theorem \ref{theo4point2overgamma} that for $\alpha_1$ close to $Q$ (but smaller than $Q$)
\begin{equation}\label{sketch:BPZ2}
\mathcal{T}_{-\frac{2}{\gamma}} (z) = C_\gamma(\alpha_1-\frac{_2}{^\gamma}, \alpha_2,\alpha_3)  |F_{-}(z)|^2  + R(\alpha_1) C_\gamma(2Q-\alpha_1-\frac{_2}{^\gamma}, \alpha_2,\alpha_3)  |F_{+}(z)|^2.
\end{equation}
In order to exploit \eqref{sketch:BPZ2}, we need to use in a crucial way analyticity of $S$ defined by  \eqref{gluinglemmaintro} (gluing lemma) along with  item  \textbf{ARC4} on the reflection coefficient. Indeed, thanks to item  \textbf{ARC4}, we have
\begin{equation*}
R(\alpha_1) C_\gamma(2Q-\alpha_1-\frac{_2}{^\gamma}, \alpha_2,\alpha_3)= \tilde{B}(\alpha_1) R(\alpha_1+\frac{2}{\gamma}) C_\gamma(2Q-\alpha_1-\frac{_2}{^\gamma}, \alpha_2,\alpha_3)
\end{equation*} 
where $\tilde{B}$ was defined by \eqref{Fundrelation2dual}. Thanks to the gluing lemma, we have
\begin{equation*}
  R(\alpha_1+\frac{2}{\gamma}) C_\gamma(2Q-\alpha_1-\frac{_2}{^\gamma}, \alpha_2,\alpha_3)=C_\gamma(\alpha_1+\frac{_2}{^\gamma}, \alpha_2,\alpha_3) .
\end{equation*} 
Therefore we can rewrite \eqref{sketch:BPZ2} equivalently as 
 \begin{equation}\label{sketch:BPZ2again}
\mathcal{T}_{-\frac{2}{\gamma}} (z) = C_\gamma(\alpha_1-\frac{_2}{^\gamma}, \alpha_2,\alpha_3)  |F_{-}(z)|^2  +  \tilde{B}(\alpha_1) C_\gamma(\alpha_1+\frac{_2}{^\gamma}, \alpha_2,\alpha_3)  |F_{+}(z)|^2.
\end{equation}
From \eqref{sketch:BPZ2again}, we can derive the shift equation \eqref{3pointconstanteqintrodual}  along the same lines as we derived the shift equation  \eqref{3pointconstanteqintro}  from \eqref{remindshiftgamma}.

Once  the shift equations  \eqref{3pointconstanteqintro} and \eqref{3pointconstanteqintrodual} are proved, the proof of Theorem \ref{theoremDOZZ} is again a consequence of Liouville's theorem on periodic functions.  In order to prove Theorem \ref{theoremDOZZ}, we suppose that $\gamma^2\not\in \Q$, the other case resulting from a continuity argument in $\gamma$. For $\gamma^2\not\in \Q$, because $C_\gamma^{{\rm DOZZ}}(\alpha_1,\alpha_2,\alpha_3)$ satisfies the same shift equations  \eqref{3pointconstanteqintro} and \eqref{3pointconstanteqintrodual}, this implies by application of Liouville 's theorem in the variable $\alpha_1$ that the ratio $\frac{C_\gamma(\alpha_1,\alpha_2,\alpha_3)}{C_\gamma^{{\rm DOZZ}}(\alpha_1,\alpha_2,\alpha_3)}$ is independent of $\alpha_1$. By symmetry of the $\alpha_1,\alpha_2,\alpha_3$ variables and recursive use of Liouville 's theorem (in the variable $\alpha_2$ and then the variable $\alpha_3$), we deduce that the quotient $\frac{C_\gamma(\alpha_1,\alpha_2,\alpha_3)}{C_\gamma^{{\rm DOZZ}}(\alpha_1,\alpha_2,\alpha_3)}$ only depends on $\gamma$. We identify this constant by using the ``$3\to2$"-bridge and the knowledge of $R$. This argument is detailed in Section \ref{proofDOZZ}.

\section{Analytic Continuation of Liouville Correlation Functions}\label{sec:analytic}
In this section we study the analytic continuation of the unit volume correlations \eqref{unitvol}. These are defined for real weights 
$\bm{\alpha}=(\alpha_1,\dots,\alpha_N)$ satisfying the extended Seiberg bounds
\begin{equation}\label{ThextendedSeibergbounds1}
U_N:=\{\bm{\alpha}\in\R^N:\frac{_1}{^\gamma}(2Q-\sum_{k=1}^N\alpha_k)
< \frac{_4}{^{\gamma^2}} \wedge \min_{1 \leq k \leq N}  \frac{_2}{^\gamma}(Q-\alpha_k),\ \ \forall k:\alpha_k<Q\}.
\end{equation}
We will prove:

\begin{theorem}\label{analytic_hyperbolic}
Fix $N\geq 3$ and distinct points $z_1,\dots,z_N\in \C^N$. The unit volume correlation function  \eqref{unitvol} extends to an analytic function of $\bm{\alpha}$ defined in a complex neighborhood of $U_N$ in $\C^N$.
\end{theorem}

\proof 
By M\"obius invariance we may assume $|z_i|>2$ and $|z_i-z_j|>2$. 
We use  \eqref{Z00} to write the unit volume correlation functions as the limit  
 \begin{align}\label{Z000}
   \langle \prod_{k=1}^N  V_{\alpha_k}(z_k)  \rangle_{\mathrm{uv}}& =  2   \gamma^{-1} \prod_{k=1}^N  |z_k|_+^{-4\Delta_{\alpha_k}}  \lim_{r\to \infty}F_r(\bm{\alpha})
   \end{align}
   where
 \begin{align}\label{Z000}
    F_r(\bm{\alpha})=
   \E  \left [  \prod_{k=1}^N e^{{\alpha_k} X_{r}(z_k)- \frac{\alpha_k^2}{2}
  \E  X_{r}(z_k)^2}  
   M_\gamma(\C_r) ^{-s}    \right ]
\end{align}
and  $\C_r:=\C\setminus \cup_{k=1}^N B(z_k,e^{-r})$. $F_r$ is  defined for all $\bm{\alpha}\in\C^N$ and   is complex differentiable in $\alpha_i$, hence defines an entire function in the $\alpha_i$. We show that there is an open $V\subset\C^N$ containing $U_N$ s.t. $F_r$  converges uniformly on compacts of $V$. Note that this is nontrivial since for $\alpha_k=a_k+ib_k$ we have
$$
|e^{{\alpha_k} X_{r}(z_k)- \frac{\alpha_k^2}{2}
    \E[ X_r(z_k)^2]} |= e^{a_kX_{r}(z_k)- \frac{a_k^2}{2}
    \E [X_r(z_k)^2]} e^{ \frac{b_k^2}{2}
  \E  X_r(z_k)^2} 
$$
and $e^{ \frac{b_k^2}{2}
  \E  X_{r}(z_k)^2} \propto e^{ \frac{b_k^2}{2}r}$ blows up as $r\to \infty$.
    
 By Remark \ref{indebm}, we know that    $t\in\R_+\to B_{r+t}^k:=X_{r+t}(z_k)-X_{r}(z_k)$ are mutually independent Brownian motions and  they are independent of  $\sigma\{X(x);x\in\C_r\}$. Hence 
\begin{align*}
    F_{r+1}(\bm{\alpha})-F_{r}(\bm{\alpha})=
   \E  \left [  \prod_{k=1}^N e^{{\alpha_k} X_{r+1}(z_k)- \frac{\alpha_k^2}{2}
  \E [ X_{r+1}(z_k)^2]}  
  ( M_\gamma(\C_{r+1}) ^{-s}  -M_\gamma(\C_{r}) ^{-s}  )  \right ].
\end{align*}

 Now we apply the Cameron-Martin theorem as in \eqref{Z1} to the real parts of the vertex insertions  to get
 \begin{align}\label{Z000}
   |  F_{r+1}(\bm{\alpha})-F_{r}(\bm{\alpha})|\leq Ce^{(r+1)\sum_{k=1}^N\frac{b_k^2}{2}}
 |  \E   (\int_{\C_{r+1}} f_r (x)M_\gamma(d^2x)) ^{-s}  -\E(\int_{\C_{r}} f_r(x) M_\gamma(d^2x) ^{-s}   ) |
\end{align}
 where 
$f_r(x)=e^{\sum_{k=1}^N\gamma a_kG_{r+1}(x,z_k)} $
and we have defined $G_{r+1}(z,z'):=\E[X(z)X_{r+1}(z')]$. We get from \eqref{hatGformula}
$$f(x):=\sup_rf_r(x)\leq C\prod_k\big(\frac{|x|_+|z_k|_+}{|x-z_k|}\big)^{\gamma\alpha_k}.
$$
We need to estimate the difference of expectations in \eqref{Z000}. Let
\begin{align*}
Y_r:=\int_{\C_{r+1}\setminus \C_r} f_r(x) M_\gamma(d^2x)   
\end{align*}
and set also $Z_r:=\int_{\C_r} f_r(x) M_\gamma(d^2x)$. Then
 \begin{align*}
  & |  \E   (\int_{\C_{r+1}} f_r(x) M_\gamma(d^2x)) ^{-s}  -\E(\int_{\C_{r}} f_r(x) M_\gamma(d^2x)) ^{-s}   ) |= |  \E  ( (Z_r+Y_r) ^{-s}  -Z_r ^{-s}   ) |\\&\leq  \E [\ind_{Y_r<\epsilon} |  (Z_r+Y_r) ^{-s}  -Z_r^{-s} |]+  \E [\ind_{Y_r\geq\epsilon} |  (Z_r+Y_r) ^{-s}  -Z_r^{-s} |  ]
   \end{align*}
where $\epsilon>0$ will be fixed later. The first expectation on the RHS is bounded by
 \begin{align*}
   \E \ind_{Y_r<\epsilon} |  (Z_r+Y_r) ^{-s}  -Z_r^{-s} |\leq C\epsilon\sup_{t\in[0,1]}\E(Z_r+tY_r) ^{-\Re s-1} \leq C\epsilon
         \end{align*}
uniformly in $r$. The last bound follows by noting that for $-\Re s-1>0$ the expectation is bounded uniformly in $r$ by $C\E (\int f(x) M_\gamma(d^2x))^{-\Re s-1}$ which is finite due to \eqref{ThextendedSeibergbounds} whereas for $-\Re s-1<0$ we may bound it for example by  $C\E (\int_{\C_1\setminus\C_2} M_\gamma(d^2x))^{-\Re s-1}$ which is finite as well. 

For the second expectation we use the H\"older inequality 
\begin{align*}
   \E \ind_{Y_r\geq\epsilon} |  Z_{r+1}^{-s}  - Z_r^{-s} |\leq C\P(Y_r\geq\epsilon)^{1/p}((\E Z_{r+1}^{-q\Re s})^{1/q} +(\E Z_{r}^{-q\Re s})^{1/q}).
         \end{align*}
Taking $q>1$ s.t. $-q\Re(s)< \min_j\tfrac{2}{\gamma}(Q-\alpha_j)\wedge \tfrac{4}{\gamma^2}$ we may bound  the two expectations uniformly  in $r$ as in the previous paragraph and then using the Markov inequality we get
\begin{align*}
   \E \ind_{Y_r\geq\epsilon} |  Z_{r+1}^{-s}  - Z_r^{-s} |  \leq C\epsilon^{-m/p}( \E Y_r^m)^{1/p}.
      \end{align*}
It remains to bound $\E Y_r^m$ for suitable $m>0$.  We note that  $\C_{r+1}\setminus\C_r=\cup_iA _r^i$ where
 $A_r^i$ is the annulus centred at $z_i$ with radii $e^{-r-1},e^{-r}$. Then  we obtain for $m<\frac{4}{\gamma^2}$
\begin{align}\label{ymombound}
\E Y_r^m\leq C\E(\sum_k\int_{A_r^k} f(x)M_\gamma(d^2x))^m\leq C\max_ke^{-r(\gamma (Q-a_k) m-\frac{\gamma^2m^2}{2})}:=Ce^{-r\theta}
\end{align}
where in the second step we used the estimate \eqref{annulus}. Now, let us fix $\bm{a}^0\in U_N$. Then we can find $m>0$ and $\delta>0$  s.t. $\theta>0$ for all $\bm{a}$ with $\min_k|a_k-a^0_k|\leq\delta$.   Hence, for $\bm{\alpha}\in\C^M$ with
$\alpha_k=a_k+ib_k$  and $\epsilon>0$
 \begin{align*}
   | F_{r+1}(\bm{\alpha})-F_{r}(\bm{\alpha})|\leq Ce^{(r+1)\sum_{k=1}^N\frac{b^2_k}{^2}}(\epsilon+\epsilon^{-m/p}e^{-\frac{r}{p}\theta}).
   \end{align*}
Taking $\epsilon=e^{-\eta r}$ with $\eta=\frac{\theta}{p+m}$ we then have 
 \begin{align*}
   | F_{r+1}(\bm{\alpha})-F_{r}(\bm{\alpha})|\leq Ce^{-(\eta- \sum_{k=1}^N\frac{b^2_k}{^2})r}.
   \end{align*}
Hence, $F_{r}(\bm{\alpha})$ converges uniformly in a ball around $\bm{a}^0$ in $\C^N$.\qed

\section{Tail estimates for Multiplicative Chaos} \label{tailestimates} 
In this section, we prove the tail estimates needed in this paper and that involve the unit volume reflection coefficient defined in \eqref{defunitR}.

\subsection{Tail estimate around one insertion}

Let $|z|>2$ and consider the random variable  
\begin{equation*}
W:=\int_{B(z,1)}\frac{F(x')}{|x'-z|^{\gamma \alpha}}M_\gamma(d^2x')
\end{equation*}
for $F$ bounded and $C^1$ in a neighborhood of $z$. We assume $\frac{\gamma}{2}< \alpha < Q$ and define auxiliary quantities  $\beta=( \frac{2}{\gamma} (Q-\alpha) +\frac{2}{\gamma^2}) \wedge  \frac{4}{\gamma^2}$ and  $\bar\eta$ by $(1-\bar\eta) \beta= \frac{2}{\gamma}(Q-\alpha)+\bar\eta$. Hence $\bar\eta$ is strictly positive. With these definitions we have 
\begin{lemma}\label{simpletail}
For all $\eta<\bar\eta$ and for some constant $C(z)$, we have
\begin{equation*}
 | \P(  W  >x ) -   |z|^{4\alpha(\alpha-Q)}F(z)^{\frac{2}{\gamma} (Q-\alpha)} \frac{\bar{R} (\alpha) }{ x^{\frac{2}{\gamma} (Q-\alpha) }}    |  \leq \frac{C(z)}{x^{\frac{2}{\gamma} (Q-\alpha)  +\eta}  }.
\end{equation*}
\end{lemma} 

\proof
We will write the integral in polar coordinates of $B(z,1)$. Define
\begin{equation*}
N= \frac{1}{2 \pi}\int_{0}^{2\pi}  X (z+ e^{i\theta}) d\theta .  
\end{equation*}
Then
\begin{equation*}
B_s:= \frac{1}{2\pi} \int_{0}^{2\pi}   (X(z+e^{-s}e^{i \theta})  -X (z+ e^{i\theta}) )d \theta
\end{equation*}
is a Brownian motion with $B(0)=0$ and we may decompose the field $X$ as 
\begin{equation*}
X(z+x')= N +   B_{-\ln|x|}+ Y_z(x')
\end{equation*} 
where $Y_z$ is a lateral noise centered around $z$ given by
\begin{equation*}
Y_z(x')=  X(z+x')- \frac{1}{2\pi} \int_{0}^{2\pi}   X (z+|x'| e^{i \theta}) d \theta   .
\end{equation*}
 Notice that $Y_z$ has same distribution as the lateral noise $Y$ (centered around $0$), that $Y_z$ and  $ B $ are independent and  $N$ is independent of  $ B$.  We have
\begin{equation}\label{boundY}
| \E[Y_z(x')N]  |    =|\ln|z+x'|-\ln|z||\leq C |x'|
\end{equation}
(since $Y_z$ lives in the space of distributions, $\E[Y_z(x')N] $ is defined for all smooth function $f$ by the relation $\E [ (\int_{\C}f(x') Y_z(x') \, d^2x') N] = \int_{\C}f(x') \E[ Y_z(x')N] \,d^2x'$ ) and the variance of $N$ is 
$$\E[N^2]= 2 \ln |z|.
$$
Hence, we get the following decomposition into independent components 
\begin{equation}\label{decompind}
X(z+x')= \Big(1+ \frac{\E[Y_z(x')N]}{\E[N^2]} \Big)N + B_{-\ln|x'|}  + (Y_z(x') - \frac{\E[Y_z(x')N]}{\E[N^2]} N)   .
\end{equation} 
We introduce a variable $\bar{N}$ distributed as $N$ but independent of 
$N, B, Y_z$. We can rewrite \eqref{decompind} as the following equality in distribution:
\begin{equation}\label{decompinduseful} 
X(z+x')=  \Big(1+ \frac{\E[Y_z(x')N]}{\E[N^2]} \Big)\bar{N} +  B_{-\ln|x'|} + (Y_z(x') - \frac{\E[Y_z(x')N]}{\E[N^2]} N)  .
\end{equation}
Plugging this relation into the expression of $W$, we get 
$$W\overset{law}{=}e^{\gamma \bar{N} -\frac{\gamma^2}{2}\E[\bar N^2]}\int_{B(0,1)}u_z(x')e^{\gamma B_{-\ln |x'|}+(\frac{\gamma^2}{2}-\gamma\alpha)\ln|x'|}e^{\gamma Y_z(x')-\frac{\gamma^2}{2}\E[Y_z(x')^2]}\,d^2x'$$
for some (random) function $u_z$ such that (using \eqref{boundY} and $C^1$-regularity of  $F$  around $z$)
\begin{equation*}
|u_z(x')- \frac{F(z)}{ |z|^4 } 
|  \leq  C( 1+ e^{C |N|+C |\bar{N}|}  ) |x'|   .
\end{equation*}
We may thus write $W=W_1+W_2$ in distribution with
\begin{align}\label{Zdecomp}
 W_1 = & e^{\gamma \bar{N} -\frac{\gamma^2}{2}\E[\bar N^2]}\frac{F(z)}{ |z|^4 }\int_0^\infty e^{\gamma (B_s-(Q-\alpha)s)}Z_s ds \\
|W_2| \leq & C (1+e^{C (\bar{N}+N)} ) \int_0^\infty e^{\gamma (B_s-(Q-\alpha+\frac{1}{\gamma})s)}Z_s ds. 
\end{align}
and  $Z$, $B$ and $\bar{N}$ independent.

Recall now the Williams decomposition Lemma \ref{lemmaWilliams}. Let  $m=\sup_{s\geq 0}(B_s-(Q-\alpha+\frac{1}{\gamma})s)$ and let  $L_{-m}$ be the largest $s$ s.t. $ \mathcal{B}^{\alpha}_{-s}= -m$. Then
\begin{align}\nonumber
&\int_0^\infty e^{\gamma (B_s-(Q-\alpha+\frac{1}{\gamma})s)}Z_s ds\stackrel{law}=e^{\gamma m}\int_{-L_{-m}}^\infty e^{\gamma\caB_s^{\alpha-\frac{1}{\gamma}}}Z_{s+L_{-m}} ds\\&\stackrel{law}=e^{\gamma m}\int_{-L_{-m}}^\infty e^{\gamma\caB_s^{\alpha-\frac{1}{\gamma}}}Z_{s}ds\leq e^{\gamma m}\int_{-\infty}^\infty e^{\gamma\caB_s^{\alpha-\frac{1}{\gamma}}}Z_s ds\label{will1}
\end{align}
where we used stationarity of the process $Z_s$.

For all $0<p< ( \frac{2}{\gamma} (Q-\alpha) +\frac{2}{\gamma^2}) \wedge  \frac{4}{\gamma^2}=\beta$, we have  
\begin{equation}\label{Z'tail}
\P( |W_2|\geq x) \leq Cx^{-p}.
\end{equation}

Indeed, for all $p_1,q_1>1$ with $\frac{1}{p_1}+\frac{1}{q_1}=1$ we have by using H\"older and \eqref{will1} that
\begin{equation}\label{Z'tail}
\P( |W_2|\geq x) \leq \frac{1}{x^p}  \E[ |W_2|^p ] \leq \frac{C}{x^p}   \E[ e^{C (\bar{N}+N)pp_1}   ]^{1/p_1}   \E[ (  e^{\gamma m}\int_{-\infty}^\infty e^{\gamma\caB_s^{\alpha-\frac{1}{\gamma}}}Z_s ds   )^{p q_1}]^{1/q_1} \leq \frac{C}{x^p} 
\end{equation}
provided $q_1$ is sufficiently close to $1$ and where we used Lemma  \ref{defmomR} which requires  $p< \frac{4}{\gamma^2}$.

We first prove an upper bound for $\P( W> x)$. From 
\eqref{Z'tail} we get for $\eta\in (0,1)$:
\begin{align*}
 &    \P (     W >x  )=   \P (     W_1+W_2 >x  )   \leq   \P (
 W_1>x -x^{1 -\eta}  ) +C x^{-p(1-\eta)}  .
\end{align*}
Proceeding as in \eqref{will1} we get
\begin{align*}
\P(W_1>x -x^{1 -\eta}  ) \leq  \P ( e^{\gamma \bar{N} -\frac{\gamma^2}{2}\E[\bar N^2]}\frac{F(z)}{ |z|^4 }e^{\gamma M}\int_{-\infty}^\infty e^{\gamma \caB^\alpha_s}Z_s ds  >x -x^{1 -\eta}  ) 
\end{align*}
where  $M=\sup_{s\geq 0}(B_s-(Q-\alpha)s)$. In view of Lemma \ref{defmomR}, we have for all $\epsilon>0$ that  
\begin{equation*}
 \P ( e^{\gamma \bar{N} -\frac{\gamma^2}{2}\E[\bar N^2]}\frac{F(z)}{ |z|^4 }\int_{-\infty}^\infty e^{\gamma \caB^\alpha_s}Z_s ds  >x -x^{1 -\eta}  ) \leq \frac{C}{x^{\frac{4}{\gamma^2}-\epsilon}}
\end{equation*}
hence \eqref{tailofmax} implies for all $\epsilon>0$
\begin{align*}
 \P (      W >x  )  & \leq  e^{ (2(Q-\alpha)^2-\gamma(Q-\alpha))\E[\bar N^2]}\big(\frac{F(z)}{ |z|^4 }\big)^{\frac{2}{\gamma} (Q-\alpha)}\frac{ \bar{R} (\alpha) }{ (x-x^{1 -\eta} )^{\frac{2}{\gamma} (Q-\alpha) }} + C x^{-p(1-\eta)}\\&\leq |z|^{-4\alpha(Q-\alpha)}F(z)^{\frac{2}{\gamma} (Q-\alpha)}\frac{ \bar{R} (\alpha) }{ x^{\frac{2}{\gamma} (Q-\alpha) }}+Cx^{-\frac{2}{\gamma} (Q-\alpha) -\eta}+C x^{-p(1-\eta)} + Cx^{-\frac{4}{\gamma^2}+\epsilon}
 \end{align*}
 for $0<p<\beta$, for some constant $C$ which may depend on $\eta$ and $\epsilon$.
 Recall that we defined  $\bar\eta>0$ by $(1-\bar\eta) \beta= \frac{2}{\gamma}(Q-\alpha)+\bar\eta$.   We conclude
 \begin{align}\label{upb}
 \P (      W >x  )   
 \leq |z|^{-4\alpha(Q-\alpha)}F(z)^{\frac{2}{\gamma} (Q-\alpha)}\frac{ \bar{R} (\alpha) }{ x^{\frac{2}{\gamma} (Q-\alpha) }}+Cx^{-\frac{2}{\gamma} (Q-\alpha) -\eta}
  \end{align} 
 for all $\eta<\bar\eta$. 

 \vskip 2mm

Now, we consider the lower bound. We have
\begin{align} \label{easyineq}
 &  \P (   W  >x  )   \geq   \P( W_1 >x+ x^{1-\bar\eta}  ) -\P(W_2<- x^{1-\bar\eta})\geq \P( W_1 >x+ x^{1-\bar\eta}  ) -Cx^{-\frac{2}{\gamma} (Q-\alpha) -\eta}
\end{align}
for all $\eta<\bar\eta$. We define for all $C>0$ the random variable
\begin{equation*}
W(L_{-C})= e^{\gamma \bar{N} -\frac{\gamma^2}{2}\E[\bar N^2]}\frac{F(z)}{ |z|^4 }e ^{\gamma C} \int_{-L_{-C}}^\infty  e^{\gamma \mathcal{B}^{\alpha}_s} Z_{s} ds
\end{equation*}
 and by the  Williams decomposition we get as in \eqref{will1} 
\begin{equation*}
W_1 \overset{(law)}{=} e^{\gamma \bar{N} -\frac{\gamma^2}{2}\E[\bar N^2]}\frac{F(z)}{ |z|^4 }e ^{\gamma M} \int_{-L_{-M}}^\infty  e^{\gamma \mathcal{B}^{\alpha}_s} Z_{s} ds=W(L_{-M})
\end{equation*} 
where  $M=\sup_{s \geq 0} (B_s-(Q-\alpha)s)$ and $M$, $\mathcal{B}^{\alpha}$ and $Z_s$ are independent.

Let $\eta'$  be such that $(1-\eta') \frac{4}{\gamma^2} = \frac{2}{\gamma} (Q-\alpha)+\eta'$. One has $\eta' \geq \bar\eta$. 
Consider the event $\caE$ defined by 
$$e^{\gamma \bar{N} -\frac{\gamma^2}{2}\E[\bar N^2]}\frac{F(z)}{ |z|^4 }\int_{-L_{-M}}^\infty  e^{\gamma \mathcal{B}^{\alpha}_s} Z_{s} ds<x^{1-\eta'}.
$$
We have trivially
$$
\P( W_1 >x+ x^{1-\bar\eta}  )\geq \P( \{W_1 >x+ x^{1-\bar\eta}\}  \cap \caE ).
$$
Under  $\{W_1 >x+ x^{1-\bar\eta}\}\  \cap\ \caE $ we have $e^{\gamma M}\geq |x|^{\eta'}$. Indeed, if $e^{\gamma M}<|x|^{\eta'}$ then under $\caE$ we get $W_1<x$ which is impossible.  Thus  $M\geq -\frac{\eta'}{\gamma}\ln|x|$ whereby $L_{-M}\geq L_{-\frac{\eta'}{\gamma}\ln|x|}$ and hence $W(L_{-\frac{\eta'}{\gamma}\ln|x|})\leq W(L_{-M})$.  We conclude 
\begin{align}\nonumber
&\P( W_1 >x+ x^{1-\bar\eta}  )\geq \P( \{W(L_{-\frac{\eta'}{\gamma}\ln|x|}) >x+ x^{1-\bar\eta}\}  \cap  \caE )\\&\geq \P( W(L_{-\frac{\eta'}{\gamma}\ln|x|}) >x+ x^{1-\bar\eta})-Cx^{-(1-\eta')\frac{4}{\gamma^2}+\epsilon}\nonumber\\&\geq
 |z|^{-4\alpha(Q-\alpha)}F(z)^{\frac{2}{\gamma} (Q-\alpha)}
 \E[   ( \int_{-L_{-\frac{\eta'}{\gamma} \ln x}}^\infty  e^{\gamma \mathcal{B}^{\alpha}_s} Z_s ds  )^{\frac{2}{\gamma} (Q-\alpha)}]  (x+x^{1-\bar\eta})^{-\frac{2}{\gamma} (Q-\alpha) } - \frac{C}{x^{\frac{2}{\gamma} (Q-\alpha)+\eta'-\epsilon} } \label{lb1}
\end{align}
for all $\epsilon>0$ where in the second step we used Lemma  \ref{defmomR}. 

We claim now that 
\begin{equation}\label{inequapproxR}
  \E[   ( \int_{-\infty}^\infty  e^{\gamma \mathcal{B}^{\alpha}_s} Z_s ds  )^{\frac{2}{\gamma} (Q-\alpha)}]- \E[   ( \int_{-L_{-\frac{\eta'}{\gamma} \ln x}}^\infty  e^{\gamma \mathcal{B}^{\alpha}_s} Z_s ds  )^{\frac{2}{\gamma} (Q-\alpha)}]   \leq C{x^{-\eta'}}.
\end{equation}
Combined with   \eqref{lb1} and  \eqref{easyineq} this yields
 \begin{align}\nonumber
 \P (      W >x  )  & \geq   |z|^{-4\alpha(Q-\alpha)}F(z)^{\frac{2}{\gamma} (Q-\alpha)}
 \frac{ \bar{R} (\alpha) }{ (x-x^{1 -\eta} )^{\frac{2}{\gamma} (Q-\alpha) }} -Cx^{-\frac{2}{\gamma} (Q-\alpha) -\eta}\\&\geq |z|^{-4\alpha(Q-\alpha)}F(z)^{\frac{2}{\gamma} (Q-\alpha)}\frac{ \bar{R} (\alpha) }{ x^{\frac{2}{\gamma} (Q-\alpha) }}-Cx^{-\frac{2}{\gamma} (Q-\alpha) -\eta}\label{lpb2}
  \end{align}
for all $\eta<\bar\eta$.  \eqref{lpb2} and  \eqref{upb} then finish the proof.

It remains to prove \eqref{inequapproxR}.  By Remark \ref{fundamentalremark}, the process $\hat{B}_s^{\alpha}$ defined for $s \leq 0$ by the relation $\hat{B}^{\alpha}_s= \mathcal{B}^{\alpha}_{s-L_{-\frac{\eta'}{\gamma} \ln x} }+ \frac{\eta'}{\gamma} \ln x$ is independent from everything and distributed like $(\mathcal{B}^{\alpha}_s)_{s \leq 0}$. We can then write
\begin{equation*}
 \int_{-\infty}^\infty  e^{\gamma \mathcal{B}^{\alpha}_s} Z_s ds  =A +{x^{-\eta'}} B
\end{equation*}
where 
\begin{equation*}
A=\int_{-L_{-\frac{\eta'}{\gamma} \ln x}}^\infty  e^{\gamma \mathcal{B}^{\alpha}_s} Z_s ds\quad \text{ and }\quad 
B= \int_{-\infty}^0 e^{\gamma \hat{B}^{\alpha}_s}Z_{s- L_{-\frac{\eta'}{\gamma} \ln x}} ds.
\end{equation*}
We now distinguish two cases: $\frac{2}{\gamma} (Q-\alpha) \leq 1$ and $\frac{2}{\gamma} (Q-\alpha)  >1$. 
\vskip 1mm
\noindent  \underline{First case} $\frac{2}{\gamma} (Q-\alpha) \leq 1$. We use $(1+u)^{  \frac{2}{\gamma} (Q-\alpha) }-1 \leq \frac{2}{\gamma} (Q-\alpha) u$ for  $u \geq 0$ to bound
\begin{align*}
  \E   [   ( A+ {x^{-\eta'}} B )^{   \frac{2}{\gamma} (Q-\alpha) }- A^{   \frac{2}{\gamma} (Q-\alpha) }  ] 
   \leq \frac{_2}{^\gamma} (Q-\alpha) {x^{-\eta'}}  \E   [  {B}{A^{\frac{2}{\gamma} (Q-\alpha) -1}}  ] .
\end{align*}
By H\"older's inequality with $p \in (1,\frac{4}{\gamma^2})$, we get 
\begin{equation*}
 \E   [  {B}{A^{\frac{2}{\gamma} (Q-\alpha) -1}}  ] 
 \leq \E[   B^p  ]^{1/p}   \E[ A^{q(\frac{2}{\gamma} (Q-\alpha) -1)}]^{\frac{_1}{^q}} <\infty
\end{equation*}
since $B$ is equal in distribution to  $\int_{-\infty}^0 e^{\gamma \hat{B}^{\alpha}_s} Z_{s} ds$ and $A \geq \int_0^\infty e^{\gamma \mathcal{B}^{\alpha}_s} Z_s ds $ which has negative moments of all order by Lemma \ref{defmomR}.

\vspace{0.1 cm}

\noindent \underline{Second case} $\frac{2}{\gamma} (Q-\alpha_2) > 1$.  Let $p:=\frac{2}{\gamma} (Q-\alpha) $.  By triangle inequality we have 
\begin{align*}
&  \E   [   ( A+ {x^{-\eta'}} B )^{  p }- A^{   p}  ]  \leq   \left ( (\E   [   A^{p}   ])^{1/ p}  + {x^{-\eta'}}(\E [   B^{p}   ])^{1/ p}    \right )^{p} -    \E[A^{   p}  ]   \\
& \leq   \left (( \E   [   A^{p}   ])^{1/ p}  + {C}{x^{-\eta'}}   \right )^{p} -    \E[A^{   p}  ] \leq {C}{x ^{-\eta'}} \E[A^{   p}  ]^{1-1/ p}    \leq {C}{x ^{-\eta'}} 
\end{align*}
where again we used that $A$ and $B$ have moments of order $p$.
\qed

\begin{remark}
A simple variation of the proof yields the result \eqref{tailofmax1}.
\end{remark}

\subsection{Tail estimate around two insertions}

Let for $i=2,3$\footnote{The indices $2,3$ occur in the applications of this estimate in the main text.}
\begin{equation*}
W_i:=\int_{B(z_i,1)}\frac{F_i(x')}{|x'-z_i|^{\gamma\alpha_i}}M_\gamma(d^2x').
\end{equation*}
We will suppose that $|z_2| \geq 2$, $|z_3| \geq 2$ and $|z_2-z_3| \geq 3$ so that the balls $B_i=B(z_i,1)$ are well separated. We denote by  $\bar\eta_2$ and $\bar\eta_3$ the exponents occurring in the tail estimates of Lemma \ref{simpletail} applied to $W_2$ and $W_3$. Set
$$
\tilde\eta_2=\bar\eta_2\wedge  \frac{_1}{^\gamma} (Q-\alpha_3)\wedge \hf,\ \ \tilde\eta_3=\bar\eta_3\wedge  \frac{_1}{^\gamma} (Q-\alpha_2)\wedge \hf.
$$
Then we have
\begin{lemma}\label{doubletail}
For all $\beta<\bar\beta:= (\frac{2}{\gamma} (Q-\alpha_2)+\tilde\eta_2) \wedge ( \frac{2}{\gamma} (Q-\alpha_3)+\tilde\eta_3)$
\begin{equation*}
 |  \P( W_2+W_3>x ) -  \sum_{i=2}^3 |z_i|^{4\alpha_i(\alpha_i-Q)}F_i(z_i)^{\frac{2}{\gamma} (Q-\alpha_i)} \frac{\bar{R} (\alpha_i) }{ x^{\frac{2}{\gamma} (Q-\alpha_i) }}    |    \leq {C}x^{-\beta} .
\end{equation*}
\end{lemma}

\begin{remark}\label{remarkdoubletail}
The above theorem is useful when $\bar\beta >  \frac{2}{\gamma} (Q-\alpha_2) \vee  \frac{2}{\gamma} (Q-\alpha_3)$. This is the case when $\alpha_2$ and $\alpha_3$ are sufficiently close to each other.

\end{remark}

\begin{remark}\label{doubletailmoments}
The proof of Theorem \ref{doubletail} is based on the fact that the two variables $W_2$ and $W_3$ are ``nearly'' independent. Along the same lines as the proof of Lemma  \ref{doubletail}, one can in fact show that for all $p_2,p_3>0$ there exists some constant $C>0$ such that
\begin{equation*}
\E[   W_2^{p_2} W_3^{p_3}  ]  \leq C \E[   W_2^{p_2}  ]   \E[   W_3^{p_3}  ].
\end{equation*}
\end{remark}

\vspace{0.1 cm}

\proof
The strategy here is to apply the previous lemma with one insertion. 
We start with the upper bound. We have
\begin{align}\label{simpleid}
 & \P( W_2+W_3>x )\leq  \P(  W_2+W_3>x , W_2> \frac{_x}{^2})+    \P(  W_2+W_3>x , W_3> \frac{_x}{^2}).
\end{align}

The variables $W_2$ and $W_3$  are nearly independent as we now argue. We consider the circle of radius $\frac{3}{2}$ centered at $z_2$. By the Markov property of the GFF, we have the following decomposition inside $B(z_2,\frac{3}{2})$
\begin{equation*}
X(x')= \tilde{X}(x') + \mathcal{P} (X)(x')
\end{equation*}
where $\mathcal{P} (X)(x')$ is the Poisson kernel of the ball $B(z_2, \frac{3}{2})$ applied to $X$ and $\tilde{X}$ is a GFF with Dirichlet boundary conditions on $B(z_2, \frac{3}{2})$ independent of $X$ on the outside of $B(z_2, \frac{3}{2})$. On the smaller ball $B(z_2,1)$, the process $\mathcal{P} (X)(x')$ is a smooth Gaussian process hence for all $p >0$
\begin{equation*}
\E[  e^{p \sup_{  |x'-z_2| \leq 1  }  \mathcal{P} (X)(x') } ] < \infty.
\end{equation*}
We set $H=\sup_{  |x'-z_2| \leq 1  }  \mathcal{P} (X)(x') $. Of course, we have
\begin{equation*}
W_2 \leq e^{\gamma H }\tilde W_2
\end{equation*}
where $\tilde W_2$ is computed with the chaos measure of $\tilde X$.  $\tilde W_2, W_3$ have moments less than orders $\frac{2}{\gamma} (Q-\alpha_i)$ respectively \cite[Lemma A.1]{DKRV} so that
 for all $u,v>0$ and all $\epsilon'>0$ 
\begin{align*}
& \P(     W_2> u  , W_3  >  v    )  \leq \P(      e^{\gamma H } \tilde W_2 >u  ,W_3   >  v    )   \\  
& \leq \frac{1}{u^{\frac{2}{\gamma} (Q-\alpha_2)  -\epsilon'}}    \E[  \tilde W_2^{   \frac{2}{\gamma} (Q-\alpha_2)  -\epsilon' }    ]  \E[         e^{(2(Q-\alpha_2) -\gamma \epsilon'  ) H  }  1_{ W_3  >  v}     ]  \\
& \leq \frac{1}{u^{\frac{2}{\gamma} (Q-\alpha_2)  -\epsilon'}}    \E[   \tilde W_2^{   \frac{2}{\gamma} (Q-\alpha_2)  -\epsilon' }    ]  \E[         e^{p(2(Q-\alpha_2) -\gamma \epsilon'  ) H  }  ]^{1/p}  \P (W_3  >  v)^{1/q}     \\
& \leq \frac{C}{u^{\frac{2}{\gamma} (Q-\alpha_2)  -\epsilon'}  v^{  ^{\frac{1}{q}(\frac{2}{\gamma} (Q-\alpha_3)  -\epsilon')} } }   
\end{align*}
for all $p,q>1$ such that $\frac{1}{p}+\frac{1}{q}=1$.  By taking  $q$ close to $1$ we conclude  
\begin{equation}\label{nearindep}
\P(    W_2 > u  , W_3   >  v    )  \leq \frac{C}{u^{\frac{2}{\gamma} (Q-\alpha_2)  -\epsilon}  v^{{\frac{2}{\gamma} (Q-\alpha_3)  -\epsilon} } }
\end{equation}
for all $\epsilon>0$.  Therefore, exploiting \eqref{nearindep} we have for all $\epsilon>0$ 
\begin{align*}
 \P(    W_2+W_3 >x , W_2 > \frac{x}{2} ) & \leq   \P( W_2+ W_3  >x ,W_2 > \frac{x}{2}   , W_3 \leq \sqrt{x})+   \P(    W_2 > \frac{x}{2}  ,W_3 > \sqrt{x})   \\
& \leq   \P(  W_2>x- \sqrt{x} )+  \frac{C}{x^{\frac{2}{\gamma} (Q-\alpha_2)}  x^{ \frac{1}{\gamma} (Q-\alpha_3)  -\epsilon }  } .
\end{align*}
We get a similar bound by interchanging $2$ and $3$. Inserting to \eqref{simpleid} we  obtain
\begin{align*}
 & \P(   W_2+W_3  >x )  \leq \P(  W_2>x- \sqrt{x} )+  \P( W_3  >x- \sqrt{x} )  +     \frac{C}{x^{\frac{2}{\gamma} (Q-\alpha_2)}  x^{ \frac{1}{\gamma} (Q-\alpha_3)  -\epsilon }  } +  \frac{C}{x^{\frac{2}{\gamma} (Q-\alpha_3)}  x^{ \frac{1}{\gamma} (Q-\alpha_2)  -\epsilon }  }  
\end{align*}
and then we use Lemma  \ref{simpletail}  on one insertion.

Now, we proceed with the lower bound. We have, exploiting \eqref{nearindep}, that for all $\epsilon>0$
\begin{align*}
  \P(  W_2+W_3 >x ) & \geq \P(  \{ W_2 >x \} \cup\{ W_3  >x\} )  \geq \P(   W_2 >x  )+ \P(   W_3>x )   - \P(    W_2  >x  , W_3 >x )   \\ 
&  \geq \P(   W_2 >x  )+ \P(    W_3 >x )   - \frac{C}{ x^{\frac{2}{\gamma} (Q-\alpha_2)    + \frac{2}{\gamma} (Q-\alpha_3)-\epsilon  }}
\end{align*}
and then we use again  Lemma  \ref{simpletail}.

\qed

\section{Proof of Lemma    \ref{defR} on the reflection coefficient}\label{sec:proofreflection}

Recall the definitions of the reflection coefficients $\bar{R}$ and $R$ in \eqref{deffullR} and \eqref{defunitR}. For later purposes we prove a more general result than Lemma    \ref{defR}, which we state now
\begin{proposition}\label{defR1}Let $\tfrac{\gamma}{2}<\alpha_2\leq\alpha_3<Q$. Then
\begin{equation*}
\lim_{\alpha_1\searrow \alpha_3-\alpha_2}(\alpha_1-\alpha_3+\alpha_2) C_\gamma (\alpha_1, \alpha_2, \alpha_3)=
a {R}(\alpha_3)
\end{equation*}
where $a=2$ if  $\alpha_2<\alpha_3$ and $a=4$ if $\alpha_2=\alpha_3$.

\end{proposition}

\proof 
We use the formulas \eqref{firstdefstructure} and \eqref{Z1} to write
\begin{equation*}
C_\gamma (\alpha_1, \alpha_2, \alpha_3)=
2 \mu^{-s} \gamma^{-1}\Gamma(s)
 \prod_{i < j}{|z_i-z_j|^{-2\Delta_{ij}-\alpha_i \alpha_j}}\E
  \left (  \int_{\C}  F(x,{\bf z}) M_\gamma(d^2x)  \right )^{-s} 
\end{equation*}
where 
\begin{equation}\label{coulomb}
F(x,{\bf z})=\prod_{k=1}^3 \left ( \frac{ |x|_+}{|x-z_k|}  \right )^{\gamma \alpha_k} .
\end{equation}
We take  $z_1=0$ and $|z_2|, |z_3|>2$ with $|z_2-z_3| >2$. Let $\alpha_1=\alpha_3-\alpha_2+\epsilon$ (with $\epsilon>0$) so that $-s=\tfrac{2}{\gamma}(Q-\alpha_3)-\tfrac{\epsilon}{\gamma}$. Write $F$ as $F_\epsilon$ to denote its explicit dependence on $\epsilon$. We need to study the limit 
\begin{equation*}
\lim_{\epsilon\downarrow 0}\E   \left (  \int_{\C}  F_\epsilon(x,{\bf z}) M_\gamma(d^2x)  \right )^{\tfrac{2}{\gamma}(Q-\alpha_3)-\tfrac{\epsilon}{\gamma}} .
\end{equation*}
Consider first the case $\alpha_2=\alpha_3=\alpha$. Set $A_i=B(z_i,1)$ for $i=2,3$ and the complement $A_c=(A_2\cup A_3)^c$.  
Let $W_{i,\epsilon}=\int_{A_i} F_\epsilon(x,{\bf z}) M_\gamma(d^2x)$ for $i=2,3$ and $W_{c,\epsilon}=\int_{A_c} F_\epsilon(x,{\bf z}) M_\gamma(d^2x)$
so that 
$$ \int_\C F_\epsilon(x,{\bf z}) M_\gamma(d^2x)=W_{c,\epsilon}+W_{2,\epsilon}+W_{3,\epsilon}.
$$
Assume first $\frac{2}{\gamma}(Q-\alpha) \leq 1$. Then
\begin{align*}
   \E[   (  W_{2,\epsilon}+W_{3,\epsilon} )  ^{\frac{2}{\gamma}(Q-\alpha)  -\frac{\epsilon}{\gamma}}]  &   \leq \E[   (  W_{c,\epsilon}+W_{2,\epsilon}+W_{3,\epsilon}  )^{\frac{2}{\gamma}(Q-\alpha)  -\frac{\epsilon}{\gamma}}]    \\
   & \leq   \E[     W_{c,\epsilon}^{\frac{2}{\gamma}(Q-\alpha)  -\frac{\epsilon}{\gamma}}]+ \E[( W_{2,\epsilon}+W_{3,\epsilon})^{\frac{2}{\gamma}(Q-\alpha)  -\frac{\epsilon}{\gamma}}].
\end{align*}
In order to treat the second expectation above, we apply the double tail estimate Lemma \ref{doubletail} (with $F_2(x)=\frac{|x|^{\gamma(\alpha_2+\alpha_3)}}{|x-z_3|^{\gamma\alpha_3}}$ and $F_3(x)=\frac{|x|^{\gamma(\alpha_2+\alpha_3)}}{|x-z_2|^{\gamma\alpha_2}}$) to get
\begin{align*}
   \P &(  W_{2,\epsilon}+W_{3,\epsilon}  >x)\\
   =&\Big(|z_2|^{4\alpha_2(\alpha_2-Q)}F_2(z_2)^{\frac{2}{\gamma}(Q-\alpha_2)}\bar{R}(\alpha_2)x^{-\frac{2}{\gamma}(Q-\alpha_2)} +|z_3|^{4\alpha_3(\alpha_3-Q)}F_3(z_3)^{\frac{2}{\gamma}(Q-\alpha_3)}\bar{R}(\alpha_3)x^{-\frac{2}{\gamma}(Q-\alpha_3)}  \Big)(1+\caO(x^{-\eta}))   \\
   =&2|z_2-z_3|^{-2 \alpha(Q-\alpha)}  \bar{R}(\alpha) x^{-\frac{2}{\gamma}(Q-\alpha)}(1+\caO(x^{-\eta}))
\end{align*}
for $\eta>0$, uniformly in $\epsilon$.

Since $F_\epsilon(x,{\bf z})\leq C(1+|x|^{-\gamma\epsilon}\mathbf{1}_{\{|x|\leq1\}})$ on $A_c$ we deduce that
 $\epsilon \E[     W_{c,\epsilon}^{\frac{2}{\gamma}(Q-\alpha)  -\frac{\epsilon}{\gamma}}]$ converges to $0$ as $\epsilon\to 0$. 
 Indeed, by Jensen,
\begin{multline*}
 \E[     W_{c,\epsilon}^{\frac{2}{\gamma}(Q-\alpha)  -\frac{\epsilon}{\gamma}}]\leq  \E[     W_{c,\epsilon} ]^{\frac{2}{\gamma}(Q-\alpha)  -\frac{\epsilon}{\gamma}}  =\Big(\int_{A_c}F_\epsilon(x,{\bf z})|x|_+^{-4}\,d^2x\Big)^{\frac{2}{\gamma}(Q-\alpha)  -\frac{\epsilon}{\gamma}}\\
 \leq C\Big(\int_{A_c}(1+|x|^{-\gamma\epsilon}\mathbf{1}_{\{|x|\leq1\}})|x|_+^{-4}\,d^2x\Big)^{\frac{2}{\gamma}(Q-\alpha)  -\frac{\epsilon}{\gamma}}
 \end{multline*}
 and this quantity is obviously bounded in $\epsilon$. We deduce
\begin{equation*}
\lim_{\epsilon\to 0}\epsilon  \E[  
( \int_\C F_\epsilon(x,{\bf z}) M_\gamma(d^2x))^{\frac{2}{\gamma}(Q-\alpha)  -\frac{\epsilon}{\gamma}}] 
=4(Q-\alpha)  |z_2-z_3|^{-2 \alpha(Q-\alpha)}  \bar{R}(\alpha)  
\end{equation*}
and then (note: we know that the $z_i$-dependence has to drop out!)
\begin{equation*}
\lim_{\epsilon\to 0}\epsilon C_\gamma (\epsilon, \alpha,\alpha) = \mu^{\frac{2}{\gamma} (Q-\alpha)}    \frac{_{8(Q-\alpha)}}{^\gamma}  \Gamma (-\frac{_{2(Q-\alpha)}}{^\gamma} ) \bar{R}(\alpha)    =   4 R(\alpha).  \\
\end{equation*}
If $\frac{2}{\gamma}(Q-\alpha) >1 $ we have by triangle inequality and $\epsilon$ small enough so that $p=\frac{2}{\gamma}(Q-\alpha)  -\frac{\epsilon}{\gamma}>1$
\begin{align*}
&   [ \E   (  W_{2,\epsilon}+W_{3,\epsilon}   )  ^p] ^{1/p}
   \leq[ \E  ( W_{c,\epsilon}+W_{2,\epsilon}+W_{3,\epsilon}   )  ^p] ^{1/p}
   \leq  [\E   ( W_{c,\epsilon}  )  ^p] ^{1/p}+ [\E   (  W_{2,\epsilon}+W_{3,\epsilon}   )  ^p] ^{1/p}
\end{align*}
and
we can conclude similarly as the previous case.

The case $\alpha_2<\alpha_3$ is similar: we use the tail estimate Lemma \ref{simpletail} around the $\alpha_3$ insertion. The difference of a factor of two results from the sum over two insertions in the double tail estimate.\qed
\begin{remark} \label{unitvolrel}For the unit volume quantities defined in \eqref{unitvolumethreepoint}   and \eqref{defunitR}   we get
\begin{equation*}
\lim_{\alpha_1\downarrow \alpha_3-\alpha_2}(\alpha_1-\alpha_3+\alpha_2) \bar C_\gamma (\alpha_1, \alpha_2, \alpha_3)=
a \tfrac{2(Q-\alpha_3)}{\gamma}{\bar R}(\alpha_3).
\end{equation*}
\end{remark}

\section{The BPZ equations and algebraic relations}\label{bpzeq}

This section is devoted to the study of the small $z$ asymptotics of the  four point functions $\caT_{-\frac{\gamma}{2}}$ and $\caT_{-\frac{2}{\gamma}}$  leading to the proof of \eqref{Fundrelation3} and \eqref{sketch:BPZ1}. The proof of the latter is the technical core of the paper and the key input in the probabilistic identification of the reflection coefficient.

\subsection{Fusion without reflection} As mentioned in Section \ref{BPZ equations} the relation \eqref{Fundrelation3} was proven in \cite[Theorem 2.3]{KRV} with the assumption $\frac{1}{\gamma}+\gamma<\alpha_1+\frac{\gamma}{2}<Q$ or in other words $\frac{\gamma}{2}+\frac{1}{\gamma}<\alpha_1<\frac{2}{\gamma}$. This interval is non empty if and only if $\gamma^2<2$. In this section we will remove this constraint. The reason for the restriction $\frac{\gamma}{2}+\frac{1}{\gamma}<\alpha_1$ was the following. In order to prove \eqref{Fundrelation3}, one must perform the asymptotic expansion of $\caT_{-\frac{\gamma}{2}}(z)$ around $z \to 0$ \eqref{Fundrelation300} as explained in section \ref{BPZ equations}.  In the case $\frac{1}{\gamma}+\frac{\gamma}{2}<\alpha_1$,  the exponent $2(1-c)$ which is equal to $\gamma(Q-\alpha_1)$ is strictly less than $1$ hence there are no polynomial terms in $z$ and $\bar z$ in the expansion \eqref{Fundrelation300} to that order (such terms are present in the small $z$ expansion of $|F_-(z)|^2$). In the case $\alpha_1<\frac{1}{\gamma}+\frac{\gamma}{2}$, the asymptotic expansion of $\caT_{-\frac{\gamma}{2}}(z)$ around $0$ is more involved. Nonetheless, we prove here:

\begin{theorem}\label{theo4point}
 We assume the Seiberg bounds for $(-\frac{\gamma}{2},\alpha_1,\alpha_2,\alpha_3)$, i.e.  $\sum_{k=1}^3 \alpha_k>2Q+\frac{\gamma}{2}$ and $\alpha_k<Q$ for all $k$.  
If  $\frac{\gamma}{2}< \alpha_1< \frac{2}{\gamma}$ then  
 \begin{equation}\label{theo4pointexpressionBPZ1}
\mathcal{T}_{-\frac{\gamma}{2}}(z)= C_\gamma(\alpha_1-\tfrac{\gamma}{2}, \alpha_2,\alpha_3)  |F_{-}(z)|^2  - \mu \frac{\pi}{  l(-\frac{\gamma^2}{4}) l(\frac{\gamma \alpha_1}{2})  l(2+\frac{\gamma^2}{4}- \frac{\gamma \alpha_1}{2}) } C_\gamma(\alpha_1+\tfrac{\gamma}{2}, \alpha_2,\alpha_3)  |F_{+}(z)|^2
\end{equation}
and the relation \eqref{3pointconstanteqintro} holds.
\end{theorem}

\proof
 Let first $\gamma^2<2$.  \eqref{theo4pointexpressionBPZ1} was proven  in \cite[Theorem 2.3]{KRV} in the case $\frac{\gamma}{2}+\frac{1}{\gamma}<\alpha_1<\frac{2}{\gamma}$. 
 This result extends to the interval  $\frac{\gamma}{2}<\alpha_1<\frac{2}{\gamma}$  by analyticity. Indeed, for fixed $\gamma\in (0,\sqrt{2})$, the interval  $\frac{\gamma}{2}+\frac{1}{\gamma}<\alpha_1<\frac{2}{\gamma}$ is non empty. Furthermore, by Theorem \ref{analytic_hyperbolic} both sides of eq. \eqref{theo4pointexpressionBPZ1} are analytic in $\alpha_1$ (with other parameters fixed) in a neighborhood of the interval $\frac{\gamma}{2}< \alpha_1< \frac{2}{\gamma}$ seen as a subset of $\C$. Uniqueness of analytic continuation thus establishes \eqref{theo4pointexpressionBPZ1} for $\gamma^2<2$. $\gamma^2=2$ is obtained by continuity in $\gamma$ (see Remark \ref{analingamma} on this).

Let now $\gamma^2>2$ 
and $\frac{\gamma}{2}<\alpha_1<\frac{2}{\gamma}$: 
The proof of \eqref{theo4pointexpressionBPZ1} follows from the study of the function $\mathcal{T}_{-\frac{\gamma}{2}}(z)$ as $z$ tends to $0$. More precisely, by the discussion in subsection \ref{BPZ equations}, it suffices to show that one has the following expansion as $z$ goes to $0$
\begin{align}\label{expansionmathcalT}
\mathcal{T}_{-\frac{\gamma}{2}}(z)=C_\gamma(\alpha_1-\tfrac{\gamma}{2}, \alpha_2,\alpha_3)+ C z+ \bar{C} \bar{z}-  \frac{\mu\pi C_\gamma(\alpha_1+\frac{\gamma}{2}, \alpha_2,\alpha_3) }{  l(-\frac{\gamma^2}{4}) l(\frac{\gamma \alpha_1}{2})  l(2+\frac{\gamma^2}{4}- \frac{\gamma \alpha_1}{2}) } |z|^{\gamma (Q-\alpha_1)}+ o( |z|^{\gamma (Q-\alpha_1)})
\end{align}
where $C$ is some constant. Thus  by \eqref{Tdefi0} we need to study the function \eqref{defcalS} with $\alpha_0=-\frac{\gamma}{2}$. To streamline notation let us set 
\begin{equation}\label{defineK}K(z,x)=\frac{|x-z|^{\frac{\gamma^2}{2}}|x|_+^{\gamma (\sum_{k=1}^3\alpha_k-\frac{\gamma}{2})}}{|x|^{\gamma \alpha_1}|x-1|^{\gamma \alpha_2}}
\end{equation}
and for any Borel set $B\subset \C$
\begin{equation}\label{def:K}
\mathcal{K}_B(z):=\int_{B}K(z,x)\,M_\gamma(d^2x).
\end{equation}
Then  $\mathcal{S}_{-\frac{\gamma}{2}}(z)=\caK_\C(z)$ where $\mathcal{S}_{-\frac{\gamma}{2}}$ was defined in \eqref{defcalS}. We will also write $\mathcal{K}(z)$ for $\mathcal{K}_\C(z)$. We set $p:=\frac{1}{\gamma}(\sum_{k=1}^3 \alpha_k-\frac{\gamma}{2}-2Q)$. Then
\begin{align}\label{Tdefi000}
  \mathcal{T}_{-\frac{\gamma}{2}}(z) &= 
  2 \mu^{-p}  \gamma^{-1}\Gamma(p)  \E[\mathcal{K}(z)^{-p}] .
 \end{align} 
A Taylor expansion yields the relation
\begin{equation}\label{expansionmomentp}
   \E[\mathcal{K}(z)^{-p}]  =\E[\mathcal{K}(0)^{-p}]   +z\partial_z\E[\mathcal{K}(z)^{-p}]_{|z=0} +\bar{z}\partial_{\bar{z}}\E[\mathcal{K}(z)^{-p}]_{|\bar{z}=0} +\mathcal{R}(z)  
  \end{equation}
where $\mathcal{R}(z)$ is a remainder term whose expression appears below \eqref{defreminder}  (not to be confused with the reflection coefficient!).
First, notice that the term $\partial_z\E[\mathcal{K}(z)^{-p}]_{|z=0}$ is well defined. Indeed we have
\begin{align*}
\partial_z\E[\mathcal{K}(z)^{-p}]_{|z=0} =-p\frac{\gamma^2}{4}\E\Big[\int_\C\frac{1}{x}K(0,x)M_\gamma(d^2x)\mathcal{K}(0)^{-p-1}\Big].
\end{align*}
Split the   $x$-integral over $\C$ in two parts, over ${B_{1/2}} $ and over ${B_{1/2}^c} $ where in this section we use the notation $B_{r}=B(0,r)$. Then
\begin{align*}
\E \Big[\int_{B_{1/2}^c}\frac{1}{|x|}K(0,x)M_\gamma(d^2x)\mathcal{K}(0)^{-p-1}\Big]\leq 2\E\, \mathcal{K}(0)^{-p}<\infty
\end{align*}
 as GMC measures possess negative moments of all orders (see Appendix, subsection \ref{chaosestimates}). For the integral over  $B_\hf$ we use the Cameron-Martin theorem (Corollary \ref{coro:Girsanov} in the Appendix) to get
\begin{align}\nonumber
&\E\Big[\int_{B_{1/2}}\frac{1}{|x|}K(0,x)M_\gamma(d^2x)\mathcal{K}(0)^{-p-1}\Big]\leq C\E \Big[ \int_{B_{1/2}}|x|^{-1-\gamma\alpha_1+\tfrac{\gamma^2}{2}}M_\gamma(d^2x)\mathcal{K}(0)^{-p-1}\Big]\\&= C\int_{B_{1/2}}|x|^{-1-\gamma\alpha_1+\tfrac{\gamma^2}{2}}\E\Big[\big(\int_\C K(0,u)e^{\gamma^2G(x,u)}M_\gamma(d^2u)\big)^{-p-1}\Big]\,d^2x.\label{bhalfboun}
\end{align}
To bound the last expectation we note that  the integrand in the $u$-integral is  bounded away from $0$  for $x\in B_\hf$ and   $u\in B(3,1)$. This ball is far away from the singularities, hence on $u\in B(3,1)$ the kernel $K(0,u)e^{\gamma^2G(x,u)}$ is bounded from below  away from $0$ . Thus  
\begin{align*}
\E\Big[\big(\int_\C K(0,u)e^{\gamma^2G(x,u)}M_\gamma(d^2u)\big)^{-p-1}\Big]\leq C\E[M_\gamma(B(3,1))^{-p-1}]<\infty
\end{align*}
as the measure $M_\gamma$ possesses moments of negative order (see Appendix, subsection \ref{chaosestimates}). The final integral in \eqref{bhalfboun} converges  as the constraint $\alpha_1<\frac{2}{\gamma}$ guarantees that $1+\gamma\alpha_1-\tfrac{\gamma^2}{2}<3-\tfrac{\gamma^2}{2}<2$ since $\gamma^2>{2}$. The same argument shows that $\bar{z}\partial_{\bar{z}}\E[(\mathcal{K}_0(z))^{-p}]_{|\bar{z}=0}$ is well defined.

It remains to investigate the remainder $\mathcal{R}(z)$ which by the Taylor integral formula is given by
\begin{equation}\label{defreminder}
\mathcal{R}(z)=
\int_0^1 (1-t) \Big(z^2\partial^2_{z}\E[\mathcal{K}(tz)^{-p}]+2z\bar{z}\partial_{z}\partial_{\bar{z}}\E[\mathcal{K}(tz)^{-p}]  +\bar{z}^2\partial^2_{\bar{z}}\E[\mathcal{K}(tz)^{-p}] \Big)\,dt\equiv \mathcal{R}_1(z)+\mathcal{R}_2(z)+\mathcal{R}_3(z).
\end{equation}

Expression \eqref{defreminder} consists of  the three terms   $\mathcal{R}_i(z)$, $i=1,2,3$. The first term in \eqref{defreminder} is given by 
$$\mathcal{R}_1(z)=\int_0^1 (1-t)z^2\partial^2_{z}\E[\mathcal{K}(tz)^{-p}] dt=r_1(z)+p_1(z)
$$
with 
\begin{align*}
r_1(z):=&-p\frac{_{\gamma^2}}{^4}(\frac{_{\gamma^2}}{^4}-1)z^2\int_0^1(1-t)\E\Big[\int_\C\frac{1}{(x-tz)^2}K(tz,x)M_\gamma(d^2x)(\mathcal{K}(tz))^{-p-1}\Big]\,dt\\
p_1(z):=&p(p+1)\frac{_{\gamma^4}}{^{16}}z^2\int_0^1(1-t)\E\Big[\Big(\int_\C\frac{1}{(x-tz)}K(tz,x)M_\gamma(d^2x)\Big)^2(\mathcal{K}(tz))^{-p-2}\Big]\,dt.
\end{align*}
The terms $p_2,r_2,p_3,r_3$ are defined similarly with respect to $\mathcal{R}_i(z)$, $i=2,3$.  We get the expansion of $\mathcal{R}(z)$ around $z=0$ thanks to the following two lemmas (whose proof is postponed right after):
\begin{lemma}\label{lemma1int}
The following holds
\begin{equation*}
r_1(z)+r_2(z)+r_3(z)=
- pA  |z|^{ \gamma(Q-\alpha_1)}+o(|z|^{ \gamma(Q-\alpha_1)})
\end{equation*}
where \begin{align*}
A=
\frac{\pi}{  l(-\frac{\gamma^2}{4}) l(\frac{\gamma \alpha_1}{2})  l(2+\frac{\gamma^2}{4}- \frac{\gamma \alpha_1}{2}) } \E\Big[  \Big(\int_\C K(0,u)e^{\gamma^2G(0,u)}M_\gamma(d^2u)\Big)^{-p-1}\Big] . 
\end{align*}
\end{lemma}

 \begin{lemma}\label{lemma2int}
 The following holds 
 \begin{equation*}
 p_1(z)+p_2(z)+p_3(z)= o(|z|^{\gamma(Q-\alpha_1)}).
 \end{equation*}
 \end{lemma}

From these two lemmas, one can deduce \eqref{expansionmathcalT}. Indeed, since 
\begin{equation*}
 2 \mu^{-p}  \gamma^{-1} p \Gamma(p)  A= \frac{\pi  \mu C_\gamma(\alpha_1+\frac{\gamma}{2}, \alpha_2,\alpha_3) }{   l(-\frac{\gamma^2}{4}) l(\frac{\gamma \alpha_1}{2})  l(2+\frac{\gamma^2}{4}- \frac{\gamma \alpha_1}{2})  } 
 \end{equation*} 
 the above two lemmas imply that
\begin{equation*}
2 \mu^{-p}  \gamma^{-1}  \Gamma(p) \mathcal{R}(z)=  -   \frac{\pi  \mu C_\gamma(\alpha_1+\frac{\gamma}{2}, \alpha_2,\alpha_3) }{   l(-\frac{\gamma^2}{4}) l(\frac{\gamma \alpha_1}{2})  l(2+\frac{\gamma^2}{4}- \frac{\gamma \alpha_1}{2})  }     |z|^{ \gamma(Q-\alpha_1)} + o(|z|^{\gamma(Q-\alpha_1)})
\end{equation*}
which yields \eqref{expansionmathcalT} thanks to the fact that $\mathcal{T}_{-\frac{\gamma}{2}}(0) =C_\gamma(\alpha_1-\tfrac{\gamma}{2}, \alpha_2,\alpha_3)$  and using \eqref{Tdefi000} and \eqref{expansionmomentp}. \qed

\medskip
\noindent
\emph{Proof of Lemma \ref{lemma1int}}.
We first study the $r_1$ term. The term  $r_1$ is analysed in the same way as a similar term in the proof of \cite[Lemma 4.5]{KRV} so we will be brief. First as above we want to restrict the $x$-integral to the ball $B_{1/2}$. The integral over $B_{1/2}^c $ produces an $\caO(z^2)$ contribution: indeed for $z$ small enough (say $|z|\leq 1/4$)
\begin{align*}
\E \Big[\int_{B_{1/2}^c}\frac{1}{|x-tz|^2}K(tz,x)M_\gamma(d^2x)\mathcal{K}(tz)^{-p-1}\Big]\leq 16\E\, \mathcal{K}(tz)^{-p}<\infty
\end{align*}
uniformly in $|z|\leq 1/4$ and $t\in [0,1]$. This can be seen by restricting  the integral in $\mathcal{K}(tz)$ to a fixed ball away from singularities and then use the fact that  GMC measures possess negative moments of all orders (see Appendix, subsection \ref{chaosestimates}).

Now we can focus on the $x$-integral over the ball $B_{1/2}$.  By  the Cameron-Martin theorem and  the change of variables $x\to ytz$ we get 
\begin{align*}
r_1(z)=-p\frac{_{\gamma^2}}{^4}(\frac{_{\gamma^2}}{^4}-1)&|z|^{2+\frac{\gamma^2}{2}-\gamma\alpha_1}\int_0^1(1-t)t^{\frac{\gamma^2}{2}-\gamma\alpha_1}\int_{B_{\frac{1}{2t|z|}}}\frac{|y-1|^{\frac{\gamma^2}{2}}}{(y-1)^2|y|^{\gamma\alpha_1}|yzt-1|^{\gamma\alpha_2}}\\ &\times\E\Big[  \Big(\int_\C K(tz,u)e^{\gamma^2G(tyz,u)})M_\gamma(d^2u)\Big)^{-p-1}\Big]\,dtd^2y+\caO(z^2).
\end{align*}
Dominated convergence theorem then implies 
\begin{align*}
r_1(z)=-p\frac{_{\gamma^2}}{^4}(\frac{_{\gamma^2}}{^4}-1)&|z|^{ \gamma(Q-\alpha_1)}\int_0^1(1-t)t^{\frac{\gamma^2}{2}-\gamma\alpha_1}dt\int_{\C}\frac{|y-1|^{\frac{\gamma^2}{2}}}{(y-1)^2|y|^{\gamma\alpha_1}}\,d^2y\\&\times \E\Big[  \Big(\int_\C K(0,u)e^{\gamma^2G(0,u)}M_\gamma(d^2u)\Big)^{-p-1}\Big]+o(|z|^{ \gamma(Q-\alpha_1)}).
\end{align*}
Using eq. \eqref{gammaz} in the Appendix to the  $y$ integral  finally yields
\begin{align}\label{equivR1}
r_1(z)=
- \frac{p}{4}   \frac{\frac{\gamma^2}{2}-{\gamma \alpha_1}}{ \frac{\gamma^2}{2}-\gamma \alpha_1+1  } A  |z|^{ \gamma(Q-\alpha_1)}+o(|z|^{ \gamma(Q-\alpha_1)})
\end{align}
with
\begin{align*}
A=
\frac{\pi}{  l(-\frac{\gamma^2}{4}) l(\frac{\gamma \alpha_1}{2})  l(2+\frac{\gamma^2}{4}- \frac{\gamma \alpha_1}{2}) } \E\Big[  \Big(\int_\C K(0,u)e^{\gamma^2G(0,u)}M_\gamma(d^2u)\Big)^{-p-1}\Big]  .
\end{align*}
Since $r_3$ in \eqref{defreminder} equals $\bar r_1$ it is also  given by \eqref{equivR1}. Finally $r_2$ yields 
\begin{align*}
r_2(z)=
- \frac{p}{2}   \frac{\frac{\gamma^2}{2}-{\gamma \alpha_1}+2}{ \frac{\gamma^2}{2}-\gamma \alpha_1+1  } A  |z|^{ \gamma(Q-\alpha_1)}+o(|z|^{ \gamma(Q-\alpha_1)}).
\end{align*}
Altogether we then get
\begin{align*}
r_1(z)+r_2(z)+r_3(z)=
- pA  |z|^{ \gamma(Q-\alpha_1)}+o(|z|^{ \gamma(Q-\alpha_1)}).
\end{align*}

\qed

 \medskip

\noindent
\emph{Proof of Lemma \ref{lemma2int}}.
We will prove that
 $p_1(z)$ is a $o(|z|^{\gamma(Q-\alpha_1)})$, the  argument for $p_2, p_3$ is similar. We will bound the expectation occurring  in  $p_1(z)$ so let us denote
 \begin{align}\label{idefin}
I(B,tz):=\E\Big[\Big(\int_{B}\frac{1}{|x-tz|}K(tz,x)M_\gamma(d^2x)\Big)^2(\int_{\C}K(tz,u)M_\gamma(d^2u))^{-p-2}\Big].
\end{align} 
We shall  prove
  \begin{align}\label{fundest}
I(\C,tz)\leq C|tz|^{\gamma(Q-\alpha_1)-2+\eta}
\end{align}
with $\eta>0$ which proves our claim since \eqref{fundest} implies 
\begin{equation*}
|p_1(z)|\leq p(p+1)\frac{_{\gamma^4}}{^{16}}|z|^2 \int_0^1(1-t)  I(\C,tz) \, dt \leq C  |z|^{\gamma(Q-\alpha_1)+\eta}  \int_0^1(1-t)  t^{  \gamma(Q-\alpha_1)-2+\eta } \, dt    \leq C |z|^{\gamma(Q-\alpha_1)+\eta}  
\end{equation*}
where we used the fact that the  $t$ integral converges at $0$ as $\gamma(Q-\alpha_1)-2=\frac{\gamma^2}{2}-2>-1$. 

We can now put $t=1$ and we will bound $I(\C,z)$ for $z$ small.  For  $z$ small enough $\frac{1}{|x-z|}$ is bounded in $B_\hf^c$ and we have $I(B_\hf^c,z) \leq C \E[ \mathcal{K}(z)^{-p}  ] \leq C$.  Since $I(\C,z)\leq 2(I(B_\hf,z)+I(B^c_\hf,z))$  it suffices to bound $I(B_\hf,z)$. 
 
Next we bound $I(A,z)$ where $A$ is   the annulus centered at origin with radii $L|z|$ and $\frac{1}{2}$ and $L>1$ will be chosen later. First, we use Jensen's inequality in the normalized measure $1_A(x)K(z,x)M_\gamma(d^2x)$ to get
\begin{align*}
I(A,z)\leq 
\E\Big[\int_A\frac{1}{|x-z|^2}K(z,x) M_\gamma(d^2x) \, \mathcal{K}_A(z)^{-p-1}\Big].
\end{align*}
Up to an additive independent Gaussian random variable, the restriction of $X$ to $ B_{\frac{1}{2}}$ satisfies a continuum version of the FKG inequality  (see section \ref{chaosestimates} in the appendix)  and therefore
\begin{align*}
&\E\Big[\int_A\frac{1}{|x-z|^2}K(z,x) M_\gamma(d^2x) \, \mathcal{K}_A(z)^{-p-1}\Big]\leq C \E\Big[\int_A\frac{1}{|x-z|^2}K(z,x) M_\gamma(d^2x)\Big] \,\E [\mathcal{K}_A(z)^{-p-1}]\
\\
\leq &C\int_{A}\frac{|x-z|^{\frac{\gamma^2}{2}-2}}{|x|^{\gamma\alpha_1}}  d^2x\leq C|z|^{\frac{\gamma^2}{2}-\gamma\alpha_1}\int_{|y|>L}\frac{|y-1|^{\frac{\gamma^2}{2}-2}}{|y|^{\gamma\alpha_1}} d^2y\leq C |z|^{-2}|z|^{\gamma(Q-\alpha_1)}L^{-\gamma(\alpha_1-\frac{\gamma}{2})} 
\end{align*}
where the last integral was convergent due to  $\alpha_1>\frac{\gamma}{2}$. This fits to \eqref{fundest} provided we take $L=|z|^{-\delta}$ with $\delta>0$.

We are left with estimating $I(B_{L|z|},z)$. Let us first consider the part not too close to the singularity at $z$: set  $S:=B_{L|z|}\setminus B(z, |z|^{1+\epsilon})$ for some $\epsilon>0$, to be fixed later.   We have
\begin{align*}
\E\Big[\Big(\int_S\frac{1}{|x-z|}K(z,x)M_\gamma(d^2x)\Big)^2\mathcal{K}(z)^{-p-2}\Big] \leq |z|^{-2-2\epsilon}\E[\mathcal{K}_S(z)^2\mathcal{K}(z)^{-p-2}].
\end{align*}
Then we get, for  $r\in (0,2)$ using the fact that $\mathcal{K}_S(z)\leq \mathcal{K}(z)$
\begin{align*}
\E[\mathcal{K}_S(z)^2\mathcal{K}(z)^{-p-2}]\leq \E[\mathcal{K}_S(z)^r(\mathcal{K}(z))^{-p-r}]\leq C(\E\mathcal{K}_S(z)^{qr})^{1/q}
\end{align*}
where in the second step we used  H\"older inequality and bounded the negative  GMC moment again by a constant. Finally, since $|x-z|\leq  2|Lz|$ on $S$ we get 
\begin{align*}
[\E(\mathcal{K}_S(z)^{qr})]^{1/q}\leq C|Lz|^{\frac{\gamma^2}{2}r}[\E[(\int_{B_{L|z|}}|x|^{-\gamma\alpha_1}M_\gamma(d^2x))^{qr}]]^{1/q}\leq C|Lz|^{\gamma(Q-\alpha_1+\frac{\gamma}{2})r-\hf\gamma^2qr^2}
\end{align*}
where the last estimate comes from the estimates of subsection \ref{chaosestimates} in the Appendix. Here we need to assume  that $rq<\frac{4}{\gamma^2}\wedge \frac{2}{\gamma}(Q-\alpha_1)$.
Notice that since we assume $\frac{\gamma}{2}<\alpha_1$ then $\frac{4}{\gamma^2}>\frac{2}{\gamma}(Q-\alpha_1)$ so that given $q$ we need to have $0<r< \frac{2}{\gamma q}(Q-\alpha_1)$. The optimal choice for $r$ is  $r^\star=\frac{\frac{\gamma}{2}+Q-\alpha_1}{\gamma q}$ (this is less than $ \frac{2}{\gamma q}(Q-\alpha_1)$ for $\alpha_1<\frac{2}{\gamma}$), in which case 
$$ \E[\mathcal{K}_S(z)^{rq}]^{1/q}\leq C|Lz|^{\frac{1}{2q}(\frac{\gamma}{2}+Q-\alpha_1)^2}.
$$
Gathering everything we conclude 
\begin{align*}
\E\Big[\Big(\int_S\frac{1}{|x-z|}K(z,x)M_\gamma(d^2x)\Big)^2\mathcal{K}(z)^{-p-2}\Big] \leq  &CL^{\frac{1}{2q}(\frac{\gamma}{2}+Q-\alpha_1)^2} |z|^{-2-2\epsilon+\frac{1}{2q}(\frac{\gamma}{2}+Q-\alpha_1)^2}.
\end{align*}
We can now fix $\delta, q,\epsilon$. First notice that $\frac{1}{2}(\frac{\gamma}{2}+Q-\alpha_1)^2-\gamma(Q-\alpha_1)=\frac{1}{2}(Q-\alpha_1-\frac{\gamma}{2})^2>0$. Hence  choosing  $q$ sufficiently close to $1$ and  then $\epsilon<\epsilon(q)$ and finally $\delta<\delta(\epsilon)$, $I(S,z)$ can be bounded by \eqref{fundest}.

\medskip

We are thus left with proving $I(B,z)\leq C|z|^{\gamma(Q-\alpha_1)-2+\eta} $ where $B:=B(z, |z|^{1+\epsilon})$. An application of the Cameron-Martin theorem (in fact we use recursively Corollary \ref{coro:Girsanov} in the Appendix) gives
\begin{align}\label{fuck}
I(B,z)=& \int_{B^2}\frac{K(z,x)K(z,x')e^{\gamma^2G(x,x')}}{(x-z)(x'-z)}\E\Big[\Big(\int_\C K(z,u)e^{\gamma^2G(x,u)+\gamma^2G(x',u)}M_\gamma(d^2u)\Big)^{-p-2}\Big]\,d^2xd^2x'\nonumber\\
\leq& C \int_{B^2}\frac{|x-z|^{\frac{\gamma^2}{2}-1}|x'-z|^{\frac{\gamma^2}{2}-1} }{|x|^{\gamma\alpha_1}|x'|^{\gamma\alpha_1}|x-x'|^{\gamma^2} } \E\Big[\Big(\int_{B_\hf} \frac{|u-z|^{\frac{\gamma^2}{2}}}{|u|^{\gamma\alpha_1}|u-x|^{\gamma^2}|u-x'|^{\gamma^2}}  M_\gamma(d^2u)\Big)^{-p-2}\Big]\,d^2xd^2x'
\end{align}
where in the upper bound we restricted the $u$ integral to $B_\hf$. By a change of variables $x=zy,x'=zy'$ this becomes
\begin{align}\label{fuck1}
I(B,z)\leq& 
C|z|^{2-2\gamma\alpha_1}\int_{B(1,|z|^{\epsilon})^2}\frac{|y-1|^{\frac{\gamma^2}{2}-1}|y'-1|^{\frac{\gamma^2}{2}-1} }{|y-y'|^{\gamma^2} }A(y,y',z)\,d^2yd^2y'
\end{align}
with
\begin{align*}
A(y,y',z)=\E\Big[\Big(\int_{B_\hf} \frac{|u-z|^{\frac{\gamma^2}{2}}}{|u|^{\gamma\alpha_1}|u-yz|^{\gamma^2}|u-y'z|^{\gamma^2}}  M_\gamma(d^2u)\Big)^{-p-2}\Big].
\end{align*}
Note that the only potential divergence in the $y,y'$ integral is at $y=y'$ since $\gamma^2>2$. Hence we need to study how $A(y,y',z)$ vanishes on the diagonal. The behaviour of $A(y,y',z)$ as $y\to y'$ is controlled  by the fusion rules (see \cite{KRV}). In the case at hand we have four insertions, located at $0,zy,zy',z$ that are all close to each other as $z\to 0$. Fusion estimates have been proven in \cite{KRV} in the case of three insertions. A simple adaptation of that proof to the case of  4 insertions is stated in Lemma \ref{fusion4} of the Appendix. The estimate for $A(y,y',z)$ depends on the relative positions of the four insertions. In our case we have  $|zy-z|\vee|zy'-z|\vee|zy-zy'|\ll |z|\wedge|zy|\wedge|zy'|$. This means that the insertions $zy,zy',z$ will merge together way before merging with $0$. We will partition the integration region in \eqref{fuck1} acording  to the relative positions of these three points or equivalently the relative positions of $y,y',1$. By symmetry in $y,y'$ we have then three integration regions in \eqref{fuck1}  to consider 
\vskip 2mm
\noindent $\bullet$  Let  $A_1:=\{|y-1|\leq |y'-1| \leq|y-y'|\}$. Then on $B(1,|z|^{\epsilon})^2\cap A_1$ we have by Lemma \ref{fusion4} (applied with $y_1=z,y_2=zy',y_3=zy,y_4=0$)
 $$A(y,y',z)\leq C|1-y'|^{\frac{1}{2}(\frac{3\gamma}{2}-Q)^2} 
  |z|^{\frac{1}{2}(\frac{3\gamma}{2}+\alpha_1-Q)^2  }.$$
 Since $2-2\gamma\alpha_1+\frac{1}{2}(\frac{3\gamma}{2}+\alpha_1-Q)^2 =-2+\gamma(Q-\alpha_1)+\frac{1}{2}(\alpha_1-\frac{2}{\gamma})^2$ we get
 \begin{align}\label{ibo}
 I(B(1,|z|^\epsilon)\cap A_1,z)\leq &C|z|^{-2+\gamma(Q-\alpha_1)+\frac{1}{2}(\alpha_1-\frac{2}{\gamma})^2} \int_{B(1,
 |z|^\epsilon)^2}\frac{1}{
 |y-y'|^{2-\frac{1}{2}(\frac{3\gamma}{2}-Q)^2} } \,d^2yd^2y' .
 \end{align}
The integral is convergent if $(\frac{3\gamma}{2}-Q)^2>0$ which is the case if $\gamma^2\neq 2$. 
 \vskip 2mm

 \noindent $\bullet$  Let  $A_2:=\{|y-1| \leq|y-y'|\leq |y'-1|\}$. Then on $B(1,|z|^\epsilon)^2\cap A_2$ we have  by Lemma \ref{fusion4}

 $$A(y,y',z)\leq C|y-y'|^{\frac{1}{2}(\frac{3\gamma}{2}-Q)^2} 
  |z|^{\frac{1}{2}(\frac{3\gamma}{2}+\alpha_1-Q)^2  }.$$
   Hence we also end up with the bound \eqref{ibo} with $A_1$ replaced by $A_2$ (since $\frac{\gamma^2}{2}-1>0$).
 
 \vskip 2mm

  \noindent $\bullet$  Let  $A_3:=\{ |y-y'|\leq |y-1|\leq |y'-1|\}$. Then on $B(1,|z|^\epsilon)^2\cap A_3$ we have  by Lemma \ref{fusion4} (applied with $y_1=zy,y_2=z,y_3=zy',y_4=0$)

  $$A(y,y',z)\leq C|y-y'|^{\frac{1}{2}(2\gamma -Q)^2} |y-1|^{1-\frac{3\gamma^2}{4}}
  |z|^{\frac{1}{2}(\frac{3\gamma}{2}+\alpha_1-Q)^2  }.$$ 
 since $\frac{1}{2}(\frac{3\gamma}{2} -Q)^2-\frac{\gamma^2}{8}-\frac{1}{2}(2\gamma -Q)^2=1-\frac{3\gamma^2}{4}$.
  Hence
   \begin{align*}
 I(B\cap A_3)\leq &C|z|^{-2+\gamma(Q-\alpha_1)+\hf(\alpha_1-\frac{2}{\gamma})^2} \int_{B(1,
 |z|^\epsilon)^2}\frac{|y-1|^{-\frac{\gamma^2}{4}}}{
 |y-y'|^{\gamma^2-\frac{1}{2}(2\gamma -Q)^2} } \,d^2yd^2y' .
 \end{align*} 
 The integral converges since $\gamma^2-\frac{1}{2}(2\gamma-Q)^2=4-\frac{1}{2} Q^2<2$.\qed

\subsection{Fusion with reflection} 
In this section we uncover the probabilistic origin of the reflection relation \eqref{rrel}, \eqref{rrel1}. Notice that the restriction $\alpha_1<\frac{2}{\gamma} $ in Theorem \ref{theo4point} comes from the second three point structure constant $C_\gamma(\alpha_1+\tfrac{\gamma}{2}, \alpha_2,\alpha_3)$  in the expression  \eqref{theo4pointexpressionBPZ1}: this condition is required in order that the first weight $\alpha_1+\tfrac{\gamma}{2}$ is consistent with condition $\alpha_1+\tfrac{\gamma}{2}<Q$ of the Seiberg bound \eqref{TheSeibergbounds}. We prove the following extension of Theorem \ref{theo4point} to the case $\alpha_1>\frac{2}{\gamma}  $.

\begin{theorem}\label{theo4point1}
Let $\sum_{k=1}^3 \alpha_k> 2 Q+\frac{\gamma}{2} $ and $\alpha_k<Q$ for all $k$. There exists $\eta>0$ s.t. if $Q-\alpha_1<\eta$ then 
\begin{equation}\label{BPZ1}
\mathcal{T}_{-\frac{\gamma}{2}}(z) = C_\gamma(\alpha_1-\frac{\gamma}{2}, \alpha_2,\alpha_3)  |F_{-}(z)|^2  + R(\alpha_1) C_\gamma(2Q-\alpha_1-\frac{\gamma}{2}, \alpha_2,\alpha_3)  |F_{+}(z)|^2.
\end{equation}
\end{theorem}

\proof 

By the discussion in subsection \ref{BPZ equations}, it suffices to show that one has the following expansion as $z$ goes to $0$
\begin{align}\label{expansionmathcalTwithR}
\mathcal{T}_{-\frac{\gamma}{2}}(z)=C_\gamma(\alpha_1-\tfrac{\gamma}{2}, \alpha_2,\alpha_3)+  R(\alpha_1) C_\gamma(2Q-\alpha_1-\frac{\gamma}{2}, \alpha_2,\alpha_3) |z|^{\gamma (Q-\alpha_1)}+ o( |z|^{\gamma (Q-\alpha_1)}).
\end{align}
Note that since now $\gamma(Q-\alpha_1)<1$ we need a Taylor expansion only to $0$-th order. We use the notations introduced in the proof of \eqref{theo4pointexpressionBPZ1}. Recall that 
\begin{equation*}
K(z,x)=\frac{|x-z|^{\frac{\gamma^2}{2}}|x|_+^{\gamma (\sum_{k=1}^3\alpha_k-\frac{\gamma}{2})}}{|x|^{\gamma \alpha_1}|x-1|^{\gamma \alpha_2}}
\end{equation*}
and for any Borel set $B\subset \C$
\begin{equation*}
\mathcal{K}_B(z)=\int_{B}K(z,x)\,M_\gamma(d^2x).
\end{equation*}
Recall also that we write $\mathcal{K}(z)$ for $\mathcal{K}_\C(z)$, use the notation $B_{r}=B(0,r)$ and set $p=\frac{1}{\gamma}(\sum_{k=1}^3 \alpha_k-\frac{\gamma}{2}-2Q)$. Since
\begin{align}\label{Tdefi000bis}
  \mathcal{T}_{-\frac{\gamma}{2}}(z) &=
  2 \mu^{-p}  \gamma^{-1}\Gamma(p)  \E[\mathcal{K}(z)^{-p}] 
 \end{align}
and $ \mathcal{T}_{-\frac{\gamma}{2}}(0)=  C_\gamma(\alpha_1-\tfrac{\gamma}{2}, \alpha_2,\alpha_3)$, in order to get \eqref{expansionmathcalTwithR}, it suffices to prove that  
\begin{equation*}
\E\big[\mathcal{K}(z)^{-p}\big] -\E\big[ \mathcal{K}(0)^{-p}\big]=\frac{1}{2}\mu^{p}\gamma  \Gamma(p)^{-1}R(\alpha_1) C_\gamma(2Q-\alpha_1-\frac{\gamma}{2}, \alpha_2,\alpha_3)|z|^{\gamma(Q-\alpha_1)}+o(|z|^{\gamma(Q-\alpha_1)}).
\end{equation*}
The leading asymptotics will result from   the integral defining $\caK$ in a small ball at the origin. Let us denote $B:=B_{ |z|^{1-\xi}}=B(0, |z|^{1-\xi}) $ with  $\xi\in (0,1)$ to be fixed later.   We define 
\begin{equation}\label{T1T2def}T_1:= \E\big[   \mathcal{K}_{B^c}(z)^{-p}\big]     -  \E\big[ \mathcal{K}(0)^{-p}\big]\quad\text{ and }\quad T_2:=  \E\big[ \mathcal{K}(z)^{-p}\big]-\E\big[   \mathcal{K}_{B^c}(z)^{-p}\big]   
\end{equation}
so that  
\begin{equation}\label{T1T2split}
\E\big[\mathcal{K}(z)^{-p}\big]     -  \E\big[ \mathcal{K}(0)^{-p}\big]=T_1+T_2.
\end{equation}

We then get the desired result thanks to the following two lemmas (where $\xi$ will be fixed in the proof of the two lemmas):

\begin{lemma}\label{lemmaintbis1}
The following holds 
\begin{equation*}
T_1=o(|z|^{\gamma(Q-\alpha_1)}).
\end{equation*}
\end{lemma}

\begin{lemma}\label{lemmaintbis2}
The following holds 
\begin{equation*}
T_2= \frac{1}{2}\mu^{p}\gamma  \Gamma(p)^{-1}R(\alpha_1) C_\gamma(2Q-\alpha_1-\frac{\gamma}{2}, \alpha_2,\alpha_3)|z|^{\gamma(Q-\alpha_1)}+o(|z|^{\gamma(Q-\alpha_1)}).
\end{equation*}
\end{lemma}

\qed

\noindent
\emph{Proof of Lemma \ref{lemmaintbis1}}
By interpolation, we get
\begin{align}\nonumber
|T_1|& \leq  p\int_0^1\E\Big[|\mathcal{K}_{B^c}(z)- \mathcal{K}(0)|
\Big( t\mathcal{K}_{B^c}(z)+(1-t)\mathcal{K}(0)\Big)^{-p-1}\Big]dt\\&\leq C\E\Big[|\mathcal{K}_{B^c}(z)- \mathcal{K}(0)|
\Big(\mathcal{K}_{B^c}(0)\Big)^{-p-1}\Big] \label{interpol}
\end{align}
where we used  $\mathcal{K}_{B^c}(z)\geq C  \mathcal{K}_{B^c}(0)$ and $ \mathcal{K}(0)\geq \mathcal{K}_{B^c}(0)$ since  $|x-z|^{\frac{\gamma^2}{2}}\geq C|x|^{\frac{\gamma^2}{2}}$ on $B^c$.  Since $\mathcal{K}(0)=\mathcal{K}_B(0)+\mathcal{K}_{B^c}(0)$ we obtain
$|T_1| \leq 
C(A_1+A_2)$
where 
\begin{equation*}
A_1 = \E[   \mathcal{K}_{B} (0) \mathcal{K}_{B^c}(0) ^{-p-1}  ]\quad \text{ and }\quad 
A_2 = \E[    | \mathcal{K}_{B^c}(z)-\mathcal{K}_{B^c}(0) | \mathcal{K}_{B^c}(0) ^{-p-1}    ].
\end{equation*}
 Using Cameron-Martin theorem, we get for $A_1$ that
\begin{align*}
 A_1   &   \leq C \int_{|x| \leq |z|^{1-\xi}} |x|^{\frac{\gamma^2}{2}-\gamma \alpha_1}\E\Big[\Big(\int_{|u|>|z|^{1-\xi}}K(0,u)|u-x|^{-\gamma^2}M_\gamma(du)\Big)^{-p-1}\Big] \,d^2x.
\end{align*}
Since $|u-x|\leq 2|u|$ for $|u|> |z|^{1-\xi}$ we may bound the expectation by
\begin{align}\label{aaaa}
 \E\Big[\Big(\int_{|u|>|z|^{1-\xi}}|u|^{-\gamma\alpha_1-\frac{\gamma^2}{2}}M_\gamma(d^2u)\Big)^{-p-1}\Big] \leq C|z|^{(1-\xi)\frac{(\alpha_1+\frac{\gamma}{2}-Q)^2}{2}}
\end{align}
where we used the GMC estimate \eqref{freezingannulus} in the Appendix. 
We conclude that
\begin{align}\label{A1bound}A_1\leq C|z|^{(1-\xi)\big(\frac{(\alpha_1+\frac{\gamma}{2}-Q)^2}{2}+ \gamma(Q-\alpha_1)\big)}.
\end{align}
Hence $A_1=o(|z|^{\gamma(Q-\alpha_1)})$ if e.g. $\xi<\frac{1}{2}$ and $\eta$ is small enough.

 Next we bound $A_2$.
Let $A$ be the annulus $A:=\{x\in\C; |z|^{1-\xi}\leq |x|\leq 1/2\}$. We can split the numerator in $A_2$ into $ | \mathcal{K}_{B^c_{1/2}}(z)-\mathcal{K}_{B^c_{1/2}}(0) | $ and $ | \mathcal{K}_{A}(z)-\mathcal{K}_{A}(0) |$ by means of the triangular inequality. On $B^c_{1/2}$ we can use $||x-z|^{\frac{\gamma^2}{2}}- |x|^{\frac{\gamma^2}{2}}| \leq C |x|^{\frac{\gamma^2}{2}} |z| $ to get
\begin{align*}
\E[  | \mathcal{K}_{B^c_{1/2}}(z)-\mathcal{K}_{B^c_{1/2}}(0) |  \mathcal{K}_{B^c}(0) ^{-p-1}    ]& \leq C  |z|  \E[   \mathcal{K}_{B^c_{1/2}}(0)  \mathcal{K}_{B^c}(0) ^{-p-1}     ]  \leq C  |z|  \E[     \mathcal{K}_{B^c}(0) ^{-p}     ]  \leq C |z| .
\end{align*}

Finally, using  $||x-z|^{\frac{\gamma^2}{2}}- |x|^{\frac{\gamma^2}{2}}| \leq C |x|^{\frac{\gamma^2}{2}-1} |z| $ on $A$ and then applying Cameron-Martin, we get
\begin{align*}
\E[ | \mathcal{K}_{A}(z)-\mathcal{K}_{A}(0) | \mathcal{K}_{B^c}(0) ^{-p-1}    ] & \leq C |z| \int_{A} {|x|^{\frac{\gamma^2}{2}-1-\gamma\alpha_1}}  \E\Big[\Big(\int_{|u|>|z|^{1-\xi}}K(0,u)|u-x|^{-\gamma^2}M_\gamma(d^2u)\Big)^{-p-1}\Big]  d^2 x.
\end{align*}
Since $|z|^{1-\xi}\leq |x|$ we can bound
\begin{align*}
\int_{|u|>|z|^{1-\xi}}K(0,u)|u-x|^{-\gamma^2}M_\gamma(d^2u)\geq C\int_{|u|>|x|}|u|^{-\gamma\alpha_1-\frac{\gamma^2}{2}}M_\gamma(d^2u)
\end{align*}
and then the GMC estimate  \eqref{freezingannulus} in the Appendix gives
\begin{align}\label{A2bound}
\E[ | \mathcal{K}_{A}(z)-\mathcal{K}_{A}(0) | \mathcal{K}_{B^c}(0) ^{-p-1}    ] 
& \leq C |z| \int_{A} {|x|^{\frac{\gamma^2}{2}-1-\gamma\alpha_1+\frac{1}{2} (\alpha_1+\frac{\gamma}{2}-Q)^2}}d^2x\\&\leq C|z|^{\xi+  (1-\xi)\big(\gamma(Q-\alpha_1)+\frac{1}{2} (\alpha_1-\frac{2}{\gamma})^2\big)}=o(|z|^{\gamma(Q-\alpha_1)})
\end{align}
for $\xi<\hf$ and $\eta$ small enough (since $\alpha_1-\frac{2}{\gamma}>\frac{\gamma}{2}-\eta$). Hence $T_1=o(|z|^{\gamma(Q-\alpha_1)})$.

\qed

\bigskip

\noindent
\emph{Proof of Lemma \ref{lemmaintbis2}}
First we show that it suffices to restrict $\mathcal{K}(z)$  to the complement of the annulus  $A_{h}:=\{x\in\C; e^{-h}|z|\leq |x|\leq |z|^{1-\xi}\}$ where $h>0$ is fixed: it will serve as a buffer zone to decorrelate the regions $\{x\in\C; |x|\leq e^{-h}|z|\}$ and $\{x\in\C; |x|> |z|^{1-\xi}\}$. Interpolating as in \eqref{interpol} we deduce
\begin{align*}
 | \E[  \mathcal{K}(z)^{-p}    -   \mathcal{K}_{A_h^c}(z)^{-p}   ]  | &  \leq \E[     \mathcal{K}_{A_h }(z) \mathcal{K}_{B^c}(0)^{-p-1}  ]  .
\end{align*} 
Using the Cameron-Martin theorem we get
\begin{align}\label{T2prelim0}
 | \E[  \mathcal{K}(z)^{-p}    -   \mathcal{K}_{A_h^c}(z)^{-p}   ]  | &  \leq C\int_{A_h }|x-z|^{\frac{\gamma^2}{2}}|x|^{-\gamma\alpha_1}\E\Big[\Big(\int_{B^c} K(0,u)|u-x|^{-\gamma^2}M_\gamma(du)   \Big)^{-p-1}\Big]  d^2x.
\end{align} 
The  expectation  was  estimated in \eqref{aaaa} so that we get
 \begin{align}\label{T2prelim}| \E[  \mathcal{K}(z)^{-p}    -   \mathcal{K}_{A_h^c}(z)^{-p}   ]  | \leq C|z|^{(1-\xi)\big(\gamma(Q-\alpha_1)+\frac{1}{2}(Q-\alpha_1-\frac{\gamma}{2})^2\big)}.
 \end{align} 
For $\xi<\hf$  and $\eta$ small this yields 
\begin{equation}\label{FirststepT2}
| \E[  \mathcal{K}(z)^{-p}    -   \mathcal{K}_{A_h^c}(z)^{-p}   ]  | =o(|z|^{\gamma(Q-\alpha_1)}).
\end{equation} 

\medskip  
 
Therefore, we just need to evaluate the quantity $\E[ \mathcal{K}_{A_h^c}(z)^{-p} ]-\E[ \mathcal{K}_{B^c}(z)^{-p} ]$ where we recall the definitions $B^c=\{|x|\geq |z|^{1-\xi}\}$ and $A_h^c=B^c\cup B_{ e^{-h}|z|}$. Hence $\mathcal{K}_{A_h^c}(z)= \mathcal{K}_{B^c}(z)+ \mathcal{K}_{B_{ e^{-h}|z|} }(z)$. We use the polar decomposition of the chaos measure  introduced in Section \ref{Reflection coefficient}. Let   ${|z|}=e^{-t}$. Then 
\begin{align*}
 \mathcal{K}_{B^c} (z)  &=  \int_0^{2\pi}\int_{-\infty}^{(1-\xi)t}  e^{\gamma (B_s-  (Q-\alpha_1)s)}  \frac{|e^{-s+i\theta}- z|^{\frac{\gamma^2}{2}}}{|1-e^{-s+i\theta}|^{\gamma \alpha_2}} (e^{-s(\gamma(\alpha_1+\alpha_2+\alpha_3-\frac{\gamma}{2})}\vee 1) N_\gamma(dsd\theta):= \mathcal{K}^1\\
 \mathcal{K}_{B_{ e^{-h}|z|} }(z)  &= 
   \int_0^{2\pi}\int_{t+h}^\infty  e^{\gamma (B_s-  (Q-\alpha_1)s)}  \frac{|e^{-s+i\theta}- z|^{\frac{\gamma^2}{2}}}{|1-e^{-s+i\theta}|^{\gamma \alpha_2}}  N_\gamma(dsd\theta):= \mathcal{K}^2.
\end{align*}
The lateral noises $Y$ which enter the definition of $N_\gamma(dsd\theta)$ in $\mathcal{K}^1$ and $\mathcal{K}^2$ are weakly correlated. Indeed, from \eqref{covlateral} we get 
\begin{equation}\label{majkah}
-e^{-\xi t}\leq E[Y(s,\theta)Y(s',\theta')] \leq 2e^{-\xi t}.
\end{equation}
for all $s<(1-\xi)t,s'>t+h$ and $\theta,\theta'\in[0,2\pi]$. Define then the process
$$P(s,\theta):= Y(s,\theta)\mathbf{1}_{\{s<(1-\xi)t\}} + Y(s,\theta)\mathbf{1}_{\{s>t+h\}}.$$
Let 
 $\widetilde{Y}$ be independent of everything with the same law as $Y$ and define the process
$$\tilde P(s,\theta):= Y(s,\theta)\mathbf{1}_{\{s<(1-\xi)t\}} + \tilde Y(s,\theta)\mathbf{1}_{\{s>t+h\}}.$$ 
Then we get
 \begin{equation}\label{Compkernels}
 \E[\tilde P(s,\theta)\tilde P(s',\theta')]-e^{-\xi t}\leq \E[P(s,\theta)P(s',\theta')]\leq\E[\tilde P(s,\theta)\tilde P(s',\theta')]+2e^{-\xi t}.
 \end{equation}
 Let $N$ be  a unit normal variable independent of everything. Then inequality \eqref{Compkernels} implies  that the covariance of $P+e^{-\frac{1}{2}\xi t}N$ dominates the covariance of $\tilde{P}$ and the covariance of $ \tilde{P}+\sqrt{2}e^{-\frac{1}{2} \xi t}N$ dominates the covariance of $P$. Therefore, we get 
 by Kahane's convexity inequality (see \cite[Theorem 2.1]{review}) with the convex function $x\in\R_+\mapsto x^{-p}$ (applied to $(P+e^{-\frac{1}{2}\xi t}N, \tilde{P})$ and $ (\tilde{P}+\sqrt{2}e^{-\frac{1}{2} \xi t}N,P)$) that there exists some $C>0$ such that 
\begin{align}\label{SecondstepT2}
e^{-C|z|^{\xi}}\E[  ( \mathcal{K}^1+   \tilde{\mathcal{K}}^2)^{-p}] \leq \E[  ( \mathcal{K}^1+   \mathcal{K}^2)^{-p}]\leq e^{C|z|^{\xi}}
\E[  ( \mathcal{K}^1+  %
\tilde{\mathcal{K}}^2)^{-p}]
\end{align} 
where $\tilde{\mathcal{K}}^2$  is computed with $\tilde Y$ instead of $Y$.  Let  
\begin{align}\label{epsilondefinition}
\beta
:=e^{\gamma B_{t+h}-\gamma (Q-\alpha_1)(t+h)- \frac{\gamma^2}{2}t}   .
\end{align} 
Then by  the Markov property of Brownian motion
\begin{align}\label{refKSzh} \tilde{\mathcal{K}}^2=& 
 \beta\int_0^{2\pi}\!\!\int_{0}^\infty  e^{\gamma (\widetilde{B}_s-  (Q-\alpha_1)s)}  \frac{|e^{-s-h+i\theta}- \tfrac{z}{|z|}|^{\frac{\gamma^2}{2}}}{|1-|z|e^{-s-h+i\theta}|^{\gamma \alpha_2}}   \tilde N_\gamma(d(h+t+s),d\theta)
 \end{align}
where $\widetilde{B}$ is a Brownian motion independent of everything and $\tilde{N}$ is the measure associated to $\tilde{Y}$. Moreover, by stationarity of $\tilde Y$ and its independence of everything we  may replace $\tilde{N}_\gamma(d(h+t+s),d\theta)$ by $\tilde{N}_\gamma(ds,d\theta)$. As a consequence
\begin{align}\label{uandlbounds}
 \E[  ( \mathcal{K}^1+  \beta c_- \mathcal{K}^3)^{-p}]\leq \E[  ( \mathcal{K}^1+   \tilde{\mathcal{K}}^2)^{-p}] \leq \E[  ( \mathcal{K}^1+  \beta c_+ \mathcal{K}^3)^{-p}]
\end{align} 
where 
\begin{align*}
\mathcal{K}^3=\int_{0}^\infty  e^{\gamma (\widetilde{B}_s-  (Q-\alpha_1)s)} \widetilde{Z}_{s}ds
\end{align*} 
with $\widetilde{Z}_{s}= \int_0^{2\pi}  e^{\gamma \tilde{Y}(s,\theta)-\frac{\gamma^2 \E[  (\tilde{Y}(s,\theta))^2 ]} {2}} d\theta$ (recall that this is a slight abuse of notation as $\widetilde{Z}_{s}$ is not a function but a distribution) and 
\begin{equation}\label{defCL}
c_\pm:=\frac{(1\mp e^{-h})^{\frac{\gamma^2}{2}}}{(1\pm |z|e^{-h})^{\gamma\alpha_2}}.
\end{equation}
By the Williams  path decomposition Lemma \ref{lemmaWilliams} and \eqref{willaplied} we deduce 
\begin{equation}\label{k3will}
\mathcal{K}^3\stackrel{law}=e^{\gamma M}\int_{-L_{-M}}^\infty e^{\gamma \mathcal{B}_s^{\alpha_1}} \widetilde{Z}_{s}\,ds
\end{equation}
where we recall $M=\sup_s(\tilde{B}_s-  (Q-\alpha_1)s)$ and $ L_{-M}$ is the last time  $\mathcal{B}_s^{\alpha_1}$ hits $-M$ (along the negative axis).
Thanks to \eqref{SecondstepT2} and \eqref{uandlbounds}, we want to study a lower bound on $\E[  ( \mathcal{K}^1+  \beta c_- \mathcal{K}^3)^{-p}]$ and an upper bound on $\E[  ( \mathcal{K}^1+  \beta c_+ \mathcal{K}^3)^{-p}]$ to conclude.

\subsubsection{Lower bound on $\E[  ( \mathcal{K}^1+  \beta c_- \mathcal{K}^3)^{-p}]$} Let us use the notation $J_A=\int_{-L_{-A}}^\infty e^{\gamma  \mathcal{B}_s^{\alpha_1}} \widetilde{Z}_sds$ and $J$ for $J_{\infty}$.
 We have
\begin{equation*}
 \E[  ( \mathcal{K}^1+  \beta c_- \mathcal{K}^3)^{-p}]\geq  \E[  ( \mathcal{K}^1+  \beta c_-e^{\gamma M} J)^{-p}].
\end{equation*} 
Using the fact  that $M$ has  exponential law with parameter $2(Q-\alpha_1)$ (and therefore the law of $e^{\gamma M}$ has density $\frac{_{2(Q-\alpha_1)}}{^\gamma}  v^{-1-\frac{2}{\gamma}(Q-\alpha_1)}\,dv $ on $[1,\infty)$), we get by first integrating over $M$ that
\begin{align*}
& \E[  ( \mathcal{K}^1+  \beta c_- \mathcal{K}^3)^{-p}]- \E[  ( \mathcal{K}^1)^{-p}] = \frac{_{2(Q-\alpha_1)}}{^\gamma}\E\Big[  \int_{1}^\infty   \Big(  ( \mathcal{K}^1+  \beta c_-  v
J)^{-p}  - (\mathcal{K}^1)^{-p}  \Big ) v^{-1-\frac{2}{\gamma}(Q-\alpha_1)}\,dv \Big ]\\
= &\frac{_{2(Q-\alpha_1)}}{^\gamma} c_-^{\frac{2}{\gamma}(Q-\alpha_1)}\E\Big[(\beta J)^{\frac{2}{\gamma}(Q-\alpha_1)}(\mathcal{K}^1)^{-p-\frac{2}{\gamma}(Q-\alpha_1)} \int_{\frac{\beta c_-J}{ \mathcal{K}^1}}
^\infty    \Big((1 +  w)^{-p}  -  1 \Big)w^{-1-\frac{2}{\gamma}(Q-\alpha_1)}\,dw \Big ] \\
\geq &\frac{_{2(Q-\alpha_1)}}{^\gamma} c_-^{\frac{2}{\gamma}(Q-\alpha_1)}\E\Big[(\beta J)^{\frac{2}{\gamma}(Q-\alpha_1)}(\mathcal{K}^1)^{-p-\frac{2}{\gamma}(Q-\alpha_1)} \Big ]\int_0
^\infty    \Big((1 +  w)^{-p}  -  1 \Big)w^{-1-\frac{2}{\gamma}(Q-\alpha_1)}\,dw\\
= &\frac{_{2(Q-\alpha_1)}}{^\gamma} c_-^{\frac{2}{\gamma}(Q-\alpha_1)}\frac{\Gamma(-\frac{2}{\gamma} (Q-\alpha_1)) \Gamma (p+\frac{2}{\gamma} (Q-\alpha_1) )}{\Gamma (p)}\E[ J^{\frac{2}{\gamma}(Q-\alpha_1)}]\E\Big[\beta^{\frac{2}{\gamma}(Q-\alpha_1)}(\mathcal{K}^1)^{-p-\frac{2}{\gamma}(Q-\alpha_1)} \Big ]
\end{align*} 
where in the second step we made a change of variables  $w=\frac{\beta c_-J}{ \mathcal{K}^1} v$ and for the lower bound we took the integration over $w\geq 0$. In the last step we  used   Lemma \ref{lemmaintegral} in the Appendix  to compute the integral
 and we also used the independence of $J$ from everything. We end up with
\begin{align}\label{ThirdstepT2}
& \E[  ( \mathcal{K}^1+  \beta c_- \mathcal{K}^3)^{-p}]- \E[  ( \mathcal{K}^1)^{-p}]\geq W
c_-^{\frac{2}{\gamma}(Q-\alpha_1)}\E\Big[\beta^{\frac{2}{\gamma}(Q-\alpha_1)}(\mathcal{K}^1)^{-p-\frac{2}{\gamma}(Q-\alpha_1)} \Big ]
\end{align} 
where we have have set 
 $$W:= 
\mu^{-\frac{2}{\gamma} (Q-\alpha_1)}R(\alpha_1) \frac{  \Gamma (p+\frac{2}{\gamma} (Q-\alpha_1) )}{\Gamma (p)}$$
and $R(\alpha_1)= \mu^{\frac{2}{\gamma} (Q-\alpha_1)} \Gamma (-\tfrac{2(Q-\alpha_1)}{\gamma})\  \frac{2(Q-\alpha_1)}{\gamma}  \E[   J^{\frac{2}{\gamma}(Q-\alpha_1)}] $ is the reflection coefficient defined in \eqref{deffullR}. We point out that $W$ is negative because of the term  $\Gamma (-\tfrac{2(Q-\alpha_1)}{\gamma})$ appearing in the expression of the reflection coefficient. 
 The remaining expectation  can be computed thanks to the Cameron-Martin theorem applied to the term $\beta^{\frac{2}{\gamma}(Q-\alpha_1)}$ (defined via \eqref{epsilondefinition}). More precisely, we have (using $t=\ln \frac{1}{|z|}$)
 \begin{equation*}
 \beta^{\frac{2}{\gamma}(Q-\alpha_1)}= e^{2 (Q-\alpha_1) B_{t+h}-2 (Q-\alpha_1)^2(t+h)}e^{- \gamma(Q-\alpha_1) t} = |z|^{\gamma(Q-\alpha_1)} e^{ 2 (Q-\alpha_1) B_{t+h}-2 (Q-\alpha_1)^2(t+h)  }
 \end{equation*}
 and we get by applying the Cameron-Martin theorem to the term $e^{ 2 (Q-\alpha_1) B_{t+h}-2 (Q-\alpha_1)^2(t+h)  }$
 \begin{align}\label{FourthstepT2}
\E\Big[\beta^{\frac{2}{\gamma}(Q-\alpha_1)}(\mathcal{K}^1)^{-p-\frac{2}{\gamma}(Q-\alpha_1)} \Big ]=|z|^{\gamma(Q-\alpha_1)}\E\big[ \widehat{ \mathcal{K}}_{B^c}(z)^{-p-\frac{2}{\gamma}(Q-\alpha_1)}\big]
\end{align} 
 where we defined for $D\subset \C$
 \begin{equation}\label{defKhat}
\widehat{ \mathcal{K}}_{D}(z):=\int_{D}\frac{|x-z|^{\frac{\gamma^2}{2}}}{|x|^{\gamma(2Q-\alpha_1)}|x-1|^{\gamma\alpha_2}}|x|_+^{\gamma (2Q-\alpha_1-\frac{\gamma}{2}+\alpha_2+\alpha_3)}M_\gamma(d^2x)
\end{equation}
and in the case $D=\C$ we will write $\widehat{ \mathcal{K}}(z)$ for $\widehat{ \mathcal{K}}_{\C}(z)$. 
Next, we claim
\begin{align}\label{hatKbound}
\E\big[ \widehat{ \mathcal{K}}_{B^c}(z)^{-p-\frac{2}{\gamma}(Q-\alpha_1)}\big]-\E\big[ \widehat{ \mathcal{K}}_{B^c}(0)^{-p-\frac{2}{\gamma}(Q-\alpha_1)}\big]=o( |z|^{\gamma(Q-\alpha_1)}).
\end{align} 
Indeed, the LHS is just $T_1$ in \eqref{T1T2def} computed with a larger $p$ and $|x|^{\gamma\alpha_1}$ replaced by $|x|^{\gamma(2Q-\alpha_1)}$. It is readily checked from \eqref{A1bound} and \eqref{A2bound} that \eqref{hatKbound} holds.

In view of \eqref{FirststepT2}, \eqref{SecondstepT2}, \eqref{uandlbounds},  \eqref{ThirdstepT2},  \eqref{FourthstepT2}, \eqref{hatKbound} we have shown that
\begin{align*}
\E[ & \mathcal{K}(z)^{-p} ]-\E[ \mathcal{K}_{B^c}(z)^{-p} ]\\
\geq  &e^{
-C|z|^{\xi }}\Big( |z|^{\gamma(Q-\alpha_1)} c_-^{\frac{2}{\gamma}(Q-\alpha_1)}W  \E\big[ \widehat{ \mathcal{K}}(0)^{-p-\frac{2}{\gamma}(Q-\alpha_1)}\big]+o( |z|^{\gamma(Q-\alpha_1)})\Big)
  -(1-e^{-C
|z|^{\xi }})\E[ \mathcal{K}_{B^c}(z)^{-p} ].
\end{align*} 
The second term on the RHS is  $\caO(|z|^\xi)=o( |z|^{\gamma(Q-\alpha_1)})$ {\it provided} we take $\xi>\gamma(Q-\alpha_1)$ (this is the condition that fixes $\xi$) so that 
recalling \eqref{defCL}, we deduce
\begin{align*}
\liminf_{z\to 0}|z|^{-\gamma(Q-\alpha_1)}\Big(\E[ &  \mathcal{K}(z)^{-p} ] -\mathcal{K}_{B^c}(z)^{-p} ]\Big) \geq    (1+e^{-h})^{ \gamma (Q-\alpha_1)}W  \E\big[ \widehat{ \mathcal{K}}(0)^{-p-\frac{2}{\gamma}(Q-\alpha_1)}\big].
\end{align*} 
Since $h$ is arbitrary, it can be chosen arbitrarily large so as to get 
\begin{align*}
 \liminf_{z\to 0}|z|^{-\gamma(Q-\alpha_1)} \Big(\E[ \mathcal{K}(z)^{-p}  -\mathcal{K}_{B^c}(z)^{-p} ]\Big) 
&  \geq   W  \E\big[ \widehat{ \mathcal{K}}(0)^{-p-\frac{2}{\gamma}(Q-\alpha_1)}\big]\\
&=\frac{1}{2}\mu^{p}\gamma  \Gamma(p)^{-1}R(\alpha_1) C_\gamma(2Q-\alpha_1-\frac{\gamma}{2}, \alpha_2,\alpha_3)
\end{align*}
where we have used the definition \eqref{expression3pointstruct} of the structure constants in the last equality. This is the desired lower bound.

\subsubsection{Upper bound on $\E[  ( \mathcal{K}^1+  \beta c_+ \mathcal{K}^3)^{-p}]$}  For the upper bound we go back to the formula \eqref{k3will} where we need to face the integration region lower value $L_{-M}$.  For $A>0$ fixed, we consider first the quantity
$$
L(z):=\E\Big[   \Big( \mathcal{K}^1  +  \beta c_+ e^{\gamma M}\int_{-L_{-M}}^\infty e^{\gamma \mathcal{B}_s^{\alpha_1}} \widetilde{Z}_{s} ds\Big)^{-p}  - (\mathcal{K}^1)  ^{-p}\Big)\mathbf{1}_{\{M\leq A\}} \Big ]$$ and we want to show that $L(z)=o(|z|^{\gamma(Q-\alpha_1)})$. Using the inequality $|(x+y)^{-p}- x^{-p}| \leq p x^{-p-1}y$ for $x,y>0$ we get
\begin{equation*}
|L(z)|\leq p c_+\E\big[  \beta e^{\gamma M}\int_{-L_{-M}}^\infty e^{\gamma \mathcal{B}_s^{\alpha_1}} \widetilde{Z}_{s}ds\,\mathcal{K}_{B^c}(z)^{-p-1}\mathbf{1}_{\{M\leq A\}}\big] \leq C e^{\gamma A}\E\big[ \int_{\R}e^{\gamma \mathcal{B}_s^{\alpha_1}} \widetilde{Z}_{s}ds\big]\E\big[  \beta \mathcal{K}_{B^c\cap B(0,\hf)}(z)  ^{-p-1} \big].
 \end{equation*}
 Recall that $\beta$ satisfies (using $t=\ln \frac{1}{|z|}$) 
 \begin{equation*}
 \beta= e^{\gamma B_{t+h}-\gamma (Q-\alpha_1)(t+h)- \frac{\gamma^2}{2}t} =  e^{\gamma (\alpha_1-\frac{2}{\gamma}) h}  |z|^{\gamma(Q-\alpha_1)} e^{\gamma B_{t+h}- \frac{\gamma^2}{2}(t +h)}  
 \end{equation*}
 and therefore using the Cameron martin theorem with $e^{\gamma B_{t+h}- \frac{\gamma^2}{2}(t +h)}   $ we get 
 \begin{align*}
   &  e^{\gamma A}\E\big[ \int_{\R}e^{\gamma \mathcal{B}_s^{\alpha_1}} \widetilde{Z}_{s}ds\big]\E\big[  \beta \mathcal{K}_{B^c\cap B(0,\hf)}(z)  ^{-p-1} \big]   \\
&   =  e^{\gamma (\alpha_1-\frac{2}{\gamma}) h}  
|z|^{\gamma(Q-\alpha_1)} e^{\gamma A}\E\big[ \int_{\R} e^{\gamma \mathcal{B}_s^{\alpha_1}} \widetilde{Z}_{s}\,ds\big]\E\big[    \Big(\int_{|z|^{1-\xi}\leq|z|\leq\hf}|x|^{-\gamma\alpha_1-\frac{\gamma^2}{2}}M_\gamma(d^2x)\Big)^{-p-1} \big]\\
&\leq  C
 e^{\gamma A}|z|^{\gamma(Q-\alpha_1)}|z|^{\frac{1-\xi}{2}(Q-\alpha_1-\frac{\gamma}{2})^2}
 \end{align*}
where the constant $C$ depends on $h$ and we used the GMC estimate \eqref{freezingannulus} and Lemma \ref{defmomR}. Hence, for $A$ and $h$ fixed, we have $L(z)=o(|z|^{\gamma(Q-\alpha_1)})$. 
	
\medskip
It remains to investigate the quantity 
\begin{align*}
U(z):=&\E\Big[   \Big( \mathcal{K}^1  +  \beta c_+ e^{\gamma M}\int_{-L_{-M}}^\infty e^{\gamma \mathcal{B}_s^{\alpha_1}} \widetilde{Z}_{s}(ds)\Big)^{-p}  - (\mathcal{K}^1)  ^{-p}\Big)\mathbf{1}_{\{M\geq A\}} \Big ]\\
\leq & \E\Big[   \Big( \mathcal{K}^1  +  \beta c_+ e^{\gamma M}J_A \Big)^{-p}  -(\mathcal{K}^1)  ^{-p} \Big)\mathbf{1}_{\{M>A\}} \Big ]
 \end{align*}
 where recall that $J_A=\int_{-L_{-A}}^\infty e^{\gamma \mathcal{B}_s^{\alpha_1}} \widetilde{Z}_{s}ds$.
Using again the law of $M$, which is exponential with parameter $2(Q-\alpha_1)$, and making the change  of variables $  \frac{\beta c_+ J_A}{ \mathcal{K}^1}e^{\gamma v}=y$ we get 
\begin{align*}
U(z)\leq & \frac{_{2(Q-\alpha_1)}}{^\gamma}\E\Big[ \int_{A}^\infty  \Big( \Big( \mathcal{K}^1 +  \beta c_+ e^{\gamma v}J_A\Big)^{-p}  -( \mathcal{K}^1)^{-p}\Big)  e^{-2(Q-\alpha_1)v}\,dv\Big ] \\
=&  \frac{_{2(Q-\alpha_1)}}{^\gamma}c_+^{\frac{2}{\gamma}(Q-\alpha_1)}\E\Big[J_A^{\frac{2}{\gamma}(Q-\alpha_1)}\beta^{\frac{2}{\gamma}(Q-\alpha_1)}\int_{e^{\gamma A}\frac{ \beta c_+J_A}{ \mathcal{K}^1}}^\infty  \Big(  (1 + y )^{-p}  -  1\Big) (\mathcal{K}^1)  ^{-p-\frac{2}{\gamma}(Q-\alpha_1)} y^{-\frac{2}{\gamma}(Q-\alpha_1)-1}\,dy\Big ].
 \end{align*}
Now we can   use Cameron-Martin as in the case of the lower bound to get that the above expectation can be rewritten as (recall \eqref{defKhat})
\begin{align*}
\E&\Big[J_A^{\frac{2}{\gamma}(Q-\alpha_1)}\beta^{\frac{2}{\gamma}(Q-\alpha_1)}\int_{e^{\gamma A}\frac{ \beta c_+J_A}{ \mathcal{K}^1}}^\infty  \Big(  (1 + y )^{-p}  -  1\Big)( \mathcal{K}^1)  ^{-p-\frac{2}{\gamma}(Q-\alpha_1)} y^{-\frac{2}{\gamma}(Q-\alpha_1)-1}\,dy\Big ]\\
=&|z|^{\gamma(Q-\alpha_1)}\E \Big[J_A^{\frac{2}{\gamma}(Q-\alpha_1)} \int_{e^{\gamma A}\frac{ \beta(|z|e^{-h})^{-2\gamma(Q-\alpha_1)} c_+J_A}{ \widehat{\mathcal{K}}_{B^c}(z)}}^\infty  \Big(  (1 + y )^{-p}  -  1\Big)\widehat{ \mathcal{K}}_{B^c}  (z)^{-p-\frac{2}{\gamma}(Q-\alpha_1)} y^{-\frac{2}{\gamma}(Q-\alpha_1)-1}\,dy\Big ].
\end{align*} 
Recalling \eqref{epsilondefinition} we have $\beta(|z|e^{-h})^{-2\gamma(Q-\alpha_1)}=e^{\gamma B_{t+h}+\gamma (Q-\alpha_1)(t+h)- \frac{\gamma^2}{2}t}= e^{ \gamma (Q-\alpha_1)h }   e^{ \gamma B_{t+h} + \gamma (Q-\alpha_1-\frac{\gamma}{2}) t  }$ and thus $\beta(|z|e^{-h})^{-2\gamma(Q-\alpha_1)}\to 0$ almost surely as $z\to 0$ (equivalently $t \to \infty$) provided $\alpha_1+\frac{\gamma}{2}>Q$ which is the case.
Dominated convergence theorem then ensures that the latter expectation converges to
\begin{align*}
& \E  [J_A^{\frac{2}{\gamma}(Q-\alpha_1)} ]\E[\widehat{ \mathcal{K}}_{B^c}(0)  ^{-p-\frac{2}{\gamma}(Q-\alpha_1)} ]\int_{0}^\infty  \big(  (1 + y )^{-p}  -  1\big)  y^{-\frac{2}{\gamma}(Q-\alpha_1)-1}\,dy  \\
& =  \frac{\Gamma(-\frac{2}{\gamma} (Q-\alpha_1)) \Gamma (p+\frac{2}{\gamma} (Q-\alpha_1) )}{\Gamma (p)}  \E  [J_A^{\frac{2}{\gamma}(Q-\alpha_1)} ]\E[\widehat{ \mathcal{K}}_{B^c}(0)  ^{-p-\frac{2}{\gamma}(Q-\alpha_1)} ], 
\end{align*}
where we have  used   Lemma \ref{lemmaintegral} in the Appendix  to compute the integral in the $y$ variable. Gathering \eqref{FirststepT2}, \eqref{SecondstepT2}, \eqref{uandlbounds} and the above considerations on the upper bound of $\E[  ( \mathcal{K}^1+  \beta c_+ \mathcal{K}^3)^{-p}]$ (and using the fact that $\xi>\gamma(Q-\alpha_1)$) we get
\begin{align*}
 & \limsup_{z\to 0}|z|^{-\gamma(Q-\alpha_1)} \Big(\E[ \mathcal{K}(z)^{-p}  -\mathcal{K}_{B^c}(z)^{-p} ]\Big)   \\
&  \leq \frac{2(Q-\alpha_1) }{\gamma}  c_+^{\frac{2}{\gamma}(Q-\alpha_1)}   \frac{\Gamma(-\frac{2}{\gamma} (Q-\alpha_1)) \Gamma (p+\frac{2}{\gamma} (Q-\alpha_1) )}{\Gamma (p)}    \E  [J_A^{\frac{2}{\gamma}(Q-\alpha_1)} ]\E[\widehat{ \mathcal{K}}_{B^c}(0)  ^{-p-\frac{2}{\gamma}(Q-\alpha_1)} ] \\
&  = \frac{1}{2} (1-e^{-h})^{\frac{2}{\gamma}(Q-\alpha_1)}  \mu^p  \Gamma (p)^{-1}    \mu^{ \frac{2}{\gamma} (Q-\alpha_1)  } \Gamma(-\frac{2}{\gamma} (Q-\alpha_1))   \frac{2(Q-\alpha_1)}{\gamma}  \E  [J_A^{\frac{2}{\gamma}(Q-\alpha_1)} ]  C_\gamma(2Q-\alpha_1-\frac{\gamma}{2}, \alpha_2,\alpha_3) .
\end{align*}
We can then conclude as for the lower bound by letting $h,A\to\infty$ since $$   \mu^{ \frac{2}{\gamma} (Q-\alpha_1)  } \Gamma(-\frac{2}{\gamma} (Q-\alpha_1))   \frac{2(Q-\alpha_1)}{\gamma}  \E  [J_A^{\frac{2}{\gamma}(Q-\alpha_1)} ]   $$ goes to $R(\alpha_1)$ as $A$ goes to infinity. 
\qed

\subsection{The 4 point function with $-\frac{2}{\gamma}$ insertion}

In this section, we prove an analogue of Theorem \ref{theo4point1} for the other degenerate insertion with weight  $-\frac{2}{\gamma}$:
\begin{theorem}\label{theo4point2overgamma}
 We assume the Seiberg bounds for $(-\frac{2}{\gamma},\alpha_1,\alpha_2,\alpha_3)$, i.e.  $\sum_{k=1}^3 \alpha_k>2Q+\frac{2}{\gamma}$ and $\alpha_k<Q$ for all $k$.  There exists  $\eta>0$ s.t. if $Q-\alpha_1<\eta$ then
\begin{equation}\label{BPZ2}
\mathcal{T}_{-\frac{2}{\gamma}} (z) = C_\gamma(\alpha_1-\frac{_2}{^\gamma}, \alpha_2,\alpha_3)  |F_{-}(z)|^2  + R(\alpha_1) C_\gamma(2Q-\alpha_1-\frac{_2}{^\gamma}, \alpha_2,\alpha_3)  |F_{+}(z)|^2.
\end{equation}
\end{theorem}

\proof The proof follows the proof of Theorem \ref{theo4point1} almost word by word and we keep the same notation with the following obvious modifications.  The function $K$ in \eqref{defineK} is replaced by
\begin{equation}\label{defineK}K(z,x)=\frac{|x-z|^2|x|_+^{\gamma (\sum_{k=1}^3\alpha_k-\frac{2}{\gamma})}}{|x|^{\gamma \alpha_1}|x-1|^{\gamma \alpha_2}}
\end{equation}
i.e. most importantly,  the factor $|x-z|^{\frac{\gamma^2}{2}}$ is replaced by $|x-z|^2$. Furthermore the exponent $p$ is now given by $p=(\alpha_1+\alpha_2+\alpha_3-\frac{2}{\gamma} -2 Q)/\gamma$ and it is positive. 

We will fix  $\eta>0$  and $\xi \in (0,1)$ so that the following conditions hold for  all $\alpha_1\in (Q-\eta,Q)$ 
\begin{align} 
\tfrac{4}{\gamma}(Q-\alpha_1)<&(1-\xi)(4-\gamma\alpha_1-2\gamma\eta)\label{defxi2}\\
\tfrac{4}{\gamma}(Q-\alpha_1)<&\xi\label{defxi2bis}.
\end{align}
Note that for $\xi=\eta=0$ \eqref{defxi2} holds since $4-\gamma Q=2-\frac{\gamma^2}{2}>0$ and therefore by continuity for small enough $\eta$ and small enough $\xi>\tfrac{4}{\gamma}\eta$ they hold as well.

As in the proof of Theorem \ref{theo4point1} we start with the splitting  \eqref{T1T2split} to $T_1$ and $T_2$ given by \eqref{T1T2def} and we first show that $T_1=o(|z|^{\frac{4}{\gamma}(Q-\alpha_1)})$. We obtain again 
$|T_1| \leq 
C(A_1+A_2)$ with the same definitions for $A_i$. 

The  Cameron-Martin bound for $A_1$ becomes
\begin{align*}
 A_1   &   \leq C \int_{|x| \leq |z|^{1-\xi}} |x|^{2-\gamma \alpha_1}\E\Big[\Big(\int_{|u|>|z|^{1-\xi}}K(0,u)|u-x|^{-\gamma^2}M_\gamma(du)\Big)^{-p-1}\Big] \,d^2x
\end{align*}
and  as the expectation is bounded by  a constant 
we conclude that
\begin{align}\label{A1bound11}A_1\leq C|z|^{(1-\xi)(4-\gamma \alpha_1)}=o(|z|^{\frac{4}{\gamma}(Q-\alpha_1)})
\end{align}
 by \eqref{defxi2}.

 Next, for $A_2$ the bound \eqref{A2bound} is replaced by
 \begin{align}\label{A2bound2} 
\E[ | \mathcal{K}_{A}(z)-\mathcal{K}_{A}(0) | \mathcal{K}_{B^c}(0) ^{-p-1}    ] 
 \leq C |z| \int_{A} |x|^{1-\gamma\alpha_1}d^2x\leq C|z|^{1+  (1-\xi)(3-\gamma\alpha_1)}=o(|z|^{\frac{4}{\gamma}(Q-\alpha_1)})
\end{align}
again by \eqref{defxi2}. Hence $T_1=o(|z|^{\frac{4}{\gamma}(Q-\alpha_1)})
$.

\bigskip
Now we proceed with $T_2$, again with the obvious changes (e.g. $\frac{\gamma^2}{2}$ in the definitions for $ \mathcal{K}^1, \mathcal{K}^2$ and $c_\pm$ replaced by 2). Hence replacing \eqref{defKhat} by
 \begin{equation*}
\widehat{ \mathcal{K}}_{D}(z):=\int_{D}\frac{|x-z|^{2}}{|x|^{\gamma(2Q-\alpha_1)}|x-1|^{\gamma\alpha_2}}|x|_+^{\gamma (2Q-\alpha_1-\frac{2}{\gamma}+\alpha_2+\alpha_3)}M_\gamma(d^2x)
\end{equation*}
 we obtain instead of \eqref{hatKbound} the bound
\begin{align}\label{hatKbound11}
\E\big[ \widehat{ \mathcal{K}}_{B^c}(z)^{-p-\frac{2}{\gamma}(Q-\alpha_1)}\big]-\E\big[ \widehat{ \mathcal{K}}_{B^c}(0)^{-p-\frac{2}{\gamma}(Q-\alpha_1)}\big]=o( |z|^{\frac{4}{\gamma}(Q-\alpha_1)}).
\end{align} 
Indeed, the LHS is  $T_1$ computed with a larger $p$ and $|x|^{\gamma\alpha_1}$ replaced by $|x|^{\gamma(2Q-\alpha_1)}$. Hence from \eqref{A1bound11} and \eqref{A2bound2} we get the bound 
\begin{align*}
\E\big[ \widehat{ \mathcal{K}}_{B^c}(z)^{-p-\frac{2}{\gamma}(Q-\alpha_1)}\big]-\E\big[ \widehat{ \mathcal{K}}_{B^c}(0)^{-p-\frac{2}{\gamma}(Q-\alpha_1)}\big]\leq C|z|^{(1-\xi)(4-\gamma (2Q-\alpha_1))}.
\end{align*} 
Since $4-\gamma (2Q-\alpha_1)=4-\gamma \alpha_1-2\gamma(Q-\alpha_1)\leq 4-\gamma \alpha_1-2\gamma\eta$, \eqref{hatKbound11} holds. The rest of the arguments for the lower and the upper bounds for $T_2$ follow then word by word. \qed

\subsection{Crossing relations}\label{sub:consequence}
Proposition \ref{propcrossing1} now follows from Theorem \ref{theo4point} as explained in Section \ref{Crossing relation}. Let us state it in the form we will apply it and also for the unit volume structure constants:
\begin{proposition}\label{propcrossing1again}
Let $\epsilon \in (\frac{\gamma}{2},\frac{2}{\gamma})$ and $\alpha,\alpha'<Q$ s.t.  $\alpha+ \alpha'+ \epsilon -\frac{\gamma}{2}>2Q$. Then  
\begin{equation} \label{Fundamentalrelation}
C_\gamma (\alpha'-\tfrac{\gamma}{2}, \epsilon, \alpha)= T(\alpha',\epsilon,\alpha)  C_\gamma (\alpha', \epsilon+\tfrac{\gamma}{2}, \alpha)
\end{equation}
where
$T$ is the given by the following formula 
\begin{equation} \label{defT0}
T(\alpha',\epsilon,\alpha)= - \mu \pi \frac{l(a) l(b)}{ l(c) l(a+b-c)}    \frac{1}{l(-\frac{\gamma^2}{4}) l(\frac{\gamma \epsilon}{2})  l(2+\frac{\gamma^2}{4}- \frac{\gamma \epsilon}{2})}
\end{equation}
where ${ a},{ b},{ c}$ are given by
\begin{align}\label{defabc}
{ a}&= \frac{_{\gamma}}{^4} (\alpha'+\alpha+\epsilon-Q- \gamma)-\hf  \quad { b}=\frac{_{\gamma}}{^4} (\alpha'-\alpha+\epsilon-Q)+\hf\quad c=1-\frac{\gamma}{2}(Q-\alpha').
\end{align}

The above relation can be rewritten under the following form for the unit volume correlations (see \eqref{unitvolumethreepoint} for the definition)
\begin{equation} \label{Fundamentalrelationanalycityunit}
\bar{C}_\gamma (\alpha'-\tfrac{\gamma}{2}, \epsilon, \alpha)= \bar{T}(\alpha',\epsilon,\alpha)  \bar{C}_\gamma (\alpha', \epsilon+\tfrac{\gamma}{2}, \alpha)
\end{equation}
where $\bar{T}$ is  given by 
\begin{equation} \label{defT}
\bar{T}(\alpha',\epsilon,\alpha) 
=  \mu^{-1}  \frac{\Gamma( \frac{1}{\gamma}(\alpha+\alpha'+{\epsilon}+\frac{\gamma}{2} -2Q))  }{\Gamma( \frac{1}{\gamma}(\alpha+\alpha'+{\epsilon}-\frac{\gamma}{2} -2Q) ) } {T}(\alpha',\epsilon,\alpha) .
\end{equation}
\end{proposition}


Along the same lines as Proposition \ref{propcrossing1again}, by exploiting Theorem \ref{theo4point2overgamma} with the $-\frac{2}{\gamma}$ insertion, one can show the two following crossing symmetry relations:

\begin{proposition}\label{2overgamma}
Let  $\alpha, \epsilon,\alpha'<Q$   with $\alpha+\alpha'+\epsilon>2Q+\frac{2}{\gamma}$. Then
 \begin{equation}\label{Firstinverserelation}
C_\gamma (\alpha-\tfrac{2}{\gamma}, \epsilon, \alpha')= \tilde T(\alpha,\epsilon,\alpha') R(\epsilon) C_\gamma (\alpha, 2Q-\epsilon-\tfrac{2}{\gamma}, \alpha')
\end{equation}
where $\tilde T$ is  given by the following formula 
\begin{equation} \label{deftildeT}
\tilde T(\alpha,\epsilon,\alpha')=  \frac{l(a) l(b)}{ l(c) l(a+b-c)}   
\end{equation}
where
\begin{align}\label{defab}
{ a}&= \tfrac{1}{\gamma} (\alpha'+\alpha+\epsilon-Q-\tfrac{4}{\gamma})-\hf  \quad { b}= \tfrac{1}{\gamma} (\alpha-\alpha'+\epsilon-Q)+\hf\quad c=1-\tfrac{2}{\gamma}(Q-\alpha).
\end{align}
\end{proposition}

{
\begin{proposition}\label{2overgammabis}
Let  $\alpha, \epsilon,\alpha'<Q$   with $\alpha+\alpha'+\epsilon>2Q+\frac{2}{\gamma}$. Then
\begin{equation} 
R(\epsilon) C_\gamma(2Q-\epsilon-\tfrac{2}{\gamma}, \alpha,\alpha' )=  
L(\epsilon,\alpha,\alpha')R(\alpha) C_\gamma(\epsilon, 2Q-\alpha-\tfrac{2}{\gamma},\alpha' )
\end{equation}
where $L$ is  given by the following formula 
\begin{equation}\label{Thekeyequationell}
L(\epsilon,\alpha,\alpha')=   \frac{l(c-1) l(c-a-b+1)}{l(c-a)l(c-b)} 
\end{equation}
 with
\begin{align} 
{ a}&= \tfrac{1}{\gamma} (\alpha'+\alpha+\epsilon-Q-\tfrac{4}{\gamma})-\hf  \quad { b}= \tfrac{1}{\gamma} (\alpha-\alpha'+\epsilon-Q)+\hf\quad c=1-\tfrac{2}{\gamma}(Q-\epsilon).
\end{align}
\end{proposition}}

\section{Proof of  
Theorem \ref{Rtheor}}\label{proofreflection}

We will suppose that $\gamma^2 \not \in \Q$. This is no restriction since the general case can be deduced from this case by continuity in $\gamma$ (Remark \ref{analingamma}). The proof of formula \eqref{keyformula} for the reflection coefficient is made of several steps as explained in Section \ref{sec:bone}:


\vspace{0.2 cm}

\noindent {\it Subsection \ref{Step1}:} We prove that $\bar{R}$ is analytic in a complex neighborhood of the  interval $(\frac{\gamma}{2},Q)$. The key to this is the crossing relation  \eqref{Fundamentalrelationanalycityunit} that allows to express $\bar{R}(\alpha)$ in terms of $C_\gamma(\alpha,\gamma,\alpha)$ (eq. \eqref{formulaforbarR}).


\vspace{0.2 cm}

\noindent {\it 
Subsection \ref{Step2}:} We prove first that $R$ satisfies the following shift equation for $\alpha$ close to $ Q$: 
\begin{equation}\label{shift1}
R(\alpha-\frac{\gamma}{2})= - \mu \pi \frac{R(\alpha)}{ l(-\frac{\gamma^2}{4}) l(\frac{\gamma\alpha}{2}-\frac{\gamma^2}{4})  l(2+\frac{\gamma^2}{4}- \frac{\gamma \alpha}{2})}.
\end{equation}
The starting point is again the crossing relation  \eqref{Fundamentalrelationanalycityunit}. Using the tail estimate Lemma \ref{doubletail} we  show that the RHS of \eqref{Fundamentalrelationanalycityunit} has two simple poles in $\epsilon$ and \eqref{shift1} follows by equating residues of both sides of  \eqref{Fundamentalrelationanalycityunit}. Next, by analyticity the relation \eqref{shift1} 
extends to a neighborhood of $\alpha \in (\gamma,Q)$. Analyticity of $\bar R$ on $(\frac{\gamma}{2},Q)$ then implies we can use  \eqref{shift1}  to extend $R$ to a neighborhood of $\R$. The extension that we also denote by $R$ is meromorphic with simple poles on the real line located at $\lbrace \frac{2}{\gamma}- \frac{\gamma}{2} \N \rbrace \cup \lbrace \frac{\gamma}{2}- \frac{2}{\gamma} \N \rbrace $.  

\vspace{0.2 cm}

\noindent {\it 
Subsection \ref{Gluing}:} We prove  the so-called gluing lemma, Lemma \ref{gluinglemma}, that uses $R$ to extend the three point structure constant to a holomorphic function in a neighborhood of $Q$. The basic input in the proof is the shift relation  \eqref{3pointconstanteqintro} proven in Theorem \ref{theo4point} and Corollary \ref{shiftcoro}, based on Theorem \ref{theo4point1}. 
\vspace{0.2 cm}

\noindent {\it 
Subsection \ref{Step3}}: We prove that $R$ satisfies the following inversion relation:
\begin{equation}\label{inversion}
R(\alpha)R(2Q-\alpha)=1.
\end{equation}
The proof is based on combining  the crossing relation Proposition \ref{2overgamma} with the gluing lemma. 

\vspace{0.2 cm}

\noindent {\it 
Subsection \ref{Step4}}:  We prove that $R$ (as a meromorphic function in a neighborhood of  $\R$) satisfies the following shift equation
\begin{equation}\label{shift2}
R(\alpha)=   - c_\gamma \frac{R(\alpha+\frac{2}{\gamma})}{l(-\frac{4}{\gamma^2})l(\frac{2 \alpha}{\gamma}) l(2+\frac{4}{\gamma^2}-\frac{2 \alpha}{\gamma})} 
\end{equation}
where $c_\gamma = \frac{\gamma^2}{4} \mu \pi R(\gamma) \not = 0$. 
Recall that from the DOZZ solution we expect that 
\begin{equation}\label{shift2a} c_\gamma={(\mu \pi l(\tfrac{\gamma^2}{4})  )^{\frac{4}{\gamma^2} }}{l(\tfrac{4}{\gamma^2})}^{-1}.
\end{equation}

\vspace{0.2 cm}

\noindent {\it 
Subsection \ref{finalproofRDOZZ}}: since $R^{{\rm DOZZ}}$ satisfies \eqref{shift1} and \eqref{shift2} with \eqref{shift2a}, we prove $R=R^{{\rm DOZZ}}$ by application of Liouville's theorem. 


\subsection{Proof of analyticity of $\bar{R}$ in the interval $(\frac{\gamma}{2},Q)$} \label{Step1}
The crossing relation  \eqref{Fundamentalrelationanalycityunit} gives for  $\alpha=\alpha'$ 
\begin{equation}\label{ssh}
\bar{C}_\gamma (\alpha-\tfrac{\gamma}{2}, \epsilon, \alpha)= \bar{T}(\alpha,\epsilon,\alpha)  \bar{C}_\gamma (\alpha, \epsilon+\tfrac{\gamma}{2}, \alpha)
\end{equation}
which holds for $\alpha<Q$ and $\epsilon \in (\frac{\gamma}{2},\frac{2}{\gamma})$ with $ \epsilon +\frac{\gamma}{2}<2 \alpha$.
From Remark \ref{unitvolrel}  we deduce 
 for  $\alpha\in (\frac{\gamma}{2}, Q)$

\begin{equation*}
\lim_{\epsilon\downarrow  \frac{\gamma}{2}}(\epsilon-\frac{_\gamma}{^2})\bar{C}_\gamma (\alpha-\tfrac{\gamma}{2}, \epsilon, \alpha)  =  
\tfrac{4(Q-\alpha) }{\gamma}  \bar{R}(\alpha).
\end{equation*}

By Theorem \ref{analytic_hyperbolic}, for $\epsilon>\frac{\gamma}{2}$, $\bar{C}_\gamma (\alpha-\frac{\gamma}{2}, \epsilon, \alpha)$ is analytic in  $\alpha\in (\frac{\gamma}{2}, Q)$ and $ \bar{C}_\gamma (\alpha, \epsilon+\frac{\gamma}{2}, \alpha)
$ is analytic in $\alpha\in (\frac{\gamma}{4}+ \frac{\epsilon}{2},Q)$. Hence the relation \eqref{ssh} holds for $\epsilon \in (\frac{\gamma}{2},\frac{2}{\gamma})$   and $\alpha\in (\frac{\gamma}{4}+ \frac{\epsilon}{2}, Q)$. Using \eqref{defT0} and \eqref{defT} a bit of  calculation gives

\begin{equation*}
\lim_{\epsilon\downarrow  \frac{\gamma}{2}}(\epsilon-\frac{_\gamma}{^2})\bar{T}(\alpha,\epsilon,\alpha) =
- \pi \tfrac{4(Q-\alpha)}{\gamma^2} 
 \frac{   l(\frac{\gamma}{2} \alpha -\frac{\gamma^2}{4}-1)}{   l(1+\frac{\gamma}{2}(\alpha-Q)) l(-\frac{\gamma^2}{4}) l(\frac{\gamma^2}{4})     }  \end{equation*}
we conclude that for   all $\alpha\in (\frac{\gamma}{2}, Q)$ 
\begin{equation}\label{formulaforbarR}
\bar{R}(\alpha)= -   \frac{_\pi }{^\gamma}   \frac{   l(\frac{\gamma}{2} \alpha -\frac{\gamma^2}{4}-1)}{   l(1+\frac{\gamma}{2}(\alpha-Q)) l(-\frac{\gamma^2}{4}) l(\frac{\gamma^2}{4})     }  
\bar{C}_\gamma (\alpha,\gamma, \alpha)
\end{equation}
which proves our claim since $\bar{C}_\gamma (\alpha,\gamma, \alpha)$ is analytic  in  $\alpha\in (\frac{\gamma}{2}, Q)$.

\subsection{Proof of the $\frac{\gamma}{2}$ shift equation \eqref{shift1}}\label{Step2}  We start again with the crossing relation \eqref{Fundamentalrelationanalycityunit} 
\begin{equation} \label{barFundamentalrelation}
\bar C_\gamma (\alpha'-\tfrac{\gamma}{2}, \epsilon, \alpha)= \bar T(\alpha',\epsilon,\alpha)  \bar C_\gamma (\alpha', \epsilon+\tfrac{\gamma}{2}, \alpha)
\end{equation}
which holds for $\epsilon \in  (\frac{\gamma}{2}, \frac{2}{\gamma})$  and  $\alpha+\alpha'+\epsilon-\frac{\gamma}{2}>2Q$ with $\alpha,\alpha'<Q$. 

By Theorem  \ref{analytic_hyperbolic}, both sides of \eqref{barFundamentalrelation} are restrictions of holomorphic functions over an open neighborhood of the intersection of the extended Seiberg bounds \eqref{ThextendedSeibergbounds} valid for each structure constant involved in each side of \eqref{barFundamentalrelation}, which thus remains valid on this set. It is rather tedious to write explicitly this set but one can check that it contains the set of values
\begin{equation} \label{alphaalpha'}
\alpha'=Q-\eta,  \ \ \alpha=\tfrac{2}{\gamma}+\eta,  \ \ \epsilon \in (2\eta,\frac{2}{\gamma})
\end{equation} 
for any $\eta\in (0,\frac{\gamma}{4})$.
 
Let us consider  both sides of \eqref{barFundamentalrelation} as a function of   $\epsilon$. From Remark \ref{unitvolrel} we obtain
\begin{equation*}
\lim_{\epsilon\downarrow 2\eta}(\epsilon-2\eta) \bar C_\gamma (\alpha'-\tfrac{\gamma}{2}, \epsilon, \alpha)=
\tfrac{4(Q-\alpha) }{\gamma}  \bar{R}(\alpha).
\end{equation*}
This indicates that $\epsilon \mapsto\bar C_\gamma (\alpha'-\frac{\gamma}{2}, \epsilon, \alpha)$ has a pole at $\epsilon=2\eta$. The extended Seiberg bounds indicates that the next pole below $\epsilon=2\eta$ is located either at $\epsilon=-2\eta$ or $\epsilon = \gamma-\frac{4}{\gamma}$. We reinforce the restriction on $\eta$ to be
\begin{equation}\label{def:etasmall}
0<\eta<(\tfrac{2}{\gamma}-\tfrac{\gamma}{2})\wedge\tfrac{\gamma}{4}\wedge \tfrac{1}{2\gamma}
\end{equation}
in order to make sure  that the next  pole is at  $\epsilon=-2\eta$ (the condition $\eta< \tfrac{1}{2\gamma}$ is just technical and makes sure the interval $(2\eta,\tfrac{1}{\gamma})$ is non empty in the argument just below). Indeed we prove:

\begin{proposition}\label{poleresult}
Let $\alpha,\alpha'$ be given by \eqref{alphaalpha'}. Then for $\eta>0$ small enough the function 
\begin{equation*}
f(\epsilon):=\bar{C}_\gamma(\alpha'-\tfrac{\gamma}{2}, \epsilon, \alpha)- \frac{\tfrac{4}{\gamma}(Q-\alpha)\bar R(\alpha)}{\epsilon-2\eta}-\frac{\tfrac{4}{\gamma}(Q-\alpha'+\tfrac{\gamma}{2})\bar R(\alpha'-\tfrac{\gamma}{2})}{\epsilon+2\eta}
\end{equation*}
extends to an analytic function in a complex neighborhood of  $\epsilon \in (-2\eta-\delta,\tfrac{1}{\gamma})$ for some $\delta>0$. 
\end{proposition}

We postpone the proof of Proposition \ref{poleresult} to the end of this subsection. 
By \eqref{barFundamentalrelation} for $\epsilon \in (2\eta,\tfrac{1}{\gamma})$ we have $f(\epsilon)=g(\epsilon) $ where
\begin{equation*} 
g(\epsilon):=\bar T(\alpha',\epsilon, \alpha)  \bar C_\gamma (\alpha', \epsilon+\tfrac{\gamma}{2}, \alpha)- \frac{\tfrac{4}{\gamma}(Q-\alpha)\bar R(\alpha)}{\epsilon-2\eta}-\frac{\tfrac{4}{\gamma}(Q-\alpha'+\tfrac{\gamma}{2})\bar R(\alpha'-\tfrac{\gamma}{2})}{\epsilon+2\eta}.
\end{equation*}
Thus by analytic continuation of $f$ obtained above, $g$ is analytic in $\epsilon$ on $(-2\eta-\delta,\tfrac{1}{\gamma})$. By Remark \ref{unitvolrel}
\begin{equation*}
\lim_{\epsilon \downarrow -2\eta} (\epsilon +2\eta)\bar C_\gamma (\alpha', \epsilon+\tfrac{\gamma}{2}, \alpha)=\tfrac{4(Q-\alpha') }{\gamma} \bar R(\alpha')
\end{equation*}
where we used $\alpha'>\alpha$. Hence, from $\lim_{\epsilon \downarrow -2\eta} (\epsilon +2\eta)g(\epsilon)=0$ we deduce
\begin{equation*}
(Q-\alpha') \bar T(\alpha',-2\eta, \alpha)  \bar R(\alpha')=(Q-\alpha'+\tfrac{\gamma}{2})\bar R(\alpha'-\tfrac{\gamma}{2}).
\end{equation*}
This is the reflection relation for unit volume reflection coefficient. Using 
 \eqref{defT} and \eqref{deffullR} a calculation gives then
\begin{equation}\label{equiv222}
 R(\alpha'-\tfrac{\gamma}{2})=T(\alpha',-2\eta, \alpha)   R(\alpha').
\end{equation}

Inserting to  \eqref{defabc}  $ \epsilon=- 2 \eta= \alpha'-\alpha- \frac{\gamma}{2}$ we get first that ${ b}
= \frac{\gamma \epsilon}{2}$ so that $l(b)=l( \frac{\gamma \epsilon}{2})$ and $ a+b-c= 1- (   2+ \frac{\gamma^2}{4} -\frac{\gamma \epsilon}{2}  )$ so that  $l(a+b-c) l(2+ \frac{\gamma^2}{4} -\frac{\gamma \epsilon}{2} ) = 1$
Therefore \eqref{defT0} becomes \begin{equation*}
T(\alpha',-2\eta,\alpha) = - \mu \pi   \frac{ l( a)  }{  l( c  ) l(-\frac{\gamma^2}{4})   }   = - \mu \pi   \frac{ l(  \tfrac{\gamma\alpha'}{2}-\tfrac{\gamma^2}{2}-1)  }{  l( 1+ \tfrac{\gamma}{2} (\alpha'-Q)  ) l(-\tfrac{\gamma^2}{4})   }  .
\end{equation*}
Using $l(x)^{-1}=l(1-x)$ \eqref{equiv222}  is the desired shift relation \eqref{shift1} (with $\alpha'=\alpha$).

We have proven  \eqref{shift1}  for $\alpha$ close to $Q$ but since by previous subsection $\bar R$ is analytic on $(\frac{\gamma}{2},Q)$ the relation \eqref{shift1} extends to $\alpha \in (\gamma,Q)$. Then we can use \eqref{shift1} to extend $R$  to a meromorphic function in a neighbourhood of $\R$ which we also denote by  $R$. Since $R^{{\rm DOZZ}}$ also satisfies \eqref{shift1}  and $0<\frac{R(\alpha)}{R^{{\rm DOZZ}}(\alpha)}< \infty$ for  $\alpha \in (\frac{\gamma}{2},Q)$ we conclude that $R$ and $R^{{\rm DOZZ}}$ have their poles and zeros located at the same places. For instance, the poles of $R$ are located at $\lbrace \frac{2}{\gamma}- \frac{\gamma}{2} \N \rbrace \cup \lbrace \frac{\gamma}{2}- \frac{2}{\gamma} \N \rbrace $.  

A useful consequence of this analytic continuation of $R$ is the following
\begin{corollary}\label{shiftcoro}
Let $\sum_{k=1}^3 \alpha_k-\frac{\gamma}{2}>2Q$ and $\alpha_k<Q$ for all $k$. There exists $\eta>0$ s.t. if $Q-\alpha_1<\eta$ then
The shift equation  \eqref{3pointconstanteqintro} holds in the form  
\begin{equation*}
R({\alpha}_1+\tfrac{\gamma}{2})C_\gamma(2Q-{\alpha}_1-\tfrac{\gamma}{2},{\alpha}_2,{\alpha}_3)
= - \frac{1}{\pi \mu}\caA(\tfrac{\gamma}{2},{\alpha}_1, \alpha_2,{\alpha}_3)C_\gamma({\alpha}_1-\tfrac{\gamma}{2},{\alpha}_2,{\alpha}_3).
\end{equation*} 
\end{corollary}
\proof 
From Theorem \ref{theo4point1}, we get
\begin{equation}\label{shiftcoro1}
\mathcal{T}_{-\frac{\gamma}{2}}(z) = C_\gamma(\alpha_1-\frac{\gamma}{2}, \alpha_2,\alpha_3)  |F_{-}(z)|^2  + R(\alpha_1) C_\gamma(2Q-\alpha_1-\frac{\gamma}{2}, \alpha_2,\alpha_3)  |F_{+}(z)|^2.
\end{equation}
Thanks to  \eqref{shift1} applied to $\alpha=\alpha_1+\tfrac{\gamma}{2}$, we get that 
\begin{equation}\label{shiftcoro2}
R(\alpha_1)= - \mu \pi \frac{R(\alpha_1+\frac{\gamma}{2})}{ l(-\frac{\gamma^2}{4}) l(\frac{\gamma \alpha_1}{2})  l(2- \frac{\gamma \alpha_1}{2})}.
\end{equation}
Plugging \eqref{shiftcoro2} into \eqref{shiftcoro1}, the result then follows from \eqref{Tsolution} and  \eqref{Fundrelation} applied to $\alpha_0=-\frac{\gamma}{2},\alpha_1,\alpha_2,\alpha_3$ and a lengthy calculation.\qed
\vskip 3mm

\noindent{\it Proof of Proposition \ref{poleresult}}
Fix points  $z_2,z_3\in\C$ such that $|z_2| \geq 2, |z_3| \geq 2$ and $|z_2-z_3| \geq 3$. From \eqref{expression3pointstruct} and \eqref{unitvolumethreepoint}, we have
\begin{equation*}
\bar C_\gamma (\alpha'-\tfrac{\gamma}{2}, \epsilon, \alpha)=G(\epsilon)\ \E\, Z_\C(\epsilon)^{1 -\frac{\epsilon}{\gamma}}
\end{equation*}
where
\begin{equation}\label{defgepsi}
G(\epsilon) = 2  
 \gamma^{-1}  \prod_{i<j}  \frac{1}{|z_i-z_j|^{\alpha_i \alpha_j+2 \Delta_{ij}}}
\end{equation}
with $z_1=0$, $(\alpha_1,\alpha_2,\alpha_3)=  (\epsilon, \alpha, \alpha'-\frac{\gamma}{2} )=(\epsilon, \frac{2}{\gamma}+\eta, \frac{2}{\gamma}-\eta)$ and for $A\subset\C$ 
\begin{equation}\label{defZ}
Z_A(\epsilon):=  \int_{A } \frac{|x|_+^{  \gamma(\epsilon+\alpha_2+\alpha_3) }}{|x|^{\gamma\epsilon}  |x-z_2|^{\gamma \alpha_2}|x-z_3|^{\gamma \alpha_3}}   M_\gamma(d^2x).
\end{equation} 
Define next
\begin{equation*}
F(\epsilon):=\E\big(Z_\C(\epsilon)^{1 -\frac{\epsilon}{\gamma}}- (Z_{B_1(z_2)}(\epsilon)+Z_{B_1(z_3)}(\epsilon) )^{1 -\frac{\epsilon}{\gamma}}\big).
\end{equation*}
The fact that $F$ is well defined is a consequence of the proof of Lemma \ref{analy1} below. Note that $Z_{B_1(z_2)}(\epsilon)$ and $Z_{B_1(z_3)}(\epsilon)$ do not depend on $\epsilon$ since for $x\in  B_1(z_2)$ or $x\in  B_1(z_3)$ we have $
|x|_+=|x|$. Hence we denote them by $Z_{B_1(z_2)}$ \textcolor{red}{and $Z_{B_1(z_3)}$}.
We start with

\begin{lemma}\label{analy1}
$F(\epsilon)$ is analytic in   a complex neighborhood of the interval $ (-2 \eta-\delta, \tfrac{1}{\gamma})$  for some $\delta>0$.
\end{lemma}

\proof Let us fix $\delta>0$ such that
\begin{equation}\label{def:delta}
2\eta+\delta<   \tfrac{4}{\gamma}-\gamma \quad \text{ and }\quad 4\eta+\delta<\gamma
\end{equation}
which is possible because of  \eqref{def:etasmall}. As in the proof of Theorem \ref{analytic_hyperbolic} we construct  $F$  as the uniform limit as $t\to\infty$  of analytic functions $F_t$ in a neighborhood of $(-2 \eta-\delta, \tfrac{1}{\gamma})$.  Let us denote $\C_t=\{z: |z|\geq e^{-t}\}$ and define
 (recall that $B_r$ stands for the ball $B_r(0)$ and $X_t$ for the circle average \eqref{circleav} with $r=e^{-t}$)
\begin{equation*}
F_t(\epsilon)= 
 \E  \left [    e^{ \epsilon X_{t} (0)-\frac{t\epsilon^2}{2} }   \left (   Z_{\C_t}(0)^{1 -\frac{\epsilon}{\gamma}}
 - (Z_{B_1(z_2)}+Z_{B_1(z_3)} )^{1 -\frac{\epsilon}{\gamma}}
  \right )     \right   ].
\end{equation*}

Let us first show that for each $t$, $\epsilon\mapsto F_t(\epsilon)$ is an analytic function of $\epsilon$ in an open neighborhood of $(-2 \eta-\delta, \tfrac{1}{\gamma})$. 
Let $R_1:=Z_{B_1(z_2)}+Z_{B_1(z_3)} $ and $R_2:=Z_{\C_t}(0)-R_1$.
By \eqref{ThextendedSeibergbounds0} and  \eqref{ThextendedSeibergbounds} $R_1$ admits moments of order $q$ for $q<\frac{2}{\gamma}(Q-\eta-\frac{2}{\gamma})$ and $R_2$ has moments of order $q$ for $q<\frac{4}{\gamma^2}$. By taking the derivative of $s \mapsto   (sR_2+R_1)^{1 -\frac{\epsilon}{\gamma}}$ we get 
\begin{equation} \label{interpol2}
\E[e^{ \epsilon X_t (0) } \big((R_2+R_1)^{{1 -\frac{\epsilon}{\gamma}}}-R_1^{{1 -\frac{\epsilon}{\gamma}}}\big)]=(1 -\frac{\epsilon}{\gamma})\int_0^1\E[e^{ \epsilon X_t(0) }R_2 (s R_2+R_1)^{{ -\frac{\epsilon}{\gamma}}}]\,ds.
\end{equation}
Let $\epsilon=\epsilon_1+i\epsilon_2$ and suppose first $\epsilon_1>0$. Since $\E |e^{p \epsilon X_t (0) }|<\infty $ for all $p<\infty$ and since chaos has negative moments  then using H\"older's inequality we can bound the integrand by $C\E[R_2^q]^{\frac{1}{q}}$ for any $q>1$. 

If $\epsilon_1<0$ we bound 
$$|\E[R_2 (s R_2+R_1)^{{ -\frac{\epsilon}{\gamma}}}]|\leq C(\E[R_2 ^{1-\frac{\epsilon_1}{\gamma}}]+\E [R_2R_1^{ -\frac{\epsilon_1}{\gamma}}]).
$$
This is finite provided  $1-\frac{\epsilon_1}{\gamma}<\tfrac{4}{\gamma^2}$ and  $ -\frac{\epsilon_1}{\gamma}<\frac{2}{\gamma}(Q-\eta-\frac{2}{\gamma})=1-\frac{2\eta}{\gamma}$ (by a slight variant of Remark \ref{doubletailmoments}). These conditions hold due to \eqref{def:delta}. Using similar bounds, one can show that the derivative (with respect to $\epsilon$) of the right hand side of \eqref{interpol2} exists and $ F_t$ is thus seen to be complex differentiable in $\epsilon$.  

\medskip
Next we show that the family $F_t$ is Cauchy in the topology of uniform convergence over compact subsets of a neighborhood of the interval $ (-2 \eta-\delta, \tfrac{1}{\gamma})$. For this we will  bound $F_{t+1}-F_t$.   First observe that because $Z_{B_1(z_2)}(\epsilon)$ and $Z_{B_1(z_3)}(\epsilon)$ are independent of $X_t(0) $ (see Remark \eqref{indebm}) these terms cancel out in
 $F_{t+1}-F_t$. Furthermore, Girsanov theorem gives 
$$  \E   e^{ \epsilon X_{t} (0)-\frac{t\epsilon^2}{2} }    Z_{\C_t}(0)^{1 -\frac{\epsilon}{\gamma}}= \E   e^{ i\epsilon_2 X_{t} (0)+\frac{t\epsilon_2^2}{2} }     Z_{\C_t}(\epsilon_1)^{1 -\frac{\epsilon}{\gamma}}.
 $$
 Hence as  in the proof of Theorem \ref{analytic_hyperbolic} we get 
 $$|F_{t+1}-F_t|\leq e^{\frac{(t+1)\epsilon_2^2}{2} } \E   | Z_{\C_{t+1}}(\epsilon_1)^{1 -\frac{\epsilon}{\gamma}}-Z_{\C_t}(\epsilon_1)^{1 -\frac{\epsilon}{\gamma}}|.$$
 
From now on, since $\epsilon_1$ is fixed  we suppress it in the notation and denote  $Z_{\C_t}(\epsilon_1)$ by $Z_{t}$. We proceed as in the proof of Theorem \ref{analytic_hyperbolic}. Let $Y_t:=Z_{t+1}-Z_{t}$. We fix $\theta>0$ and write
\begin{align*}
 \E   | Z_{t+1}^{1 -\frac{\epsilon}{\gamma}}-Z_{t}^{1 -\frac{\epsilon}{\gamma}}|
\leq  \E 1_{ Y_t \leq e^{-\theta t}}| Z_{t+1}^{1 -\frac{\epsilon}{\gamma}}-Z_{t}^{1 -\frac{\epsilon}{\gamma}}|+ \E  1_{ Y_t \geq e^{-\theta t} }| Z_{t+1}^{1 -\frac{\epsilon}{\gamma}}-Z_{t}^{1 -\frac{\epsilon}{\gamma}}|.
\end{align*}
Interpolating the first term is bounded by
 \begin{align*}
   \E \ind_{Y_t<e^{-\theta t}}|  (Z_t+Y_t) ^{1 -\frac{\epsilon}{\gamma}} -Z_t^{1 -\frac{\epsilon}{\gamma}} |\leq Ce^{-\theta t}\sup_{s\in[0,1]}\E(Z_t+sY_t) ^{ -\frac{\epsilon_1}{\gamma}}\leq Ce^{-\theta t}
  \E(Z_\C(\epsilon_1)^{ -\frac{\epsilon_1}{\gamma}}    )  .
\end{align*}
The last expectation is finite since $-\frac{\epsilon_1}{\gamma}<\frac{2}{\gamma}(Q-\eta-\frac{2}{\gamma})=1-\frac{2\eta}{\gamma}$ holds by \eqref{def:delta}.

For the second term we use in turn H\"older's inequality, with $p,q>1$ such that $\tfrac{1}{p}+\tfrac{1}{q}=1$, and the mean value theorem to get
\begin{align*}
 \E  1_{ Y_t \geq e^{-\theta t} }| Z_{t+1}^{1 -\frac{\epsilon}{\gamma}}-Z_{t}^{1 -\frac{\epsilon}{\gamma}}|&\leq [\P\big( Y_t \geq e^{-\theta t}\big)]^{1/p} [\E | Z_{t+1}^{1 -\frac{\epsilon}{\gamma}}-Z_{t}^{1 -\frac{\epsilon}{\gamma}}|^q]^{\frac{1}{q}}\\&
 \leq [\P\big( Y_t \geq e^{-\theta t}\big)]^{1/p} \sup_{s\in[0,1]}[\E Y_t^q(Z_t+sY_t) ^{ -q\frac{\epsilon_1}{\gamma}}]^{\frac{1}{q}}.
 \end{align*}
By the Markov inequality, the definition $Y_t=Z_{\C_{t+1}\setminus\C_{t}}$ and  the chaos moment estimate \eqref{annulus}  we get
 \begin{align*}
 \P\big( Y_t \geq e^{-\theta t}\big)]^{1/p}
 \leq  
 e^{-\frac{\theta m}{p}t}\E[Y_t^m]^{1/p}\leq 
 e^{\frac{1}{p}({\gamma(Q-\epsilon_1-\theta)m-\tfrac{m^2\gamma^2}{2}})}
\end{align*}
so that  we end up with the bound
\begin{equation}\label{estF}
|F_{t+1}-F_t|\leq Ce^{\frac{t\epsilon_2^2}{2} } (e^{-\theta t}+C_t(q)e^{\frac{1}{p}({\gamma(Q-\epsilon_1-\theta)m-\tfrac{m^2\gamma^2}{2}})})
\end{equation}
where we defined 
\begin{equation}\label{estF11}C_t(q)=\sup_{s\in[0,1]}[\E Y_t^q(Z_t+sY_t) ^{ -q\frac{\epsilon_1}{\gamma}}]^{\frac{1}{q}}.
\end{equation}

Now we have to optimize with respect to the free parameters $p,q,\theta,m$. We first fix  $q$ (hence $p$) to make   \eqref{estF11} finite. Let first $\epsilon_1>0$. By existence of negative moments of chaos we get for all $r>q$
\begin{equation*}
C_t(q)\leq  C(r)[\E Y_t^r]^{\frac{1}{r}}.
\end{equation*}
Hence by the chaos moment estimate \eqref{annulus} $\sup_tC_t(q)<\infty$ if $q<\frac{2}{\gamma}(Q-\epsilon_1)\wedge \frac{4}{\gamma^2}= \frac{4}{\gamma^2}$. 

If  $\epsilon_1<0$ we bound $Y_t\leq Z_{B_1}(\epsilon_1)$ and $Z_{t+1}\leq  Z_{B_1}(\epsilon_1)+ Z_{B_1^c}(\epsilon_1)$ to get
\begin{equation*}
C_t(q)\leq  [\E Y_t^qZ_{t+1}^{ -q\frac{\epsilon_1}{\gamma}}]^{\frac{1}{q}}\leq C [\E Z_{B_1}(\epsilon_1)^{q(1-\frac{\epsilon_1}{\gamma})}+\E Z_{B_1}(\epsilon_1)^{q}Z_{B_1^c}(\epsilon_1)^{ -q\frac{\epsilon_1}{\gamma}}]^{\frac{1}{q}}.
\end{equation*}
The first expectation is finite if  $q\big(1+\frac{2\eta+\delta}{\gamma}\big)<\frac{4}{\gamma^2}$ 
and  by Remark \ref{doubletailmoments} the second one is finite if  $q\frac{2\eta+\delta}{\gamma}<\frac{2}{\gamma}(Q-\eta-\frac{2}{\gamma})=1-\frac{2\eta}{\gamma}$. Due to \eqref{def:delta} we can find $q>1$ such that 
 this condition holds and hence $\sup_tC_t(q)<\infty$.

Next, we choose  $\theta>0$ such that $Q-\frac{2}{\gamma}-\theta>0$ and then  $m\in (0,1)$ small enough such that
$$\kappa:= p^{-1} ({\gamma(Q-\frac{_2}{^\gamma}-\theta)m-\tfrac{\gamma^2}{2}m^2})>0.
$$
 As we have $\epsilon_1<\frac{2}{\gamma}$  we get from \eqref{estF}
\begin{equation}\label{estF1}
|F_{t+1}-F_t|\leq Ce^{\frac{t\epsilon_2^2}{2} } (e^{-\theta t}+e^{-\kappa t}).
\end{equation}
 Hence 
  the sequence $F_t$ converges uniformly in compacts of    a neighborhood of $(-2 \eta-\delta, \tfrac{1}{\gamma})$. Finally  observe that $F(\epsilon)=\lim_{t\to \infty}F_t(\epsilon)$ for $\epsilon \in \R$ with $\epsilon \in (2\eta,Q)$. \qed

\vskip 2mm

Next we note the following simple lemma on analytic continuation of moments of random variables:

\begin{lemma}\label{analy2a} Let  $Y\geq 0$ be a random variable with a tail estimate 
\begin{equation}\label{eq:analy2a}
|\P(Y>t)-c_1t^{-\beta_1}-c_2t^{-\beta_2}|\leq c_3t^{-\beta_3}
\end{equation} 
 for some  constants $c_1,c_2,c_3> 0$  and $0<\beta_1<\beta_2<\beta_3$.  Then $s\in \C\mapsto\E Y^s$ extends to a 
meromorphic function in the strip $0<\Re s<\beta_3$  given by
$$\E[Y^s]= \frac{c_1s}{\beta_1-s}+\frac{c_2s}{\beta_2-s}+r(s)$$
where $r$ is holomorphic in $0<\Re s<\beta_3$.
\end{lemma}
\proof 
Since 
$$\E[Y^s]=s\int_0^\infty\P(Y>t)t^{s-1}\,dt
$$
and $\P(Y>t)\leq Ct^{-\beta_1}$ the mapping $s\in \C\mapsto \E[Y^s]$ is holomorphic on the set $\{s\in\C; 0<\Re(s)<\beta_1\}$.  Writing 
$$\E[Y^s]=s\int_0^\infty\big(\P(Y>t)-(c_1t^{-\beta_1}+c_2t^{-\beta_2})\mathbf{1}_{\{t\geq 1\}}\big)t^{s-1}\,dt -\frac{sc_1}{s-\beta_1}-\frac{sc_2}{s-\beta_2}$$
the claim follows as the first term on the RHS is holomorphic on the set $\{s\in\C; 0<\Re(s)<\beta_3\}$ due to the assumption \eqref{eq:analy2a}. \qed
 
\vskip 2mm

We apply this lemma to the study of the random variable $Y=Z_{B_1(z_2)}+Z_{B_1(z_3)}$. We use the tail estimate  Lemma \ref{doubletail} for $\eta$ small enough so that $\alpha_2,\alpha_3$ are both sufficiently close to each other (recall Remark \eqref{remarkdoubletail}). We have $\beta_1=1-\frac{2\eta}{\gamma}$, $\beta_2=1+\frac{2\eta}{\gamma}$ and some calculation gives
$$
c_1=|z_2|^{-4\eta Q}|z_2-z_3|^{-2(1+Q\eta+\eta^2)}\bar R(\alpha_2),\ \ c_2=|z_3|^{4\eta Q}|z_2-z_3|^{-2(1-Q\eta+\eta^2)}\bar R(\alpha_3).
$$
Then Lemma \ref{analy2a} gives
\begin{equation}\label{analy2aa}
\E(Z_{B_1(z_2)}+Z_{B_1(z_3)})^{1-\frac{\epsilon}{\gamma}}= \frac{(\gamma-\epsilon)c_1}{\epsilon-2\eta}+\frac{(\gamma-\epsilon)c_2}{\epsilon+2\eta}+r(\epsilon)
\end{equation} 
 where $r$ is analytic in a complex neighborhood of $(-2\eta-\delta, \frac{1}{\gamma})$. Since $G(\epsilon)$ in \eqref{defgepsi} is analytic in this region too we conclude by combining Lemma \ref{analy1} and \eqref{analy2aa} that
\begin{equation*}
\bar C_\gamma (\alpha'-\tfrac{\gamma}{2}, \epsilon, \alpha)= \frac{a_1}{\epsilon-2\eta}+\frac{a_2}{2\eta+\epsilon}+f(\epsilon)
\end{equation*} 
with $a_1=G(2\eta)(\gamma-\epsilon)c_1=\frac{2}{\gamma}(\gamma-2\eta)\bar R(\alpha_2)=\frac{4}{\gamma}(Q-\alpha_2)\bar R(\alpha_2)$ (note that the $z_2,z_3$ dependence has to cancel!) and  $a_2=G(-2\eta)(\gamma-\epsilon)c_2=\frac{4}{\gamma}(Q-\alpha_3)\bar R(\alpha_3)$. $f$ is analytic in a complex neighborhood of $(-2\eta-\delta, \frac{1}{\gamma})$. This completes the proof of Proposition \ref{poleresult}.\qed

\subsection{The gluing lemma}\label{Gluing}

We introduce the following condition:
\begin{equation}\label{relationexistenceBPZ}
Q+\gamma-\alpha_2-\alpha_3 < \frac{4}{\gamma} \wedge \gamma \wedge \min_{2 \leq i \leq 3} 2(Q-\alpha_i), \quad \alpha_2,\alpha_3<Q .
\end{equation}

\begin{lemma}\label{gluinglemma}
Suppose that $\alpha_2,\alpha_3$ satisfy condition \eqref{relationexistenceBPZ}. Then the function
\begin{equation*}
S(\alpha):= \begin{cases}
C_\gamma(\alpha, \alpha_2,\alpha_3), \; \text{if} \: \alpha < Q \\
R(\alpha) C_\gamma(2Q-\alpha, \alpha_2,\alpha_3), \; \text{if} \: \alpha > Q \\
\end{cases}
\end{equation*}
is the restriction on the real line of a holomorphic function $\bar S$ defined in a neighborhood of $Q$ and
given by 
\begin{equation}\label{holomorphicextension}
\bar{S}(\alpha)= - \frac{1}{\pi \mu}\caA(\tfrac{\gamma}{2},{\alpha}-\tfrac{\gamma}{2},{\alpha}_2,{\alpha}_3)C_\gamma({\alpha}-\gamma,{\alpha}_2,{\alpha}_3)
\end{equation}
where the function $\caA$ is defined in \eqref{Aformula}.

\end{lemma}

\proof  Let us first check that $\bar S$ is analytic in a neighborhood of $Q$.  By  \eqref{relationexistenceBPZ} we can find $\epsilon>0$ 
such that for all $\alpha \in[Q-\epsilon,Q+\epsilon]$
\begin{equation}\label{relationexistenceBPZepsilon}
2Q+\gamma-\alpha-\alpha_2-\alpha_3 < \frac{4}{\gamma} \wedge \gamma \wedge \min_{2 \leq i \leq 3} 2(Q-\alpha_i).
\end{equation}
By Theorem  \ref{analytic_hyperbolic} $\alpha\to C_\gamma({\alpha}-\gamma,{\alpha}_2,{\alpha}_3)$ is analytic in the region
\begin{equation}\label{relationexistenceBPZepsilon1}
2Q+\gamma-\alpha-\alpha_2-\alpha_3 < \frac{4}{\gamma} \wedge 2(Q+\gamma-\alpha)\wedge  \min_{2 \leq i \leq 3} 2(Q-\alpha_i) 
\end{equation}
which holds by \eqref{relationexistenceBPZepsilon} if $\epsilon <\frac{\gamma}{2}$.


Let first ${\alpha}\in (Q-\epsilon,Q)$. By   Theorem \ref{theo4point} the shift relation  \eqref{3pointconstanteqintro} holds in the form (take $\alpha_1+\frac{\gamma}{2}=\alpha$)
\begin{equation}\label{onemore}
C_\gamma({\alpha},\tilde{\alpha}_2,\tilde{\alpha}_3)
=- \frac{1}{\pi \mu}\caA(\tfrac{\gamma}{2},{\alpha}-\tfrac{\gamma}{2},\tilde{\alpha}_2,\tilde{\alpha}_3)C_\gamma({\alpha}-\gamma,\tilde{\alpha}_2,\tilde{\alpha}_3)
\end{equation} 
provided $2Q+\gamma-\alpha-\tilde{\alpha}_2-\tilde{\alpha}_3<0$ and $\gamma<\alpha<Q$. Thus, for $\epsilon$ small enough \eqref{onemore} holds for $\tilde{\alpha}_2,\tilde{\alpha}_3\in (Q-\epsilon,Q)$ and both sides are analytic in $\tilde{\alpha}_2,\tilde{\alpha}_3$ there. As we saw already, the RHS can be analytically continued to the values $(\tilde{\alpha}_2,\tilde{\alpha}_3)=(\alpha_2,{\alpha}_3)$.  By Theorem \ref{analytic_hyperbolic} the LHS is analytic in $\tilde{\alpha}_2,\tilde{\alpha}_3$ in a neighbourhood of
$2Q-\alpha-\tilde\alpha_2-\tilde\alpha_3 <0$. The point $(\tilde{\alpha}_2,\tilde{\alpha}_3)=(\alpha_2,{\alpha}_3)$ belongs to this region. 

Now let us turn to ${\alpha}\in (Q,Q+\epsilon)$. By Corollary \ref{shiftcoro} there exists $\eta>0$ s.t.  
\begin{equation*}
R( \alpha) C_\gamma(2Q- \alpha,\tilde{\alpha}_2,\tilde{\alpha}_3)
= - \frac{1}{\pi \mu}\caA(\tfrac{\gamma}{2}, \alpha-\tfrac{\gamma}{2},\tilde{\alpha}_2,\tilde{\alpha}_3)C_\gamma( \alpha-\gamma,\tilde{\alpha}_2,\tilde{\alpha}_3)
\end{equation*} 
provided $2Q+\gamma- \alpha-\tilde{\alpha}_2-\tilde{\alpha}_3< 0$ and $Q+\tfrac{\gamma}{2}- \alpha<\eta$. We saw above that the RHS extends to $ \alpha\in (Q,Q+\epsilon)$ and $(\tilde{\alpha}_2,\tilde{\alpha}_3)=(\alpha_2,{\alpha}_3)$.  By Theorem \ref{analytic_hyperbolic} the LHS extends to $ \alpha\in (Q,Q+\epsilon)$ and  $ \alpha-\tilde{\alpha}_2-\tilde{\alpha}_3< 0$. The point  $(\tilde{\alpha}_2,\tilde{\alpha}_3)=(\alpha_2,{\alpha}_3)$ belongs to this region. \qed

\subsection{Proof of the inversion relation \eqref{inversion}}\label{Step3}

The strategy is to combine the crossing relation Proposition \ref{2overgamma}:
\begin{equation}\label{Firstinverserelation}
C_\gamma (\alpha-\tfrac{2}{\gamma}, \epsilon, \alpha')= \tilde{T}(\alpha,\epsilon,\alpha')  R(\epsilon) C_\gamma (\alpha, 2Q-\epsilon-\tfrac{2}{\gamma}, \alpha')
\end{equation}
with the gluing Lemma \ref{gluinglemma} to obtain
\begin{equation}\label{Firstinverserelation1}
C_\gamma (\alpha-\tfrac{2}{\gamma}, \epsilon, \alpha')= \tilde{T}(\alpha,\epsilon,\alpha')   R(\epsilon) R(\alpha) C_\gamma (2Q-\alpha, 2Q-\epsilon-\tfrac{2}{\gamma}, \alpha')
\end{equation}
and then take the limit $\epsilon\to 2Q-\alpha$ and choose $\alpha'$ appropriately. 
To carry out this idea we need to check carefully the analyticity domains. Let us 
consider the following values for $\alpha,\alpha',\epsilon$:
 \begin{equation}\label{weightval}
\alpha= Q-\eta, \quad \epsilon=Q- \eta', \quad \alpha'=\tfrac{2}{\gamma}.
\end{equation}
where we will take $ |\eta|$  and $\eta'>0$  small  in what follows. 
\eqref{Firstinverserelation}  was proven in Proposition \ref{2overgamma} for $\alpha, \epsilon$ and $\alpha'$ close but strictly less than $Q$ with $\alpha+\alpha'+\epsilon>2Q+\frac{2}{\gamma}$. 
We use Theorem \ref{analytic_hyperbolic}  to extend the unit volume three point structure constant $\bar C_\gamma (\alpha-\tfrac{2}{\gamma}, \epsilon, \alpha')$  to the values \eqref{weightval}. 
The conditions in \eqref{ThextendedSeibergbounds1} become  $\eta+\eta' <\frac{4}{\gamma} \wedge (\frac{4}{\gamma}+2\eta)\wedge2\eta'\wedge\gamma$ and this gives  $\eta+\eta' <2\eta'$ as we are taking $|\eta|,\eta'$ small. Then  $\bar C_\gamma (\alpha-\tfrac{2}{\gamma}, \epsilon, \alpha')$ extends to the region $\eta'-\frac{2}{\gamma}<\eta<\eta'$. Note that $\eta'-\frac{2}{\gamma}<0$.  

Similarly, the condition for the function $ \bar C_\gamma (\alpha, 2Q-\epsilon-\tfrac{2}{\gamma}, \alpha')$ 
becomes  $\eta-\eta'< \frac{4}{\gamma} \wedge(2\eta)\wedge 2(\frac{2}{\gamma}-\eta')\wedge\gamma$
 with $\eta>0$. Since $|\eta|,\eta'$ are small this condition holds for $\eta<\eta'$. 
 In conclusion,  both unit volume structure constants extend to the region  $0<\eta<\eta'<\frac{\gamma}{2}$.
 The structure constants   $C_\gamma$ also extend to this region since the $s$-parameters in 
\eqref{unitvolumethreepoint}   are $\frac{-\eta-\eta'}{\gamma}$ and $\frac{\eta'-\eta}{\gamma}$  respectively and they do not take values in $\Z_-\cup\{0\}$. Hence \eqref{Firstinverserelation} holds in the common domain $0<\eta<\eta'<\frac{\gamma}{2}$. 

Next,  we apply the gluing Lemma \ref{gluinglemma} to the function  $  C_\gamma (\alpha, 2Q-\epsilon-\tfrac{2}{\gamma}, \alpha')$ to extend it to the region $\alpha>Q$ i.e. $\eta<0$. The condition \eqref{relationexistenceBPZ} becomes $\gamma-\eta'<\gamma\wedge 2(\frac{2}{\gamma}-\eta')$. This holds if $\eta'<\frac{4}{\gamma}-\gamma$. We conclude that the relation \eqref{Firstinverserelation1} holds for the values \eqref{weightval} with $0<-\eta<\eta'$ if $\eta'$ is small enough.

Now, we consider the limit of \eqref{Firstinverserelation1} as $\epsilon\uparrow 2Q-\alpha$ i.e. $\eta'\downarrow -\eta$. 
We get
\begin{equation*}
\lim_{\epsilon\uparrow 2Q-\alpha}(2Q-\alpha-\epsilon)C_\gamma (\alpha-\tfrac{2}{\gamma}, \epsilon, \tfrac{2}{\gamma}) 
=-2, \quad \lim_{\epsilon\uparrow 2Q-\alpha}(2Q-\alpha-\epsilon)C_\gamma (2Q-\alpha, 2Q-\epsilon-\tfrac{2}{\gamma}, \tfrac{2}{\gamma})
=2.
\end{equation*} 
Indeed, the two limits above correspond to insertions such that $s$ goes to $0$ in expression \eqref{expression3pointstruct} and can therefore be deduced from the following general fact for the unit volume structure constant defined in \eqref{unitvolumethreepoint}  and for $\alpha_1,\alpha_2,\alpha_3$ satisfying the extended Seiberg bounds \eqref{ThextendedSeibergbounds}
\begin{equation}\label{insertionto0}
 \bar{C}_\gamma(\alpha_1,\alpha_2,\alpha_3)=2\gamma^{-1}\quad \text{ if }\quad \sum_{i=1}^3\alpha_i=2Q.
\end{equation}
Also we get
\begin{equation*}
\lim_{\epsilon\uparrow 2Q-\alpha}(2Q-\alpha-\epsilon)l(a)=-\gamma, \quad \lim_{\epsilon\uparrow 2Q-\alpha}(2Q-\alpha-\epsilon)^{-1}l(b) =\tfrac{1}{\gamma}
\end{equation*}
and 
\begin{equation*}
\lim_{\epsilon\uparrow 2Q-\alpha} (2Q-\alpha-\epsilon)l(c)l(a+b-c)= l (1+\tfrac{2 \eta}{\gamma}) l(-\tfrac{2 \eta}{\gamma})=1.
\end{equation*}
We conclude that $R(2Q-\alpha)R(\alpha)=1$ for $\alpha$ close to $Q$ hence by analyticity in a neighbourhood of  $\R$. \qed

\subsection{ Proof of the $\frac{2}{\gamma}$ shift equation}\label{Step4}

We start from the following identity, obtained in Proposition \ref{2overgammabis},  for  $\epsilon, \alpha, \alpha'$ close to but strictly less than $Q$:
\begin{equation}\label{Thekeyequation}
R(\epsilon) C_\gamma(2Q-\epsilon-\tfrac{2}{\gamma}, \alpha,\alpha' )=  
L(\epsilon,\alpha,\alpha')R(\alpha) C_\gamma(\epsilon, 2Q-\alpha-\tfrac{2}{\gamma},\alpha' )
\end{equation}
where
\begin{equation}\label{Thekeyequationell}
L(\epsilon,\alpha,\alpha')=   \frac{l(c-1) l(c-a-b+1)}{l(c-a)l(c-b)} 
\end{equation}
with
\begin{align}\label{defab}
{ a}&= \tfrac{1}{\gamma} (\alpha'+\alpha+\epsilon-Q-\tfrac{4}{\gamma})-\hf  \quad { b}= \tfrac{1}{\gamma} (\alpha-\alpha'+\epsilon-Q)+\hf\quad c=1-\tfrac{2}{\gamma}(Q-\epsilon).
\end{align}
We will study \eqref{Thekeyequation} for
\begin{equation}\label{analytic}
\epsilon= \tfrac{\gamma}{2}+ \eta', \quad \alpha= \tfrac{\gamma}{2}+ \eta'', \quad \alpha'= \tfrac{2}{\gamma}+\eta 
\end{equation}
where $|\eta'|,\eta,\eta''$ will be taken small enough in what follows. 
We will   use Theorem \ref{analytic_hyperbolic}  to extend  equation \eqref{Thekeyequation} to  these values. First we have
\begin{equation}\label{Thekeyequationell00}
L(\epsilon,\alpha,\alpha')=  \frac{l(\tfrac{2\eta'}{\gamma}-\tfrac{4}{\gamma^2})l(1+\tfrac{4}{\gamma^2}- \tfrac{2\eta''}{\gamma})}{l(1+\tfrac{1}{\gamma}(\eta'-\eta''-\eta)l(\tfrac{1}{\gamma}(\eta-\eta''+\eta')}=
\frac{\eta+\eta'-\eta''}{\eta+\eta''-\eta'}  L_1(\epsilon,\alpha,\alpha')
\end{equation}
where $L_1$ is analytic around the point $ \tfrac{\gamma}{2},\tfrac{\gamma}{2}, \tfrac{2}{\gamma}$. Recalling \eqref{unitvolumethreepoint} we can then write
\begin{equation}\label{ThekeyequationA}
R(\epsilon) \bar C_\gamma(2Q-\epsilon-\tfrac{2}{\gamma}, \alpha,\alpha' )=  
L_2(\epsilon,\alpha,\alpha')R(\alpha) \bar C_\gamma(\epsilon, 2Q-\alpha-\tfrac{2}{\gamma},\alpha' )
\end{equation}
where $L_2$  is analytic around the point $( \tfrac{\gamma}{2},\tfrac{\gamma}{2}, \tfrac{2}{\gamma})$.
By \eqref{ThextendedSeibergbounds1} $\bar C_\gamma(2Q-\epsilon-\frac{2}{\gamma}, \alpha,\alpha' )$ 
extends to the region   $\eta'-\eta-\eta''<2\eta'$, $\eta'>0$ i.e. it is analytic for  $\eta',\eta,\eta''>0$.

For  $\bar C_\gamma(\epsilon, 2Q-\alpha-\tfrac{2}{\gamma},\alpha' )$  the condition in \eqref{ThextendedSeibergbounds1} becomes
$\eta''-\eta'-\eta  <  2\eta''
$ 
so that it extends to the region 
\begin{equation}\label{RHSext}
\eta+\eta''>   -\eta' .
\end{equation}
In particular the  eq. \eqref{Thekeyequation}  holds in the region  $\eta',\eta,\eta''>0$. 

Next, we want to extend  $C_\gamma(2Q-\epsilon-\frac{2}{\gamma}, \alpha,\alpha' )=C_\gamma(Q-\eta', \alpha,\alpha' )$ to $\eta'<0$ using the gluing lemma. The condition \eqref{relationexistenceBPZ} becomes
$
\gamma-\eta-\eta''<\gamma-2\eta
$ which requires $\eta<\eta''$. Therefore we get
\begin{equation}\label{Thekeyequationagain}
R(\epsilon) R(2Q-\epsilon-\tfrac{2}{\gamma})     C_\gamma(\epsilon+\tfrac{2}{\gamma}, \alpha,\alpha' )=   L(\epsilon,\alpha,\alpha')
R(\alpha) C_\gamma(\epsilon, 2Q-\alpha-\tfrac{2}{\gamma},\alpha' )
\end{equation}
for $\eta'<0$ sufficiently close to $0$. By condition \eqref{ThextendedSeibergbounds1} $C_\gamma(\epsilon+\tfrac{2}{\gamma}, \alpha,\alpha' )$ is analytic in $0\neq\eta'+\eta+\eta''<-2\eta'$. Combining this with \eqref{RHSext} and \eqref{Thekeyequationell} we conclude \eqref{Thekeyequationagain} holds in the region $-\eta'<\eta+\eta''<-3\eta'$. Therefore we may take the limit $\eta\to -\eta'$. Using also the  inversion relation $R(2Q-\epsilon-\tfrac{2}{\gamma})=R(\epsilon+\tfrac{2}{\gamma})^{-1}$ we end up with
\begin{equation}\label{TKEA}
\frac{R(\epsilon)}{ R(\epsilon+\tfrac{2}{\gamma}) }    C_\gamma(\epsilon+\tfrac{2}{\gamma}, \tfrac{\gamma}{2}+\eta'',Q-\epsilon)=  L(\epsilon,\tfrac{\gamma}{2}+\eta'',Q-\epsilon)
R(\tfrac{\gamma}{2}+\eta'') C_\gamma(\epsilon, Q-\eta'', Q-\epsilon)
\end{equation}
where we used $\alpha'=\tfrac{2}{\gamma}+\eta'=Q-\epsilon$. This identity holds in the region $0<\eta''<2\eta'$.

We will now take the limit of \eqref{TKEA} as $\eta''\to 0$. We have from equation \eqref{Thekeyequationell00}
\begin{equation}\label{Thekeyequationagain1}
  L(\epsilon,\tfrac{\gamma}{2}+\eta'',Q-\epsilon)=\frac{l(\tfrac{2}{\gamma}(\epsilon-Q))l(1+\tfrac{4}{\gamma^2}-\tfrac{2\eta''}{\gamma})}{l(\tfrac{2\epsilon}{\gamma}-\tfrac{\eta''}{\gamma}))l(-\tfrac{\eta''}{\gamma})}=-\frac{\eta''}{\gamma}\frac{l(\tfrac{2}{\gamma}(\epsilon-Q))l(1+\tfrac{4}{\gamma^2})
  }{l(\tfrac{2\epsilon}{\gamma})}+\caO({\eta''}^2)
\end{equation}
and by the first shift equation \eqref{shift1}
\begin{equation*}
R(\frac{\gamma}{2}+\eta'')= - \mu \pi \frac{R(\gamma+\eta'')}{ l(-\frac{\gamma^2}{4}) l(\frac{\gamma^2}{4}+
\frac{\gamma\eta''}{2})  l(2- \frac{\gamma \eta''}{2})}=\frac{2\mu\pi}{\gamma\eta''}\big(\frac{R(\gamma)}{ l(-\frac{\gamma^2}{4}) l(\frac{\gamma^2}{4}
)  }+\caO(\eta'')\big).
\end{equation*} 
Hence
\begin{equation*}
L(\epsilon,\tfrac{\gamma}{2}+\eta'',Q-\epsilon)
R(\tfrac{\gamma}{2}+\eta'') =\frac{\mu\pi\gamma^2}{8l(\tfrac{2\epsilon}{\gamma})}R(\gamma)l(\tfrac{2}{\gamma}(\epsilon-Q))l(1+\tfrac{4}{\gamma^2})+\caO(\eta'')
\end{equation*} 
where we used $l(x)l(-x)=-x^{-2}$. 

It remains to study the structure constants in \eqref{TKEA} as $\eta''\to 0$ using \eqref{expression3pointstruct}. We have
\begin{equation*}\label{Thekeyequationagain1}
\lim_{\eta''\to 0}\eta''   C_\gamma(\epsilon+\tfrac{2}{\gamma}, \tfrac{\gamma}{2}+\eta'',Q-\epsilon)=  2
\end{equation*}
since in \eqref{TKEA}  $s=\frac{\eta''}{\gamma}$. 
The second structure constant is dealt with by:
\begin{lemma}\label{lastone}
$\lim_{\eta''\to 0}\eta''  C_\gamma(\epsilon,Q-\eta'',Q-\epsilon)=-4$.
\end{lemma}
Hence we conclude
\begin{equation}\label{Thekeyequationagain11}
\frac{R(\epsilon)}{ R(\epsilon+\tfrac{2}{\gamma}) }   =-\frac{\mu\pi\gamma^2}{4l(\tfrac{2\epsilon}{\gamma})}R(\gamma)l(\tfrac{2}{\gamma}(\epsilon-Q))l(1+\tfrac{4}{\gamma^2})=-\frac{c_\gamma}{l(\tfrac{2\epsilon}{\gamma})l(-\tfrac{4}{\gamma^2})l(2+\tfrac{4}{\gamma^2}-\tfrac{2\epsilon}{\gamma})}
\end{equation}
with $c_\gamma=\frac{\gamma^2}{4}\mu\pi R(\gamma)$. This is the desired shift equation.

\vskip 2mm

\noindent{\it Proof of Lemma \ref{lastone}}. Using \eqref{expression3pointstruct} with $(\alpha_1,\alpha_2,\alpha_3)=  (Q-\eta'',Q-\epsilon, \epsilon)$ we have
\begin{equation*}
 C_\gamma(\epsilon,Q-\eta'',Q-\epsilon)=\tfrac{2}{\gamma}\mu^{\frac{\eta''}{\gamma}}\Gamma(-\tfrac{\eta''}{\gamma})\E
 (\int f(x)M_\gamma(d^2x))^{ \frac{\eta''}{\gamma}}
\end{equation*}
where 
$$
f(x)=   \frac{|x|_+^{\gamma (\alpha_1+\alpha_2+\alpha_3)}}{|x|^{\gamma \alpha_1} |x-1|^{\gamma \alpha_2} }  .
$$
Let $A:=\int 1_{|x|<\hf}f(x)M_\gamma(d^2x)$ and  $B:=\int 1_{|x|\geq\hf}f(x)M_\gamma(d^2x)$. By sub-additivity
\begin{equation*}
 \E A^{ \frac{\eta''}{\gamma}}\leq  (\int f(x)M_\gamma(d^2x))^{ \frac{\eta''}{\gamma}}\leq  \E A^{ \frac{\eta''}{\gamma}}+\E B^{ \frac{\eta''}{\gamma}}.
\end{equation*}
Now $\E B^p<\infty$ for some $p>0$ independent of $\eta''$. Thus $\lim_{\eta''\to 0}\E B^{ \frac{\eta''}{\gamma}}=1$ and then
\begin{equation*}
\lim_{\eta''\to 0}\eta'' C_\gamma(\epsilon+\tfrac{2}{\gamma}, \tfrac{\gamma}{2}+\eta'',Q-\epsilon)=2\lim_{\eta''\to 0} \E A^{ \frac{\eta''}{\gamma}}
.
\end{equation*}
Obviously
\begin{equation*}
\lim_{\eta''\to 0} \E A^{ \frac{\eta''}{\gamma}}=\lim_{\eta''\to 0} \E (\int_{|x| \leq 1} |x|^{-\gamma(Q-\eta'')} M_\gamma(dx))^{ \frac{\eta''}{\gamma}}
\end{equation*}
since only the neighborhood of $0$ contributes in the limit $\eta''\to 0$.
From \eqref{willaplied} we get
\begin{equation*}
r(\eta''):=\int_{|x| \leq 1} |x|^{-\gamma(Q-\eta'')} M_\gamma(dx)\stackrel{law}=e^{\gamma M_{\eta''}}\int_{-L_{-M_{\eta''}}}^\infty e^{\gamma\caB_s^{Q-\eta''}}Z_{s} ds
\end{equation*}
where $M_{\eta''}$ is the supremum of Brownian motion with drift $-\eta''$.
Then we may bound
\begin{equation*}
 \E (e^{\gamma M_{\eta''}}\int_{0}^1 e^{\gamma\caB_s^{Q-\eta''}}Z_{s} ds)^{ \frac{\eta''}{\gamma}}\leq\E r(\eta'')^{ \frac{\eta''}{\gamma}}\leq \E (e^{\gamma M_{\eta''}}\int_{-\infty}^\infty e^{\gamma\caB_s^{Q-\eta''}}Z_{s} ds)^{ \frac{\eta''}{\gamma}}
.
\end{equation*}
Let $I_{\eta''}:=\int_{-\infty}^\infty e^{\gamma\caB_s^{Q-\eta''}}Z_{s} ds$. Then by H\"older
\begin{equation*}
 \E  (e^{\gamma M_{\eta''}}I_{\eta''})^{ \frac{\eta''}{\gamma}}\leq  (\E  (e^{p\eta''M_{\eta''}})^{\frac{1}{p}}(\E (I_{\eta''})^{ \frac{q\eta''}{\gamma}})^{\frac{1}{q}}.
 \end{equation*}
Take $1/q=\sqrt{\eta''}$. Then $\limsup_{\eta''\to 0} (\E (I_{\eta''})^{ \frac{q\eta''}{\gamma}})^{\frac{1}{q}}=1$. Recalling that $\P(M_{\eta''} >v)=e^{- 2 \eta'' v}$ we then get since $p=1+\caO(\sqrt{\eta'' })$ that $\limsup_{\eta''\to 0}  (\E  e^{p\eta''M_{\eta''}})^{\frac{1}{p}}=2$. Hence 
\begin{equation*}
\limsup_{\eta''\to 0} \E\, r(\eta'')^{ \frac{\eta''}{\gamma}}\leq 2.
\end{equation*}
For the lower bound we set  $J_{\eta''}:=\int_0^1 e^{\gamma\caB_s^{Q-\eta''}}Z_{s} ds$ and use again H\"older
\begin{equation*}
\liminf_{\eta''\to 0} \E (e^{\gamma M_{\eta''}}J_{\eta''})^{ \frac{\eta''}{\gamma}}\geq \liminf_{\eta''\to 0}(\E e^{\frac{\eta'' }{p}M_{\eta''}})^p(\E (J_{\eta''})^{ -\frac{\eta''}{pq\gamma}})^{-\frac{1}{q}}=2.
 \end{equation*}
This finishes the proof.\qed

\begin{remark}
A straightforward computation yields that 
\begin{equation*}
\tfrac{\gamma^2}{4}  \mu \pi R^{{\rm DOZZ}}(\gamma)= {  (\pi  \mu   l(\tfrac{\gamma^2}{4}))^{\frac{4}{\gamma^2}} }{l(\tfrac{4}{\gamma^2})}^{-1}
\end{equation*}
and therefore we expect that $c_\gamma$ is given by \eqref{shift2a}. 
However, at this stage of the proof, we can not yet conclude this. 
$c_\gamma$ will be determined indirectly in subsection \ref{finalproofRDOZZ}. 
\end{remark}

\subsection{Proof that $R=R^{{\rm DOZZ}}$ }\label{finalproofRDOZZ}

 Let $ \psi(\alpha)=  \frac{R(\alpha)}{R^{{\rm DOZZ}}(\alpha)}$. $\psi$   is meromorphic in a neighborhood of  $\R$. Since $R$ and $R^{{\rm DOZZ}}$ obey the same $\frac{\gamma}{2}$ shift equation, the function $\psi$ is $\frac{\gamma}{2}$ periodic.  $\psi$ is strictly positive  in $(\frac{\gamma}{2},Q)$ so by periodicity $\psi$ is strictly positive on $\R$.  By the $\frac{2}{\gamma}$ shift equation, one has for all $\alpha \in \R$
\begin{equation*}
 \psi( \alpha)= C_\gamma \psi (\alpha+\frac{_2}{^\gamma})
\end{equation*}
for some constant $C_\gamma$. According to Liouville's theorem, if a continuous function $f$ say has two periods $T_1$ and $T_2$ such that $\frac{T_2}{T_1} \not\in \Q$ then $f$ is a constant function. Therefore, if $\frac{\gamma}{2}$ and  $\frac{2}{\gamma}$ are independent over the rationals (i.e. if $\gamma^2\notin \Q$) then we conclude that $C_\gamma=1$ and $\psi(\alpha)=\psi$ is constant in $\alpha$. From \eqref{defunitR} we see that $\bar R(Q)=1$ and from \eqref{deffullR} since $\Gamma(-x)x\to -1$ as $x\to 0$ we get $R(Q)=-1$. On the other hand, from \eqref{defRDOZZ} follows  $R^{{\rm DOZZ}}(Q)=-1$ hence the constant $\psi=1$. Hence $R(\alpha)=R^{{\rm DOZZ}}(\alpha)$ for all $\alpha$. The case $\gamma^2\in \Q$ follows by continuity. This concludes the proof.

\section{Proof of the DOZZ formula }\label{proofDOZZ}

We suppose that $\gamma^2\notin \Q$; the general case follows by  continuity. Let us fix $\alpha_2,\alpha_3$ in $(Q-\eta,Q)$ for $\eta$ sufficiently small and consider the function $F:\alpha_1 \mapsto C_\gamma (\alpha_1,\alpha_2,\alpha_3)$. Let us collect what we have proven about $F$. By Theorem \ref{analytic_hyperbolic} $F$ is analytic on $(2 \eta,Q)$ and by Theorem \ref{theo4point} it satisfies  
 the $\frac{\gamma}{2}$ shift equation \eqref{3pointconstanteqintro}, for $\frac{\gamma}{2}+2 \eta<\alpha_1< \frac{2}{\gamma}$. Therefore  $F$ extends to a meromorphic function on a strip of the form $\R \times (-\beta,\beta)$ with $\beta>0$ satisfying \eqref{3pointconstanteqintro}. We call this extension $F$ too. 

 
Now, using the exact expression for $R$ (or relation \eqref{shift2} with $c_\gamma= \mu \pi l(\frac{\gamma^2}{4})  )^{\frac{4}{\gamma^2} }l(\frac{4}{\gamma^2})^{-1}$)  Theorem \ref{theo4point2overgamma} can be written as   
\begin{align*}
\mathcal{T}_{-\frac{2}{\gamma}}(z)= &C_\gamma(\alpha_1-\frac{2}{\gamma},\alpha_2,\alpha_3) |F_-(z)|^2 \\
&- \frac{(\mu \pi l(\frac{\gamma^2}{4})  )^{\frac{4}{\gamma^2} }}{l(\frac{4}{\gamma^2})} \frac{R(\alpha_1+\frac{2}{\gamma})}{   l(-\frac{4}{\gamma^2}) l(\frac{2 \alpha_1}{\gamma})  l(2+\frac{4}{\gamma^2}- \frac{2 \alpha_1}{\gamma}) } C_\gamma(2Q-\alpha_1-\frac{2}{\gamma}, \alpha_2,\alpha_3) |F_+(z)|^2.
 \end{align*}
By the gluing Lemma, the extension $F$ is given in a neighborhood of $\alpha=Q$ by $F(\alpha)=R(\alpha)F (2Q-\alpha)$. Hence, one can infer from the above expression the shift equation \eqref{3pointconstanteqintrodual} for $\alpha_1 \in \R \times (-\beta,\beta)$ (same argument as the one used to derive \eqref{3pointconstanteqintro}). Hence $F$ satisfies both \eqref{3pointconstanteqintro} and  \eqref{3pointconstanteqintrodual}.

Now, we consider the function $\psi_{\alpha_2,\alpha_3}: \alpha_1 \mapsto \frac{C_\gamma (\alpha_1,\alpha_2,\alpha_3)}{C_\gamma^{{\rm DOZZ}} (\alpha_1,\alpha_2,\alpha_3)}$ in the strip $\R \times (-\beta,\beta)$. This function is holomorphic since $C_\gamma$ and $C_\gamma^{{\rm DOZZ}} $ are meromorphic with the same simple poles and zeros (which can be read off the $\frac{\gamma}{2}$ shift equation \eqref{3pointconstanteqintro}). Furthermore, $\psi_{\alpha_2,\alpha_3}$ is $\gamma$ and $\frac{4}{\gamma}$ periodic since $C_\gamma$ and $C_\gamma^{{\rm DOZZ}} $ both satisfy \eqref{3pointconstanteqintro} and   \eqref{3pointconstanteqintrodual}. According to Liouville's theorem, if a continuous function $f$ say has two periods $T_1$ and $T_2$ such that $\frac{T_2}{T_1} \not\in \Q$ then $f$ is a constant function. Therefore, using Liouville's theorem, $\psi_{\alpha_2,\alpha_3}(\alpha_1)=c_{\alpha_2,\alpha_3}$ for some constant $c_{\alpha_2,\alpha_3}$ depending on $\alpha_2,\alpha_3$. 

Since $C_\gamma$ and $C_\gamma^{{\rm DOZZ}}$ are symmetric in their arguments we obtain $\psi_{\alpha_2,\alpha_3}(\alpha_1)=\psi_{\alpha_1,\alpha_3}(\alpha_2)=\psi_{\alpha_1,\alpha_2}(\alpha_3)$ for  $\alpha_1,\alpha_2,\alpha_3\in(Q-\eta,Q)$. Hence $c_{\alpha_2,\alpha_3}$ is constant in $\alpha_2,\alpha_3$. Therefore  $C_\gamma (\alpha_1,\alpha_2,\alpha_3)=a_\gamma C_\gamma^{{\rm DOZZ}} (\alpha_1,\alpha_2,\alpha_3)$ for $\alpha_1,\alpha_2,\alpha_3$ satisfying \eqref{ThextendedSeibergbounds} with $N=3$ for some constant $a_\gamma$ (by analycity). Finally $a_\gamma=1$ since  both $C_\gamma$ and $ C_\gamma^{{\rm DOZZ}}$ satisfy Lemma \ref{defR}.\qed

\section{Appendix}

  \subsection{The Cameron-Martin theorem}
  We state here the classical Cameron-Martin theorem for the GFF X. Let $\mathcal {S}'(\C)$ be the space of tempered distributions. It is well known that $X$ lives in $\mathcal {S}'(\C)$. 
  
 \begin{theorem}{\bf Cameron-Martin theorem }\label{th:Girsanov}

Let $\mathcal{X}$ be some some Gaussian variable which is measurable with respect to the GFF $X$. Let $F$ be some bounded continuous function on $\mathcal {S}'(\C)$. Then we have the following identity
\begin{equation*}
\E[  e^{\mathcal{X}-\frac{\E[\mathcal{X}^2]}{2}}   F( (X(x))_{x \in \C}   )   ]= \E[    F( (X(x)   +\E[X(x) \mathcal{X}])_{x \in \C}   )   ].
\end{equation*}
\end{theorem}  

  In particular we get the following corollary:
  
\begin{corollary}\label{coro:Girsanov}
 Let $F$ be some bounded continuous function on $\mathcal {S}'(\C)$ and $f$ some bounded continuous function on $\C$. Then we have the following identity
\begin{equation}\label{coroGirs}
\E[   \left ( \int_{\C}  f(u) M_\gamma(d^2u)   \right )  F( (X(x))_{x \in \C}   )   ]= \int_{\C} f(u)\E[    F( (X(x)   +\E[X(x) X(u)])_{x \in \C}   )   ]  \frac{d^2u}{|u|_+^4}.
\end{equation}
\end{corollary}  
 
 \proof
 For $\epsilon>0$, if $X_\epsilon$ denotes the circle average of $X$ then by Fubini (interchanging $\E[.]$ and $\int_{\C}$)
\begin{align*}
& \E[   \left ( \int_{\C}  f(u) e^{\gamma X_\epsilon(u)-\frac{\gamma^2}{2} \E[X_\epsilon(u)^2]}\frac{d^2u}{|u|_+^4}   \right )  F( (X(x))_{x \in \C}   )   ]  \\
&  = \int_{\C}  f(u) \E[   e^{\gamma X_\epsilon(u)-\frac{\gamma^2}{2} \E[X_\epsilon(u)^2]}   F( (X(x))_{x \in \C}   )   ]  \frac{d^2u}{|u|_+^4} \\
& = \int_{\C} f(u)  \E[    F( (X(x)   +\E[X(x) X_\epsilon(u)])_{x \in \C}   )   ]   \frac{d^2u}{|u|_+^4} 
\end{align*} 
where in the last line we have used the Girsanov Theorem \ref{th:Girsanov}. In conclusion, we have
\begin{equation*}
 \E[   \left ( \int_{\C}  f(u) e^{\gamma X_\epsilon(u)-\frac{\gamma^2}{2} \E[X_\epsilon(u)^2]}\frac{d^2u}{|u|_+^4}   \right )  F( (X(x))_{x \in \C}   )   ]    = \int_{\C} f(u)  \E[    F( (X(x)   +\E[X(x) X_\epsilon(u)])_{x \in \C}   )   ]   \frac{d^2u}{|u|_+^4} .
\end{equation*} 
We conclude by letting $\epsilon$ go to $0$. 
 \qed

\subsection{Chaos estimates}\label{chaosestimates}
 
We list here estimates for chaos integrals that are used frequently in the paper. Some estimates are standard in the literature on GMC and other estimates were recently proved in \cite[section 6]{KRV}.

\subsubsection{Standard moment estimates}
We start by reviewing the standard estimates and for these estimates we refer to the review \cite{review}. For any open and bounded subset $\mathcal{O}$, the following condition holds on moments (see \cite{review}):
\begin{equation}\label{condmom}
\E[   M_\gamma(\mathcal{O})^p  ]  < \infty \quad \iff \quad p \in (-\infty,\frac{4}{\gamma^2}).
\end{equation}
Moreover, if $p \in (-\infty,\frac{4}{\gamma^2})$ and $z \in \C$ then there exists a constant $C>0$ (depending on $z$ and $p$) such that for all $\epsilon \leq 1$ (see \cite{review})
\begin{equation}  \label{momentscaling}
\E  ( M_\gamma(B(z,\epsilon))^p)\leq C\epsilon^{\gamma Qp-\frac{\gamma^2p^2}{2}}.
\end{equation}

 Let $A(z,\epsilon)$ be the annulus with radii $\epsilon,2\epsilon$ and center $z$. We get as corollary of  \eqref{momentscaling} that for $p\in [0,\frac{4}{\gamma^2})$
 
\begin{equation}  \label{annulus}
\E  \left (  \int_{A(z,\epsilon)}  |x-z|^{-\gamma\alpha}M_\gamma(d^2x)\right)^p\leq C\epsilon^{\gamma (Q-\alpha)p-\frac{\gamma^2 p^2}{2}}.
\end{equation}
For negative moments, we have for $\alpha>Q$ and $p>0$ such that $\alpha-Q<\gamma p$ (see the methods of \cite[section 6]{KRV}):
 \begin{align}\label{freezingannulus} 
 \E\Big[\Big(\int_{|x-z|>\epsilon}|x-z|^{-\gamma\alpha}M_\gamma(d^2x)\Big)^{-p}\Big] \leq C\epsilon^{\hf(\alpha-Q)^2}.
\end{align}

\subsubsection{Fusion estimate}

The following result follows from the methods of \cite[section 6]{KRV}:

\begin{lemma}\label{fusion4}
Assume $(\alpha_i)_{i=1,\dots,4}$ are real numbers satisfying $\alpha_i<Q$ and $p:=\gamma^{-1}(\sum_{i=1}^4\alpha_i-2Q) >0$. 
Consider   $y_1,y_2,y_3,y_4 \in\C$ such that 
$|y_1-y_2|\leq |y_1-y_3|\leq |y_2-y_3|\leq \min_{i\in\{1,2,3\}}|y_4-y_i|$. 


1) If $\alpha_1+\alpha_2<Q$,      $\alpha_1+\alpha_2+\alpha_3\geq Q$ and $\alpha_4\geq 0$  then
\begin{align*}
\E\Big[&\Big(\int_{B(y_1,10)} \prod_{i=1}^4|u-y_i|^{-\gamma\alpha_i}M_\gamma(d^2u)\Big)^{-p-2}]\\
\leq &C\Big(\frac{|y_1-y_3|}{|y_1-y_4|}\Big)^{\frac{1}{2}(\alpha_1+\alpha_2+\alpha_3-Q)^2}  |y_1-y_4|^{\frac{1}{2}(\alpha_1+\alpha_2+\alpha_3+\alpha_4-Q)^2  }.
\end{align*}

2) If $\alpha_1+\alpha_2> Q$, $\alpha_3\leq 0$ and $\alpha_3+\alpha_4\geq 0$ then
\begin{align*}
\E\Big[&\Big(\int_{B(y_1,10)} \prod_{i=1}^4|u-y_i|^{-\gamma\alpha_i}M_\gamma(d^2u)\Big)^{-p-2}]\\
\leq &C\Big(\frac{|y_1-y_2|}{|y_1-y_3|}\Big)^{\frac{1}{2}(\alpha_1+\alpha_2-Q)^2} \Big(\frac{|y_1-y_3|}{|y_1-y_4|}\Big)^{\frac{1}{2}(\alpha_1+\alpha_2+\alpha_3-Q)^2 -\frac{\alpha_3^2}{2}}|y_1-y_4|^{\frac{1}{2}(\alpha_1+\alpha_2+\alpha_3+\alpha_4-Q)^2 }.
\end{align*}
\end{lemma}

\subsubsection{FKG inequality}

Finally, we recall a result on log-correlated fields in dimension 2 which comes out of a construction in \cite{MRM}. Recall that From \cite{MRM}, there exists a Gaussian white noise $\mu$ on some measure space $(\mathcal{S},\nu)$ ($\nu$ is a Radon measure) and deterministic subsets $(C(x))_{|x| \leq \frac{1}{2}}$ of $\mathcal{S}$ such that the field $(\tilde{X}(x))_{|x| \leq \frac{1}{2}}$ defined by
\begin{equation}\label{deffieldtildeX}
\tilde{X}(x)= \mu(C(x))
\end{equation}
is a Gaussian field with covariance given by
\begin{equation*}
\E[ \tilde{X}(x)\tilde{X}(y) ]= \ln \frac{1}{|x-y|}+c
\end{equation*}
where $c$ is some positive constant. In particular, the construction \eqref{deffieldtildeX} implies that $\tilde{X}$ satisfies the FKG inequality; more precisely, if $F,G$ are two increasing functions in each coordinate $\tilde{X}(x)$ then 
\begin{equation*}
\E[   F((\tilde{X}(x))_{|x| \leq \frac{1}{2}}) G((\tilde{X}(x))_{|x| \leq \frac{1}{2}}) ] \geq \E[   F((\tilde{X}(x))_{|x| \leq \frac{1}{2}})  ] \E[   G((\tilde{X}(x))_{|x| \leq \frac{1}{2}}) ] .
\end{equation*}
The above continuum version of the FKG inequality can be deduced from the standard one (see \cite[section 2.2]{grimmett} for the case of countable product) by discretization and taking the limit as the mesh of discretization goes to $0$. 
Since $(\tilde{X}(x))_{|x| \leq \frac{1}{2}}$ has same distribution as $(X(x)+\sqrt{c}Y)_{|x| \leq \frac{1}{2}}$ where $Y$ is a fixed standard Gaussian independent from $X$, this implies that $(X(x)+\sqrt{c}Y)_{|x| \leq \frac{1}{2}}$ also satisfies the FKG inequality.  

\subsection{Proof of lemma \ref{defmomR}}\label{app:momR}

By symmetry, it is enough to show that
\begin{equation*}
\E[ \left ( \int_{0}^\infty  e^{   \gamma \mathcal{B}_s^\alpha  } Z_s ds  \right )^{p}]< \infty. 
\end{equation*}
Let first $p>0$. If $0<p\leq 1$ we have by subadditivity
\begin{equation*}
\E[ \left ( \int_{0}^\infty  e^{   \gamma \mathcal{B}_s^\alpha  } Z_s ds  \right )^{p}]\leq \sum_{n=1}^{\infty}  \E \Big [ \left ( \int_n^{n+1}  e^{   \gamma \mathcal{B}_s^\alpha  } Z_s ds  \right )^{p} \Big ] 
\end{equation*}
and for $1<p<\frac{4}{\gamma^2}$ by convexity
\begin{equation*}
[\E \left ( \int_{0}^\infty  e^{   \gamma \mathcal{B}_s^\alpha  } Z_s ds  \right )^{p}]^{1/p} \leq \sum_{n=1}^{\infty}  \Big [ \E \left ( \int_n^{n+1}  e^{   \gamma \mathcal{B}_s^\alpha  } Z_s ds  \right )^{p} \Big ]^{1/p} .
\end{equation*}

We set $\nu= Q-\alpha$.  The process $ \mathcal{B}_s^\alpha$ is stochastically dominated by a Brownian motion with drift $-\nu$  starting from origin and conditioned to stay below $1$ (see Appendix, subsection \ref{reminderdiff}); hence we have if $B_s$ is a standard Brownian motion starting from $0$ 
\begin{align*}
& \E \Big [ \left ( \int_n^{n+1}  e^{   \gamma \mathcal{B}_s^\alpha  } Z_s ds  \right )^{p} \Big ]    \leq C  \E \Big [ \left (   1_{B_n -\nu n \leq 1} \int_n^{n+1}  e^{   \gamma (B_s- \nu s)  }Z_s ds   \right )^{p} \Big ]    \\
& \leq C \E  [ (\int_n^{n+1} Z_s ds )^{p} ]  \E[   e^{\gamma p \sup_{s \in [n,n+1]} (B_s-B_n) }]  \E[ 1_{B_n -\nu n \leq 1} e^{\gamma p ( B_n -\nu n) }  ]   \\
& \leq C  \E[ 1_{B_n -\nu n \leq 1} e^{\gamma p ( B_n -\nu n) }  ]= Cn^{-\hf}\int_{-\infty}^{1}  e^{\gamma p  y}  e^{- \frac{(y+\nu n)^2}{2n}}dy
\end{align*} 
where we used \eqref{zmoment}. Considering separately $y<-\frac{\nu n}{2}$ and $y\in [-\frac{\nu n}{2},1]$ the last integral is seen to be exponentially small in $n$ and the claim follows. 


Let now $p=-q<0$. Set $\tau_{-1}= \inf \lbrace  s \geq 0, \; \mathcal{B}^{\alpha}_s=-1 \rbrace$. The process $\mathcal{B}^{\alpha}_{s+\tau_{-1}}+1$ is a Brownian motion with drift $-\nu$ starting from $0$ and conditioned to stay below $1$. Therefore, we have if $B_s$ is a standard Brownian motion starting from $0$ and  $\beta :=\sup_{s\geq 0}(B_s-\nu s)$
\begin{align*}
& \E[ \left ( \int_{0}^\infty  e^{   \gamma \mathcal{B}_s^\alpha  } Z_s ds  \right )^{-q}]   \leq \E[ \left ( \int_{\tau_{-1}} ^{\tau_{-1}+1}  e^{   \gamma \mathcal{B}_s^\alpha  } Z_s ds  \right )^{-q}]   = C \E[  1_{\beta\leq 1}    \left ( \int_{0} ^{1}  e^{   \gamma ( B_s-\nu s)  } Z_s ds  \right )^{-q}  ]    
\end{align*}  
since $Z_s$ is stationary and $ \mathcal{B}_s^\alpha $ is independent from $Z$. Finally, we conclude by
\begin{align*}
 & \E[ 1_{\beta\leq 1}     \left ( \int_{0} ^{1}  e^{   \gamma ( B_s-\nu s)  } Z_s ds  \right )^{-q}  ]    \leq \E[ 1_{\beta\leq 1}   e^{-q \inf_{s \in [0,1]} (B_s-\nu s)  }  ]\, \E (  \int_0^1 Z_s ds )^{-q}   < \infty  
\end{align*}
where \eqref{zmoment} was used. \qed

\subsection{A reminder on diffusions}\label{reminderdiff}

A drifted Brownian motion $(B_t+\mu t)$ with $\mu>0$ is a diffusion with generator $\mathcal{G}_\mu= \frac{1}{2} \frac{d^2}{dx^2} + \mu \frac{d}{dx}$. When seen until hitting $b>0$,  the dual process $Y_b$ of $B_t+\mu t$  is a diffusion with generator $\frac{1}{2} \frac{d^2}{dx^2} - \mu  \coth (\mu (b-x))\frac{d}{dx}$. Therefore, $b-Y_b$ has generator $\frac{1}{2} \frac{d^2}{dx^2} + \mu  \coth (\mu x)\frac{d}{dx}$ which is the generator of $(B_t+\mu t)$ conditioned to be positive. We denote this process $\mathcal{B}_t^\mu$. We also denote by $\mathcal{B}_t^0$ the standard 3d Bessel process which corresponds to the case $\mu=0$.  

We have the following comparison principle:
\begin{lemma}\label{lemmecompdiff}
There exists a probability space such that for $0 \leq \mu <\mu'$, we have  almost surely for all $t$: $\mathcal{B}_t^\mu \leq \mathcal{B}_t^{\mu'}$.
\end{lemma}

\proof
For all $x >0$, we consider the drift $\varphi_x(\mu)=  \mu  \coth ( \mu x)$ as a a function of $\mu \in [0,\infty)$. A straightforward computation yields
\begin{equation*}
\varphi_x'(\mu)= \frac{e^{4 \mu x}- 4 \mu x e^{2 \mu x} -1}{(e^{2\mu x}-1)^2}.  
\end{equation*}
Therefore $\varphi_x'(\mu) \geq 0$ since $e^{u}-ue^{\frac{u}{2}}-1 \geq 0$ for all $u \geq 0$.
\qed

We will need another comparison principle. Let $B^{\mu,A}_t$ be the drifted Brownian motion $(B_t+\mu t)$ starting from $0$ and conditioned to be above $-A$ with $A>0$.

\begin{lemma}\label{lemmecompdiff}
Let $\mu \geq 0$. There exists a probability space such that for $A>0$ we have almost surely for all $t$:  $\mathcal{B}_t^\mu\geq B^{\mu,A}_t $.
\end{lemma}

\proof 
This can also be read off the drift. Indeed, for $\mu,x$ fixed, we consider $\psi_{\mu,x} (A)=  \mu  \coth ( \mu (x+A))$. We have 
\begin{equation*}
\forall x \geq -A, \; \; \; \psi_{\mu,x}' (A) = \mu^2 (1- \coth ( \mu (x+A))^2)  \leq 0.
\end{equation*}
\qed

\subsection{Functional relations on $\Upsilon_{\frac{\gamma}{2}}$ and $R^{{\rm DOZZ}}$}
The function $\Upsilon_{\frac{\gamma}{2}}$ defined by \eqref{def:upsilon} can be analytically continued to $\C$ and it satisfies the following remarkable functional relations for $z \in \C$
\begin{equation}\label{shiftUpsilon}
\Upsilon_{\frac{\gamma}{2}} (z+\frac{\gamma}{2}) =\frac{\Gamma(\frac{\gamma}{2}z )}{\Gamma (1-\frac{\gamma}{2}z)} (\frac{\gamma}{2})^{1-\gamma z}\Upsilon_{\frac{\gamma}{2}} (z), \quad
\Upsilon_{\frac{\gamma}{2}} (z+\frac{2}{\gamma}) =\frac{\Gamma(\frac{2}{\gamma}z )}{\Gamma (1-\frac{2}{\gamma}z)} (\frac{\gamma}{2})^{\frac{4}{\gamma} z-1} \Upsilon_{\frac{\gamma}{2}} (z).
\end{equation}
The function $\Upsilon_{\frac{\gamma}{2}}$ has no poles in $\C$ and the zeros of $\Upsilon_{\frac{\gamma}{2}}$ are simple (if $\gamma^2 \not \in \Q$) and given by the discrete set $(-\frac{\gamma}{2} \N-\frac{2}{\gamma} \N) \cup (Q+\frac{\gamma}{2} \N+\frac{2}{\gamma} \N )$: for more on the function $\Upsilon_{\frac{\gamma}{2}}$ and its properties, see the reviews \cite{nakayama,Rib,Tesc1} for instance.

With definition \eqref{defRDOZZ} and a little algebra, one can show that $R^{{\rm DOZZ}}(\alpha)$ satisfies the following shift equation for all $\alpha \in \C$ 

\begin{equation}\label{shift1DOZZ}
R^{{\rm DOZZ}}(\alpha-\frac{\gamma}{2})= - \mu \pi \frac{R^{{\rm DOZZ}}(\alpha)}{ l(-\frac{\gamma^2}{4}) l(\frac{\gamma\alpha}{2}-\frac{\gamma^2}{4})  l(2+\frac{\gamma^2}{2}- \frac{\gamma \alpha}{2})}
\end{equation}

as well as the dual shift equation for all $\alpha \in \C$ 

\begin{equation}\label{shift3DOZZ}
R^{{\rm DOZZ}}(\alpha)= - \frac{(\mu \pi l(\frac{\gamma^2}{4})  )^{\frac{4}{\gamma^2} }}{l(\frac{4}{\gamma^2})} \frac{R^{{\rm DOZZ}}(\alpha+\frac{2}{\gamma})}{   l(-\frac{4}{\gamma^2}) l(\frac{2 \alpha}{\gamma})  l(2+\frac{4}{\gamma^2}- \frac{2 \alpha}{\gamma}) }.
\end{equation}

\subsection{Derivation of $R^{{\rm DOZZ}}$ from $C_\gamma^{{\rm DOZZ}}$}
Recall that the function $\Upsilon_{\frac{\gamma}{2}}$ satisfies the  shift equations \eqref{shiftUpsilon}. According to the DOZZ formula \eqref{DOZZformula}, since $\Upsilon_{\frac{\gamma}{2}}(0)=0$, we get for $\alpha> \frac{\gamma}{2}$ and using the above relations 
\begin{align*}
\epsilon C(\alpha,\epsilon,\alpha)  & \underset{\epsilon \to 0}{\sim}  4 (\pi \: \mu \:  l(\frac{\gamma^2}{4})  \: (\frac{\gamma}{2})^{2 -\gamma^2/2} )^{\frac{2 (Q -\alpha)}{\gamma}}   \frac{ \epsilon^2 \Upsilon_{\frac{\gamma}{2}}'(0)^2\Upsilon_{\frac{\gamma}{2}}(\alpha)^2}{   \epsilon^2 \Upsilon_{\frac{\gamma}{2}}'(0)^2   \Upsilon_{\frac{\gamma}{2}}(\alpha-Q) \Upsilon_{\frac{\gamma}{2}}(\alpha)}  \\
 &  =    4 (\pi \: \mu \:  l(\frac{\gamma^2}{4})  \: (\frac{\gamma}{2})^{2 -\gamma^2/2} )^{\frac{2 (Q -\alpha)}{\gamma}}   \frac{ \Upsilon_{\frac{\gamma}{2}}(\alpha)}{     \Upsilon_{\frac{\gamma}{2}}(\alpha-Q) }   \\
& =  4 (\pi \: \mu \:  l(\frac{\gamma^2}{4})  \: (\frac{\gamma}{2})^{2 -\gamma^2/2} )^{\frac{2 (Q -\alpha)}{\gamma}} \frac{\Gamma(\frac{\gamma (\alpha-Q+\frac{2}{\gamma})}{2})}{\Gamma(1- \frac{\gamma (\alpha-Q+\frac{2}{\gamma})}{2})} (\frac{\gamma}{2})^{1-\gamma (\alpha-Q+\frac{2}{\gamma})} \:  \frac{\Gamma(\frac{2 (\alpha-Q)}{\gamma })}{\Gamma(1- \frac{2 (\alpha-Q)}{\gamma})} (\frac{\gamma}{2})^{\frac{4}{\gamma} (\alpha-Q)-1}   \\ 
& =  4   ( \frac{2}{\gamma} )^{-2}(\pi \: \mu \:  l(\frac{\gamma^2}{4})  )^{\frac{2 (Q -\alpha)}{\gamma}} \frac{\Gamma(1- \frac{\gamma (Q-\alpha)}{2})}{\Gamma(\frac{\gamma (Q-\alpha)}{2})} \:  \frac{\Gamma(-\frac{2 (Q-\alpha)}{\gamma })}{\Gamma(1+ \frac{2 (Q-\alpha)}{\gamma})}  \\ 
& = - 4  (\pi \: \mu \:  l(\frac{\gamma^2}{4})  )^{\frac{2 (Q -\alpha)}{\gamma}} \frac{\Gamma(-\frac{\gamma (Q-\alpha)}{2})}{\Gamma(\frac{\gamma (Q-\alpha)}{2})} \:  \frac{\Gamma(-\frac{2 (Q-\alpha)}{\gamma })}{\Gamma(\frac{2 (Q-\alpha)}{\gamma})} \\
& = 4 R^{{\rm DOZZ}}(\alpha). 
\end{align*}

\subsection{An integral formula}
We have:

\begin{lemma}\label{lemmaintegral}
For all $p>0$ and $a \in (1,2)$ the following identity holds 
\begin{equation*}
 \int_0^{\infty}    \left (  \frac{1}{(1+v )^p}-1  \right ) \frac{1}{v^a}  dv=   \frac{\Gamma(-a+1)  \Gamma (p+a-1)}{   \Gamma (p) }.
\end{equation*}
\end{lemma}

\proof
We set $\bar{a}=-a+1$ and $\bar{b}=p+a-1$. 
We have 
\begin{align*}
 & \int_0^{1}    \left (  \frac{1}{(1+v )^p}-1  \right ) \frac{1}{v^a}  dv- \int_1^{\infty}   \frac{1}{v^a}  dv  
 = -\frac{1}{a-1} \sum_{k \geq 1} (-1)^k  \frac{(p)_k (-a+1)_k}{k! (-a+2)_k}  -\frac{1}{a-1}  \\
& = -\frac{1}{a-1} \sum_{k \geq 0} (-1)^k  \frac{(p)_k (-a+1)_k}{k! (-a+2)_k}  = \frac{1}{\bar{a}}\:  {}_{2}F_1(\bar{a}, \bar{a}+ \bar{b}, \bar{a}+1, z=-1)  .
\end{align*}
Next, we have
\begin{align*}
  \int_1^{\infty}   \frac{1}{(1+v )^p} \frac{1}{v^a}   dv  & =  \int_0^{1}   \frac{1}{(1+v )^p} v^{p+a-2}   dv  
 = \frac{1}{p+a-1}\sum_{k \geq 0} (-1)^k \frac{ (p)_k(p+a-1)_k}{k!(p+a)_k}   \\
& = \frac{1}{\bar{b}} \:  {}_{2}F_1(\bar{b}, \bar{a}+ \bar{b}, \bar{b}+1, z=-1) .
\end{align*}
Finally, we use the following formula (see \cite{hyp}):
\begin{equation*}
\bar{b}\:  {}_{2}F_1(\bar{a}, \bar{a}+ \bar{b}, \bar{a}+1, z=-1)+ \bar{a}  \:  {}_{2}F_1(\bar{b}, \bar{a}+ \bar{b}, \bar{b}+1, z=-1)= \frac{\Gamma (\bar{a}+1) \Gamma (\bar{b}+1)}{\Gamma (\bar{a}+\bar{b})}.
\end{equation*}
This yields the desired relation since $\Gamma (z+1)= z \Gamma (z)$.\qed

\subsection{Some identities}
We have the following identity for all $z$
\begin{equation*}
 \int_{\C}    \frac{|u-z|^{\frac{\gamma^2}{2}}-|u|^{\frac{\gamma^2}{2}} -  \frac{\gamma^2}{4} |u|^{\frac{\gamma^2}{2}} ( \frac{z}{u}+ \frac{\bar{z}}{\bar{u}} )  }{ |u|^{\gamma \alpha_1}   } d^2u= |z|^{\gamma (Q-\alpha_1)}  \frac{\pi}{ l(\frac{\gamma \alpha_1}{2})   l(- \frac{\gamma^2}{4})l (2- \frac{\gamma \alpha_1}{2}+ \frac{\gamma^2}{4})} .
\end{equation*}
Applying $\partial_{z}^2$ to this we get
\begin{equation*}
\frac{\gamma^2}{4} (\frac{\gamma^2}{4}-1)  \int_{\C}    \frac{|u-z|^{\frac{\gamma^2}{2}}  }{(z-u)^2 |u|^{\gamma \alpha_1}   } d^2u= \frac{\gamma (Q-\alpha_1)}{2} (\frac{\gamma (Q-\alpha_1)}{2}-1) \frac{|z|^{\gamma (Q-\alpha_1)}}{z^2}  \frac{\pi}{ l(\frac{\gamma \alpha_1}{2})   l(- \frac{\gamma^2}{4})l (2- \frac{\gamma \alpha_1}{2}+ \frac{\gamma^2}{4})} .
\end{equation*}
Hence for $z=1$ this yields
\begin{equation}\label{gammaz}
\frac{\gamma^2}{4} (\frac{\gamma^2}{4}-1)  \int_{\C}    \frac{|u-1|^{\frac{\gamma^2}{2}}  }{(1-u)^2 |u|^{\gamma \alpha_1}   } d^2u= (\frac{\gamma^2}{4} +1 -\frac{\gamma \alpha_1}{2} ) ( \frac{\gamma^2}{4}  -\frac{\gamma \alpha_1}{2} )   \frac{\pi}{ l(\frac{\gamma \alpha_1}{2})   l(- \frac{\gamma^2}{4})l (2- \frac{\gamma \alpha_1}{2}+ \frac{\gamma^2}{4})} .
\end{equation}

%
Finally, by taking the $\partial_{z\bar{z}}$ derivative, we get 
\begin{equation}\label{gammazbarz}
(\frac{\gamma^2}{4})^2  \int_{\C}    \frac{|u-1|^{\frac{\gamma^2}{2}}  }{|1-u|^2 |u|^{\gamma \alpha_1}   } d^2u= (\frac{\gamma^2}{4} +1 -\frac{\gamma \alpha_1}{2} )^2  \frac{\pi}{ l(\frac{\gamma \alpha_1}{2})   l(- \frac{\gamma^2}{4})l (2- \frac{\gamma \alpha_1}{2}+ \frac{\gamma^2}{4})} .
\end{equation}



{\small 
\begin{thebibliography}{20}

 
\bibitem{Albeverio1}
Albeverio, S., Hoegh-Krohn, R.: The Wightman axioms and the mass gap for strong interactions
of exponential type in two dimensional space-time. J. Funct. Anal. 16, 39-82 (1974).


\bibitem{Albeverio2}
S. Albeverio, G. Gallavotti and R. Hoegh-Krohn: 
Some Results for the Exponential Interaction in Two or More Dimensions, \emph{Commun. Math. Phys}., 70,187-192 (1979).

\bibitem{AGT}
L. F. Alday, D. Gaiotto, and Y. Tachikawa. Liouville Correlation Functions from Four Dimensional Gauge Theories, \emph{Lett. Math. Phys.}, \textbf{91} 167-197 (2010).

\bibitem{aru}
Aru J., Huang Y., Sun X.:  Two perspectives of the 2D unit area quantum sphere and their equivalence,  \emph{Commun. Math. Phys}., vol 356, Issue 1, pp 261Ð283.
 
 
\bibitem{baverez}
Baverez G., Wong M.D.:  Fusion asymptotics for Liouville correlation functions,  \href{https://arxiv.org/abs/1807.10207}{arXiv:1807.10207}.

\bibitem{BPZ} Belavin A.A., Polyakov A.M., Zamolodchikov A.B. : Infinite conformal symmetry in two-dimensional quantum field theory, \emph{Nuclear Physics B} {\bf 241} (2), 333-380 (1984). 


\bibitem{Ber}
Berestycki N.: An elementary approach to Gaussian multiplicative chaos, \emph{Electronic communications in Probability} {\bf 27}, 1-12 (2017).  





\bibitem{bct} E.Braaten, T.Curtright, C.Thorn. Phys.Lett., {\bf B118} (1982) 115 .

  

%
%

\bibitem{bootstrap}
Collier S., Kravchuk P., , Lin Y-H, , Yin X, Bootstrapping the Spectral Function: On the Uniqueness of Liouville and the Universality of BTZ, \href{https://arxiv.org/abs/1702.00423}{arXiv:1702.00423}.

\bibitem{ct}  
T. Curtright, C.Thorn. Phys.Rev.Lett., {\bf 48} (1982) 1309.

\bibitem{DKRV}
David F., Kupiainen A., Rhodes R., Vargas V.: Liouville Quantum Gravity on the Riemann sphere, \emph{Communications in Mathematical Physics} {\bf 342} (3), 869-907 (2016).

\bibitem{DKRV2}
David F., Kupiainen A., Rhodes R., Vargas V.: Renormalizability of Liouville Quantum Gravity at the Seiberg bound, \href{http://arxiv.org/abs/1506.01968}{arXiv:1506.01968}.


\bibitem{DoOt}
Dorn H., Otto H.-J.: Two and three point functions in Liouville theory, \emph{Nuclear Physics B} {\bf 429} (2), 375-388 (1994).     

 


\bibitem{Dub0}
Dub\'edat J.: SLE and the Free Field: partition functions and couplings, \emph{Journal of the AMS} \textbf{22} (4), 995-1054 (2009). 




\bibitem{DMS}
Duplantier B., Miller J., Sheffield: Liouville quantum gravity as mating of trees,  	\href{http://arxiv.org/abs/1409.7055}{arXiv:1409.7055}.


\bibitem{Ry1}
El-Showk, S., Paulos M.F., Poland D., Rychkov S., Simmons-Duffin  D.,  Vichi A.: Solving the 3D Ising model with the conformal bootstrap, Phys. Rev. D 86, 025022 (2012).

\bibitem{Ry2}
El-Showk, S., Paulos M.F., Poland D., Rychkov S., Simmons-Duffin  D.,  Vichi A.: Solving the 3D Ising model with the conformal bootstrap II. c-minimization and precise critical exponents , Phys. Rev. D 86, 025022 (2012).




\bibitem{fybu}
Fyodorov Y., Bouchaud J.-P.: Freezing and extreme value statistics in a Random
Energy Model with logarithmically correlated potential, \emph{J. Phys.A: Math.Theor} 41 (2008) 372001.
 
 \bibitem{FLeR}
 Fyodorov Y., Le Doussal P., Rosso A.: Statistical Mechanics of Logarithmic REM: Duality, Freezing
and Extreme Value Statistics of $1/f$ Noises generated by Gaussian Free Fields, \emph{J. Stat.Mech.} P10005 (2009).
 
 
\bibitem{gn} Gervais J.-L. and  Neveu A. Nucl.Phys., {\bf B238} (1984) 125.
 
 \bibitem{grimmett}
Grimmett G.R.: Percolation, second edition, Springer-Verlag, Berlin, Grundlehren der Mathematischen Wissenschaften,  321, 1999.

 \bibitem{HaMaWi} Harlow D., Maltz J., Witten E.: Analytic Continuation of Liouville Theory, \emph{Journal of High Energy Physics} (2011).
 
 
 \bibitem{hyp}
 \href{http://functions.wolfram.com/HypergeometricFunctions/Hypergeometric2F1/03/03/04/.}{http://functions.wolfram.com/HypergeometricFunctions/Hypergeometric2F1/03/03/04/}.
 
 \bibitem{IJK}
 Ikhlef Y., Jacobsen J.L., Saleur H.: Three-point functions in $c\leq 1$ Liouville theory and conformal loop ensembles, Phys. Rev. Lett. 116, 130601 (2016).
 
 
\bibitem{cf:Kah} 
Kahane, J.-P.: Sur le chaos multiplicatif,
  \emph{Ann. Sci. Math. Qu{\'e}bec}, \textbf{9} no.2 (1985), 105-150.
 
 \bibitem{KaraSh}
 Karatzas I. Shreve S.: Brownian motion and stochastic calculus, Springer-Verlag.
 
\bibitem{cf:KPZ} Knizhnik, V.G., Polyakov, A.M., Zamolodchikov, A.B.: Fractal structure of 2D-quantum gravity, \emph{Modern Phys. Lett A}, \textbf{3}(8), 819-826 (1988).
 
 
 
 
  \bibitem{cargese} Kupiainen A., Constructive Liouville Conformal Field Theory,  \href{https://arxiv.org/abs/1611.05243}{arXiv:1611.05243}.

 \bibitem{KRV}
 Kupiainen A., Rhodes R., Vargas V.: Local conformal structure of Liouville Quantum Gravity, \emph{Commun. Math. Phys.} (2018) online first, \href{https://arxiv.org/abs/1512.01802}{arXiv:1512.01802}. 
 
\bibitem{KRVjhep}
 Kupiainen A., Rhodes R., Vargas V.: The DOZZ formula from the path integral,  \href{https://doi.org/10.1007/JHEP05(2018)094}{J. High Energ. Phys. (2018) 2018: 94}. 

%

\bibitem{JFL} 
Le Gall J.-F.: Brownian geometry, \href{https://arxiv.org/pdf/1810.02664.pdf}{arXiv:1810.02664}.    

%
%
%


\bibitem{MO}
Maulik D. Okounkov A.: Quantum Groups and Quantum Cohomology, \href{https://arxiv.org/abs/1211.1287}{arXiv:1211.1287}.  


\bibitem{MS1}
Miller J., Sheffield S.: Liouville quantum gravity and the Brownian map I: The QLE(8/3,0) metric, \href{https://arxiv.org/abs/1507.00719}{arXiv:1507.00719 }.


\bibitem{MS2}
Miller J., Sheffield S.: Liouville quantum gravity and the Brownian map II: geodesics and continuity of the embedding, \href{https://arxiv.org/abs/1605.03563}{arXiv:1605.03563}.



\bibitem{MS3}
Miller J., Sheffield S.: Liouville quantum gravity and the Brownian map III: The conformal structure is determined, \href{https://arxiv.org/abs/1608.05391}{arXiv:1608.05391}.


\bibitem{nakayama}
Nakayama Y.: Liouville field theory: a decade after the revolution, \emph{Int.J.Mod.Phys. A} \textbf{19}, 2771-2930 (2004).

\bibitem{OPS}
 O'Raifeartaigh L., Pawlowski J.M., Sreedhar V.V.: The Two-exponential Liouville Theory and the Uniqueness of the Three-point Function, \emph{Physics Letters B} {\bf 481} (2-4), 436-444 (2000).
 
   
\bibitem{OS1}
 Osterwalder, K., Schrader, R.:
Axioms for Euclidean Green's functions I.
\emph{Commun. Math. Phys.} 31, 83-112 (1973).

\bibitem{OS2}
Osterwalder, K., Schrader, R.:
Axioms for Euclidean Green's functions II, \emph{Commun. Math. Phys.} 42, 281-305 (1975).
 
\bibitem{ostrovsky} 
Ostrovsky D.: A Review of Conjectured Laws of Total Mass of Bacry-Muzy GMC Measures on the Interval and Circle and Their Applications. \emph{Reviews in Mathematical Physics}, vol 30, Issue 10 (2018).

 \bibitem{ostrovsky2} 
Ostrovsky D.: Mellin transform of the limit lognormal distribution, \emph{Commun. Math. Phys}., vol 288, Issue 1, pp 287Ð310(2009).

 
\bibitem{Pol}
Polyakov A.M.: Quantum geometry of bosonic strings, \emph{Phys. Lett. } \textbf{103B} 207 (1981).

\bibitem{Pol1}
Polyakov A.M.: From Quarks to Strings, \href{https://arxiv.org/abs/0812.0183}{arXiv:0812.0183}. 

\bibitem{pol74}Polyakov A. M.: Nonhamiltonian approach to conformal quantum field theory, Zh.
Eksp. Teor. Fiz. 66 (1974) 23-42.
 
\bibitem{remy}
Remy G.: The Fyodorov-Bouchaud formula and Liouville Conformal Field theory, \href{https://arxiv.org/abs/1710.06897}{arXiv:1710.06897 }.

\bibitem{zhu}
Remy G., Zhu T.: The distribution of Gaussian multiplicative chaos on the unit interval, \href{https://arxiv.org/abs/1804.02942}{arXiv:1804.02942}.

 
\bibitem{Rib}
Ribault S.: Conformal Field theory on the plane, \href{https://arxiv.org/abs/1406.4290}{arXiv:1406.4290}. 


\bibitem{review} 
Rhodes R., Vargas, V.: Gaussian multiplicative chaos and applications: a review,  \href{https://projecteuclid.org/euclid.ps/1415023603}{ \emph{Probability Surveys}, Volume 11 (2014), 315-392}.

\bibitem{RV}
Rhodes R., Vargas V.: Lecture notes on Gaussian multiplicative chaos and Liouville Quantum Gravity, \href{https://arxiv.org/abs/1602.07323}{arXiv:1602.07323}.  

\bibitem{MRM}
Rhodes R., Vargas V.: Multidimensional Multifractal Random Measures,   \emph{Electronic Journal of Probability}, 2010, Vol 15, 241-258.

\bibitem{tail}
Rhodes R., Vargas V.: The Tail expansion of Gaussian multiplicative chaos and the Liouville reflection coefficient, to appear in \emph{Annals of Probability}, \href{https://arxiv.org/abs/1710.02096}{arXiv:1710.02096}.



\bibitem{SV} 
 Schiffmann O., Vasserot E.: Cherednik algebras, W-algebras and the equivariant cohomology of the moduli space of instantons on $A^2$, \emph{Publications math\'ematiques de l'IHES} {\bf 118} (1), 213-342 (2013).

\bibitem{seiberg}
Seiberg N.: Notes on Quantum Liouville Theory and Quantum  Gravity, \emph{Progress of Theoretical Physics}, suppl. 102, 1990.

\bibitem{She07}
Sheffield S.: Gaussian free fields for mathematicians, \emph{Probab. Th. Rel. Fields} \textbf{139}, 
521-541 (2007).

 






 
\bibitem{Tesc} Teschner J.: On the Liouville three point function, \emph{Phys. Lett. B} {\bf 363}, 65-70 (1995).

\bibitem{Tesc1} 
Teschner J.: Liouville Theory Revisited, Class.Quant.Grav.\emph{18}:R153-R222,(2001).

\bibitem{Tesc2}  Teschner J.: A lecture on the Liouville vertex operators, Int.J.Mod.Phys. A19S2, 436-458 (2004)

\bibitem{varg}
Vargas V.: Lecture notes on Liouville theory and the DOZZ formula, \href{https://arxiv.org/abs/1712.00829}{arXiv:1712.00829}.

\bibitem{Williams}
Williams, D.: Path Decomposition and Continuity of Local Time for One-Dimensional Diffusions, I, \emph{Proceedings of the London Mathematical Society} {\bf s3-28} (4), 738-768 (1974).

%



\bibitem{ZaZaarxiv}
Zamolodchikov A.B., Zamolodchikov A.B.:  Structure constants and conformal bootstrap
in Liouville field theory, arXiv:hep-th/9506136.

\bibitem{ZaZa}
Zamolodchikov A.B., Zamolodchikov A.B.: Conformal bootstrap
in Liouville field theory, \emph{ Nuclear Physics B}
{\bf 477} (2), 577-605 (1996).


\end{thebibliography}}
\end{document}